%% file: main.tex
\documentclass[11pt]{article}
\usepackage{amsmath,amssymb,latexsym,float,epsfig}
\usepackage{latexsym}
\usepackage{amsthm}
\usepackage{cases}
\usepackage{epsfig}
\usepackage{algpseudocode}
\usepackage{csquotes}

\usepackage{scalefnt}
\usepackage{booktabs}
\usepackage[english]{babel}
\usepackage{graphicx}
\usepackage{titlesec}
\usepackage{epsfig}
\usepackage{enumerate}
\usepackage{caption}
\usepackage{subcaption}
\usepackage{booktabs}
\usepackage[english]{babel}
\usepackage{graphicx}
\usepackage{multirow}
\usepackage{rotating}
\usepackage{fancyhdr}
\usepackage{setspace}
\usepackage[section]{placeins}
\usepackage{xcolor}

\usepackage{breqn}
\usepackage{longtable}
\usepackage{makecell}
\usepackage{array, makecell}
\usepackage{pseudocode}
\usepackage{algorithmicx}
\usepackage{algorithm}
\usepackage{url} 
\usepackage[title]{appendix}
\usepackage{ragged2e}
\usepackage{pgfplots}
\usepackage{pgfkeys}
\usepackage{tikz}
\usepackage{graphics}
\usetikzlibrary{patterns}
\usetikzlibrary{patterns.meta}
\usepackage{color}
\usepackage[
backend=biber,
style=numeric,
sorting=nyt
]{biblatex}
\newtheorem{theorem}{Theorem}
\newtheorem{corollary}[theorem]{Corollary}
\newtheorem{lemma}[theorem]{Lemma}

\usepackage[nocomma]{optidef}
\usepackage{mathtools,zref-savepos}
\newtagform{roman}[]()

\topmargin by -\headsep \textheight 8.9in \oddsidemargin 0pt
\evensidemargin \oddsidemargin \marginparwidth 0.5in \textwidth
6.5in
\begin{document}
\title{Managing Word-of-Mouth through Capacity Allocation and Advertisement}

\author{ Asya Dipkaya \thanks{asyadipkaya@gmail.com}\and Sakine Batun \thanks{sakine@metu.edu.tr} \and Bahar \c{C}avdar \thanks{bcavdar@tamu.edu}}
\date{
    $^*$$^\dagger$ \small{Department of Industrial Engineering, Middle East Technical University, Ankara, Turkey }\\%
        $^\ddagger$\small{Department of Engineering Technology and Industrial Distribution, Texas A\&M University, College Station, TX, USA }\\%

     %\today
}

%\date{}
\maketitle
\begin{abstract}
{
Advancements in service sector and growing online platforms are intensifying the information exchange between customers through (electronic) word-of-mouth (WoM). The information obtained by WoM has shown to be a dominant factor in customers' purchase decisions creating an endogenous demand structure. Service providers can monitor how their service is perceived by consumers through different methods, for example surveys. Although this requires an additional effort, the understanding and integration of endogenous demand into operational decisions offer great benefits. In this paper, we study a service system where customers are sensitive to the on-demand access to the service. Customers form a perception based on the information obtained by WoM communication and the advertisement. Depending on the type of the service environment, service capacity can be flexible or constant. We consider two types of service providers: aware firm that has complete information on endogenous demand structure, and naïve firm that has partial information. Our focus is to understand the optimal advertisement and capacity decisions, and the value of information on the underlying demand. For firms that have flexible service capacity, we show that it is optimal to employ aggressive advertisement strategies in the early stages. Myopic naive firms often misinterpret the market conditions and cease operating, where they could in fact realize profit. For the cases where the capacity is not flexible, it may not be possible to avoid negative WoM. Therefore, the service provider is forced to keep the level of advertisements more compatible with the actual quality of the service. This prevents the firm from overcrowding the system through advertisement without considering the service quality. %Therefore, WoM acts as a corrective tool for the firm.
}
\end{abstract}


\input{sections/introduction}

\input{sections/literature}

\input{sections/modelling}

\input{sections/structural_results}

\input{sections/computational_results}

\input{sections/conclusion}

\printbibliography

\newpage

\begin{appendices}
\input{sections/appendices}

\end{appendices}

\end{document}

%% file: sections/introduction.tex
\section{Introduction}
\label{intro}
Easy access to information has created many opportunities for consumers and service providers. New channels provide new communication tools for customers enabling word-of-mouth (WoM) communication, and affect consumer choices. Companies also have means to understand where they stand in the eyes of consumers through the information shared on online platforms. Studies have shown that WoM communication is highly effective in shaping consumers' opinions since the transmitted information is believed to be without commercial concern and hence is relied on (\cite{litvin2008electronic}, \cite{huete2017literature}). Online reviews, ratings and social media constitute a less traditional type of WoM, namely eWoM (electronic Word-of-Mouth), by which the information sharer broadcasts messages based on their experience with the product or service of concern, rather than having the messages transmitted from person to person. This amplifies the impact of WoM and the spread of information \cite{huete2017literature}. As a result, WoM plays an important role for consumers by helping them have a better understanding of the firms' true performances. Studies regarding consumer behavior and purchase decisions suggest that each consumer has some expectation from a product/service that gets altered as the service level changes (\cite{shapiro1982consumer}, \cite{bolton1991longitudinal}). Adopting this observation, we consider consumers who have prior beliefs or first impressions regarding a service provider's performance that get updated as consumers receive new information when WoM communication takes effect.

Having underlined the significance of the communication between consumers, it is not realistic to consider firms as passive agents. Firms often place attempts in the form of rebates, promotions or advertisement to influence consumers' purchasing decisions \cite{zeithaml1988consumer}. The intensity of such efforts determines how much, or through how many different channels, consumers are exposed to the new information provided by the firm. Therefore, we consider such activities as a factor affecting consumer perception, besides WoM. In the figure below, we represent the mechanism that determines the consumer perception.
\input{wom.tex}
The information on service quality is especially important to customers in e-commerce industry, where there are a large number of alternatives for customers to pick from. In such services (e.g., online grocery shopping services, food retailers, transportation services, call centers), the availability of the service is of essence. There is an extensive literature underlining the impact of fill rate on demand and customer retention (e.g., \cite{fitzsimons2000consumer}, \cite{heim2001operational}, \cite{garnett2002designing}, \cite{urban2005inventory}, \cite{mckinnon2007store}, \cite{adelman2013dynamic}, \cite{craig2016impact}, \cite{mandal2018stocking}). In such service environments, the likelihood of being able to receive the service when requested is an important determinant of demand, and customers associate this likelihood with the fill rate of the firm. For example, consider restaurants that offer food delivery services through delivery platforms like Uber Eats, Postmates, Grubhub, etc. When the restaurants receive a large number of orders that exceed their capacity, they do not accept new orders losing potential customers. As this becomes more frequent, it is expected that due to low fill rate, customers will divert to other options. To recover from this, firms can use tools such as advertising promising better service (e.g., higher availability rates). Managing the service capacity becomes vital under such settings. There are a number of studies in the literature suggesting different methods to respond to volatile demand by adjusting service capacity (\cite{aksin2008effective}, \cite{fuhrman1999co}, \cite{dorsch2014combining}). In many settings, those methods are applied and benefited from; however, depending on the nature of the service and the type of capacity units, flexible capacity may not always be applicable.

In this paper, we consider a service environment with two main actors: a firm that provides one type of service and a customer population that exchanges information regarding the fill rate of the firm via WoM (both traditional and eWoM). Each consumer has a perception regarding the fill rate of the firm and a utility function that takes this perception into account. Consequently, we consider a utility-based demand function. Meanwhile, firm has two actions it can use to manage demand by interfering in consumer perception:
\begin{itemize}
    \item [(i)] placing effort in advertising to modify perception by directly communicating with customers,
    \item [(ii)] adjusting the service capacity level and modify perception indirectly by encouraging positive WoM as a result of high fill rate.
\end{itemize}

We study the advertising and capacity decisions considering endogenous demand on a finite horizon. We focus on strategies to balance the traditional operational decisions with the management of customer perception through WoM. The research questions we study are as follows:
\begin{itemize}
	\item What kind of advertisement strategy should a firm adopt? Under which conditions should it advertise more aggressively?
	\item What are the characteristics of a capacity decision under endogenous demand? Is it optimal to always fulfill demand?
	\item How does the advertising policy differ between constant- and flexible-capacity systems?
	%\item Is it possible to exploit WoM? Does WoM motivate the firm to provide better service quality?
	\item What is the value of information on the endogenous demand structure?
\end{itemize}

We model the consumer utility and endogenous demand structure in detail, and study the firm's optimal policies under variable and constant service capacities. Under each of these scenarios, we compare two types of service providers:
\begin{itemize}
    \item [(i)] \emph{aware firm}, which has complete information on the underlying endogenous demand structure.
    \item [(ii)] \emph{naive firm}, which has partial information on the underlying endogenous demand structure.
\end{itemize}

We show that for firms with flexible service capacity, it is always optimal to maintain enough capacity to fulfill demand under certain cost functions. Consequently, such firms never suffer negative WoM. Therefore, they are able to seize greater portions of the market by employing advertisement strategies that start off aggressively and continue more moderately. We observe high investment on advertising and fast reputation growth in the early periods in optimal policies. On the other hand, the firms that cannot adjust the service capacity often face both positive and negative WoM. The advertisement policies become more conservative to avoid negative WoM. Therefore, the existence of the WoM communication forces firms to avoid making promises that they cannot fulfill.
%The advertisement strategy of such firms cannot be generalized, as its structure highly depends on the market parameters. 
%These firms with constant service capacity often need to come up with complex advertising policies to manage demand. Additionally, we discuss that having complete information on the endogenous demand is an important tool for the firms to make use of, especially when the service capacity is constant. In the absence of this information, firms cannot make long term plans as they cannot predict how their actions will impact their position in the market. In that case, firms make myopic decisions that result in inefficient policies. They often misjudge how challenging the market is and decide to stay out, missing out on unforeseen opportunities.

The rest of the paper is organized as follows: in Section \ref{sec:1}, we review the related literature. In Section \ref{sec:2}, we model the problem for aware and naive firms under variable and constant service capacity scenarios. In Sections \ref{sec:3} and \ref{sec:4}, we provide the structural results on the optimal policies and computational results on the value of information for the endogenous demand (or WoM communication). Finally, we provide concluding remarks in Section \ref{sec:5}.

%% file: wom.tex
\tikzset{every picture/.style={line width=0.75pt}} %set default line width to 0.75pt        
\begin{figure}[H]
\vspace{-3mm}
    \centering
        \begin{tikzpicture}[x=0.75pt,y=0.75pt,yscale=-1,xscale=1]
        
        %Straight Lines [id:da7142333325902221] 
        \draw    (93,85) -- (93,115) ;
        \draw [shift={(93,114.33)}, rotate = 270] [color={rgb, 255:red, 0; green, 0; blue, 0 }  ][line width=0.75]    (10.93,-3.29) .. controls (6.95,-1.4) and (3.31,-0.3) .. (0,0) .. controls (3.31,0.3) and (6.95,1.4) .. (10.93,3.29)   ;
        %Straight Lines [id:da47904045606467593] 
        \draw    (93,185) -- (93,155) ;
        \draw [shift={(93,155)}, rotate = 450] [color={rgb, 255:red, 0; green, 0; blue, 0 }  ][line width=0.75]    (10.93,-3.29) .. controls (6.95,-1.4) and (3.31,-0.3) .. (0,0) .. controls (3.31,0.3) and (6.95,1.4) .. (10.93,3.29)   ;
        %Straight Lines [id:da3917586250514016] 
        \draw    (135,137) -- (190,137) ;
        \draw [shift={(190,137)}, rotate = 180] [color={rgb, 255:red, 0; green, 0; blue, 0 }  ][line width=0.75]    (10.93,-3.29) .. controls (6.95,-1.4) and (3.31,-0.3) .. (0,0) .. controls (3.31,0.3) and (6.95,1.4) .. (10.93,3.29)   ;
        
        % Text Node
        \draw (70,120) node [anchor=north west][inner sep=0.75pt]   [align=left] {customer\\opinions};
        % Text Node
        \draw (35,70) node [anchor=north west][inner sep=0.75pt]   [align=left] {WoM (about firm's performance)};
        % Text Node
        \draw (55,190) node [anchor=north west][inner sep=0.75pt]   [align=left] {advertisements};
        % Text Node
        \draw (205,119) node [anchor=north west][inner sep=0.75pt]   [align=left] {updated customer\\opinions};

        \end{tikzpicture}
    \caption{Consumer Perception, WoM and Advertising}
    \label{fig:simple_wom}
\end{figure}
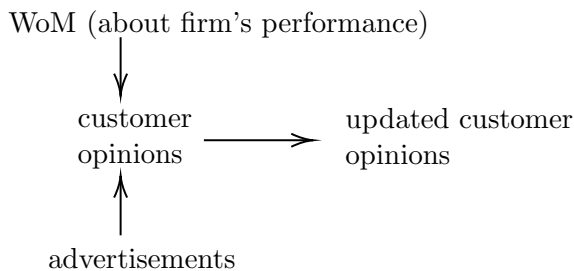

%% file: sections/literature.tex
\section{Literature Review}
\label{sec:1}
In this section, we review the related literature on WoM, consumer behavior, characteristics of the demand function and operational decision making under endogenous demand.

There are different definitions for WoM \cite{huete2017literature}. Some studies define it as uni-directional information transmission between two individuals when the receiver believes that the sender has no commercial concern \cite{daugherty2014ewom}, while some define it as bi-directional communication between customers, where the senders are without commercial concern \cite{litvin2008electronic}. In our study, we consider WoM communication occurring between multiple agents that have no commercial concern.

When it comes to eWoM, it is discussed that such online mechanisms function as a trust-builder and motivate collaboration among strangers that share information about companies, products or services \cite{dellarocas2003digitization}. Comparing eWoM with traditional WoM, eWoM's great accessibility and the faster diffusion of information stand out \cite{huete2017literature}. With motivated consumers to share their experiences, the amount of available information grows rapidly in online marketplaces. Moreover, customers tend to trust the information they gather from such platforms more than what the companies promise through promotion activities \cite{nieto2014marketing}. It is clear that eWoM is largely affecting consumer decisions, and therefore demand. As a result, companies may view this as a downside, since their influence on demand is lower than it used to be \cite{yang2017effects}. There is an increasing number of studies emphasizing that incorporating the information exchange between customers and the impacts of it in operational decision making is important \cite{mayzlin2006promotional}. \cite{cavdarWoM} is one of the first studies that integrates the WoM into an operational problem. They study a consolidated shipment policy problem considering the information exchange between different classes of customers, and present results on how to manage WoM communication through operational decisions.

In order to understand and model the consumer behavior better, we review the literature regarding consumer behavior and the relationship between consumers' experiences, expectations, beliefs and demand. Later, we combine these findings with WoM behavior to construct our model. To begin with, Shapiro (1982) suggests that each customer has some expectation from the service, which determines whether she will decide to purchase, implying that demand for the service is based on customers' expectations \cite{shapiro1982consumer}. Additionally, it is discussed that customers' expectations depend on several factors, including brand reputation, price and advertising (\cite{zeithaml1988consumer}, \cite{johnson1995rational}, \cite{dodds1991effects}, \cite{boulding1993dynamic}). Studies show that changes in the service level have a significant effect on how customers rate service quality while customers’ behavior is affected by their prior beliefs as well \cite{bolton1991longitudinal}. These results lead us to defining a measure to reflect customers' current perception of the company as well as their expectations from the service. On top of that, we consider that consumers exchange information regarding service quality by WoM, get exposed to advertisement, and update their perception accordingly. Then, we model demand based on this perception measure.

In Adelman's (2013) study on how to allocate a limited capacity among a set of customers, it is argued that customers remember company's past performance in fulfilling demand, and their purchase decisions are affected by it \cite{adelman2013dynamic}. In our study, we consider a service environment, where perception towards service quality is largely determined by company's performance in fulfilling demand. Therefore, we accommodate Adelman's observation in our demand model, with a consumer population that share information regarding the firm's fill rate performance. This fill rate-sensitive service environment underlines the importance of service capacity in fulfilling demand (e.g., number of couriers of an online grocery shopping firm). Therefore, we incorporate firm's operational decisions that affect the fill rate while modelling the endogenous demand. An example of operational decision making under endogenous demand is provided in a study carried out by Mandal et al. (2018), where the classical newsvendor problem is extended by introducing behavioral biases \cite{mandal2018stocking}. They also consider endogenous demand which depends on how probable the consumers think getting served is. In their problem, this probability depends on the inventory level, which corresponds to service capacity level in our case. They show that this probability is equal to firm's fill rate. They highlight the importance of decision making by accounting for endogenous demand, and show that such demand models motivate firms to keep the fill rate high. A similar finding is presented in the study carried out by Urban (2005) regarding inventory-dependent demand models \cite{urban2005inventory}. Urban classifies such demand models in two groups. The first group covers the models where demand is affected by the initial inventory while the second covers demand models that incorporate instantaneous inventory levels. Our model falls under the second category.

There are extensive empirical studies on the relationship between fill rate and demand, which report similar results. Craig et al. (2016) examine the relationship between service level and demand for a well-known functional apparel product supplier and show that inventory level of the supplier has a statistically significant impact on demand \cite{craig2016impact}. In parallel, Heim and Sinha (2001) establish that service availability and timeliness are strongly associated with customer retention in their study related with electronic food retailers \cite{heim2001operational}. Similar results on customers looking for substitutes or not making the purchase at all in case of stockouts are shown for different product groups as well (\cite{fitzsimons2000consumer}, \cite{mckinnon2007store}).

When considering a service environment, where fill rate is essential to maintain demand, new questions arise regarding how service capacity is managed. One common method is \emph{shared services} (insourcing), which are organizations established by companies to carry out similar business functions for different centers/units of the same company \cite{aksin2008effective}. Another method is outsourcing the service capacity entirely while Fuhrman (1999) suggests \emph{co-sourcing} as a better alternative, which is the practice of meeting some of the demand using internal resources and outsourcing the rest \cite{fuhrman1999co}. These methods, if applicable, allow firms to operate with different service capacity levels throughout a planning horizon \cite{akcsin2008call}. Note that the studies we have mentioned so far regarding the service capacity management are mostly devoted to call center optimization, which is a befitting service environment example where purchase decisions of customers depend on whether they are able to get served in a timely manner, or not \cite{garnett2002designing}. Meanwhile, Dorsch and Häckel (2014) provide a study on managing service capacity by different outsourcing options \cite{dorsch2014combining}. They consider a service firm outsourcing capacity by different contract types: (i) outsourcing a fixed capacity from a vendor, (ii) outsourcing elastic capacity and (iii) outsourcing surplus capacity. The first two options rely on long-term contracts while the third can be executed in short notice. Outsourcing surplus capacity is often observed via \emph{Business Process as a Service (BPaaS)}, in which business processes are outsourced through cloud computing services on demand. When applicable, this practice  can provide firms with great agility to meet volatile demand.

Next, we review different types of demand functions in a service environment. Huang et al. (2013) present a comprehensive review of different demand function types used in the literature \cite{huang2013demand}. They discuss that the factors that influence consumers' purchasing decisions are price, rebate, lead-time, space (area occupied on shelves), quality and advertisement. To represent the impact of the perceived quality's on consumers' purchasing decisions and capture consumer behavior better; we use a demand model based on consumer utility as a function of service quality. In the preferred demand model, Moorthy (1988) and Tirole (1988) suggest that the utility obtained by a single consumer when the service is purchased can be modeled as follows (\cite{moorthy1988product}, \cite{tirole1988theory}).
\begin{equation}
 u(p) = \theta y - p,
\end{equation}
where $y$ denotes product/service quality, $p$ denotes price, $u$ denotes utility as a function of price and $\theta$ is how much utility the consumer gains for unit quality if she receives the product/service. $\theta$ is used to define different customer types. A customer with low $\theta$ has lower quality valuation; therefore, gains less utility from unit quality. Consumers with positive net surplus make the purchase and others do not. When $\theta $ is uniformly distributed between 0 and 1, and market size is $a$, demand function is derived as follows.
\begin{equation}
d(p,\:a) = a \left (1 - F \left( \frac{p}{y} \right) \right),    
\end{equation}
where $F$ is the c.d.f. of $\theta$.

We extend this model by accounting for advertisement, WoM and customers' reaction to company's perceived reputation, which is explained in detail in the upcoming sections. We model a service environment that accounts for the above described factors by an endogenous demand structure. In a nutshell, we consider a set of consumers that have a first impression about the fill rate of the service, receive new information from different sources, update their opinions, generate demand based on utility and observe the new fill rate. In the meanwhile, firm can make capacity and advertisement decisions.

In the current literature, operational decision making under endogenous demand is mainly observed for product suppliers, rather than service providers. In our model, we capture consumer behavior in a service environment by a rigorous demand mechanism, using multiple time periods to observe the turn of events and long-term policies, which is our main contribution to the literature.

%% file: sections/modelling.tex
\section{Modeling}
\label{sec:2}
We consider a system where a company provides one type of service. Customers are faced with the choice of requesting the service from the company or looking for substitutes. When they are making this decision, their belief about the company's performance, or fill rate, plays an important role. If the company has been successful in fulfilling demand in the past, customers find it more appealing. The information on whether their service request will be fulfilled by the company is not immediately available to customers. However, they can infer the fill rate based on the available information. One source of information is WoM communication. As customers use the service provided by the company, they share their experiences through communication channels. In addition to this, company can also take an action to alter the perception of the fill rate. These actions can include advertising, promotions, etc. Therefore, customers are exposed to two sources of information: (i) historical fill rate and (ii) promises made by company regarding the future service. The former is available through WoM, the latter depends on the advertisement/promotion efforts of the company. Observing these two sources of information, customers develop an idea about the fill rate in case they request the service. Therefore, the demand is endogenous to operational decisions. Our focus is on understanding the impact of this demand structure and WoM communication on advertisement and capacity decisions in a multi-period problem.

We express consumers' collective perception towards the fill rate of the firm by the \textbf{\emph{reputation}} of the firm. In the multi-period problem, reputation is updated every period, which is then reflected on demand.

\input{endogenous_demand}
Figure \ref{fig:endo_demand} expresses the endogenous demand structure in two consecutive time periods. At the beginning of each time period $i$, the firm observes where it stands in the eyes of the customers, that is the reputation, denoted by $R_i$. Since the reputation of the company is based on the fill rate, it takes a value in $[0,1]$. Observing the reputation, the firm makes decisions on:
\begin{itemize}
    \item[(i)] advertisement, $A_i\in [0,1]$,
    \item[(ii)] service capacity, $S_i$.
\end{itemize}
We assume that these decisions take effect instantaneously at the beginning of period $i$, and the advertisement efforts improve consumer perception, creating the \emph{post-advertisement reputation} $R^{'}_i\in [0,1]$. Consequently, demand $D_i$ is realized, and the firm operates to fulfill it as much as the capacity, $S_i$ allows. For each unit of service provided, company receives revenue, which is equal to the service price, $p$. We assume that the cost of maintaining one unit of service capacity is $c$ to the company, and $p$ is strictly greater than $c$. Additionally, the firm incurs the cost of $c_A$ for the unit advertisement effort. At the end of the period, the fill rate is realized, which then impacts the next period's reputation $R_{i+1}$ through WoM. Notice that advertisement's impact on demand is direct and immediate while the choice of capacity impacts the demand indirectly through WoM, and not until the next period.

In the following part, we present the formulations for WoM communication, demand as a function of operational decisions and WoM, and the profit function.

First, we model the post-advertisement reputation in period $i$ as follows:
\begin{equation}
    \label{Eq:PostAd}
R_i^{'}=R_i+(1-R_i)A_i^{\alpha},
\end{equation}
where the parameter $\alpha > 0$ is a reflection of customers' reaction to the observed advertisement. Higher values of $\alpha$ imply that customers are more resistant to advertisement and more prone to basing their decisions on the firm's performance. Therefore, we refer to $\alpha$ as the \emph{advertisement resistance} of the customers. Due to the multiplicative component in (\ref{Eq:PostAd}), advertisement's effect diminishes at higher reputation levels.

Advertisement decision can be handled in different ways. In this study, we consider that company has two options: high level and low level advertisements, denoted by $H$ and $L$, respectively. Over an $N$-period planning horizon, the firm has to decide on an advertisement sequence such as $A_1 = H, \: A_2= L, \: A_3 = L, ...\: A_N = L$.

When customers are faced with multiple options, customer choice is commonly modeled by utility-based demand functions \cite{huang2013demand}. To model customer demand, we extend the utility model in \cite{moorthy1988product} and \cite{tirole1988theory} by integrating the company's reputation into the customer utility.

We consider a customer population of size $\Lambda$, where each customer has a \emph{taste} or \emph{reputation valuation} for the service provided by the company denoted by $\theta \sim U(0,1)$. If customers do not prefer to request the service from the company, they receive the reservation utility $u$, which is the utility they receive from the substitutes of the service. A customer with a taste $\theta$ for the service, receives the following utility in period $i$ if she makes the purchase.
\begin{equation}
    m\: \theta \: R_i^{'} - p,
    \label{Eq:utility}
\end{equation}
where $m$ is the monetary value of reputation. By purchasing the service, the customer receives the utility of $m\: \theta \: R_i^{'}$ and sacrifices $p$. We assume that utility and price are measured in the same units without loss of generality. If the expression in (\ref{Eq:utility}) is greater than the reservation utility $u$, it is worth for the customer to purchase the service. So, the customers, whose taste parameter $\theta$ satisfies the following condition, create demand:
\begin{equation}
\label{Eq:eta}
    \theta > \frac{u+p}{m R_i^{'}}
\end{equation}
Then, the fraction of customers, who create demand at period $i$ is $1 - F(\frac{u+p}{m R_i^{'}})$, where $F$ is the cumulative distribution function of $\theta$. Since $\theta$ is uniformly distributed between 0 and 1, this fraction can be written as $ 1 - \min \{1, \: \frac{u+p}{m R'_i}\}$. Accordingly, the realized demand in period $i$ can be expressed as:
\begin{equation}
    D_i = \max\left\{0, \: \Lambda \left(1 - \frac{u+p}{m R_i^{'}}\right) \right\}
\end{equation}
We assume $\frac{u+p}{m} < 1$ to avoid trivial cases. Otherwise, regardless of the company's reputation, no demand can be generated.
 
Finally, the impact of fill rate on the reputation is modeled by a linear function as follows:
\begin{equation}
\label{Eq:RepWoM}
    R_{i+1}=
    \begin{cases}
    (1-w) R_i^{'} + w F_i, & \text{if  } D_i>0\\
    R_i^{'}, & \text{otherwise}
    \end{cases}
\end{equation}
where $w$ is the weight assigned to the information on the most recent fill rate gained by WoM and $w \in [0,\:1]$. Reputation is altered based on this function at the end of each time period. $w=1$ implies that customers only use the most recent information to determine the reputation of the company, while $w=0$ implies that customers are not affected by the WoM. When there is no demand generated in the current period, the reputation of the company remains the same since there is no new information on the fill rate.

Note that reputation level increases with advertisement while the impact of the fill rate can be in both directions. Fill rate increases the reputation if $F_i > R_i^{'}$, and decreases it if $F_i < R_i^{'}$. This can be interpreted as consumers expecting firm's performance to be aligned with its reputation. If the firm performs better than expected, reputation goes up; otherwise, reputation goes down.

In period $i$, the company collects the profit $\pi_i$:
\begin{equation}
    \label{Eq:Profit}
    \pi_i=p \: \min\{S_i, D_i\}-c \: S_i-c_{A} \: A_{i}.
\end{equation}
The company's objective is to maximize the total profit $\Pi$ over an $N$-period planning horizon, where $\Pi=\sum_{i=1}^{N}\pi_i$. Capacity selection and advertisement efforts are the main decisions made by the company in each period.

In the rest of this section, we present the problem formulations for different settings. We will first address the case where the company can adjust the service capacity in each time period. In the second setting, we will consider that the service capacity is predetermined and constant over all periods. The former is associated with the cases where increasing/decreasing the capacity is possible by short-term contracts or other capacity-adjustment methods as discussed in Section \ref{sec:1}. In both settings, we analyze the problem from the perspectives of the aware firm and the naive firm. The aware firm has complete information on the endogenous demand structure and can incorporate the impact of its operational decisions on demand. On the other hand, the naive firm is not aware of the fact that the operational decisions impact the demand through WoM. We present our analysis of all models in Section \ref{sec:3}.

\subsection{Variable Capacity}
\label{subsec:varCap}
In this section, we present the problem formulation under variable capacity. The firm can adjust the capacity in each period to avoid lost sales or capacity surplus, depending on its strategy. In addition to the capacity decision, the firm also chooses between the high ($H$) and low ($L$) advertisement levels in each period. We first model the aware firm's problem.

%\subsubsection{Aware Firm Under Variable Capacity}
\textbf{Aware Firm Under Variable Capacity: }
Under variable capacity, the aware firm tries to maximize the total profit throughout a planning horizon of $N$ periods by determining advertisement level $A_i$ and service capacity $S_i$ for each period $i$. Having full information on WoM communication, in any period $i$, the firm is able to make decisions for all future periods. The problem can be formulated as a finite-horizon dynamic program, where the stages correspond to time periods, the state is expressed by the current reputation $R_i$, the actions are advertisement level $A_i$ and service capacity $S_i$, and the value function $P_i(R_i)$ represents the cumulative profit throughout periods $i,i+1,...,N$. The optimality equation can be expressed as:
\begin{equation}
  P_i(R_i) =
    \begin{dcases}
      \max_{\substack{A_i\in \{H,\:L\}, \: S_i\geq 0}} \{\pi_i + P_{i+1}(R_{i+1})\} & \text{if } i<N,\\
      \max_{\substack{A_i\in \{H,\:L\}, \: S_i\geq 0}} \{\pi_i\} & \text{if } i=N,
    \end{dcases}       
\end{equation}
\noindent where
\begin{align*}
        & && \pi_i = p \: \min \{D_i,\:S_i\} - c S_i - c_A A_i,\\
        & && R_{i+1}=
        \begin{dcases}
        (1-w) \left(R_{i} + (1-R_{i}) A_{i}^\alpha \right) + w \frac{\min\{D_i,\:S_i\}}{D_i} & \text{if} \:\: D_i>0,\\
        R_{i} + (1-R_{i}) A_{i}^\alpha & \text{if} \:\: D_i=0,
        \end{dcases}\\
        & && D_i = \max\left\{0,\: \Lambda\left(1 - \frac{u+p}{m \left(R_i + (1-R_i) A_i^\alpha \right)} \right)\right\}.
\end{align*}

The objective is to maximize the total profit over the planning horizon, $\Pi=P_1(R_1)$. The dynamic program above can be solved by backward recursion. However, the structure of the demand function and the continuous state space for the reputation result in a highly complicated formulation to derive a closed form solution. Therefore, we direct our efforts into developing structural results to solve the problem more efficiently and understand the behavior of the optimal solution. The results are presented in Section \ref{sec:3}.

\textbf{Naive Firm Under Variable Capacity: }
Under variable capacity, the naive firm can observe its reputation at the beginning of each period and expect advertisement's effect on demand. However, the firm does not know how customer perception evolves according to the underlying WoM communication, and assumes that there is some exogeneity in demand. Therefore, the naive firm cannot capture the inter-dependency between multiple periods, and follows a myopic approach. This implies single-period profit maximization, and the firm's problem in any period $i$ can be expressed as:
\begin{equation}
    \max_{\substack{A_i\in \{H,\:L\}, \: S_i\geq 0}} \{\pi_i\}
\end{equation}
\noindent where
\begin{align*}
        & && \pi_i = p \: \min \{D_i,\:S_i\} - c S_i - c_A A_i,\\
        & && D_i = \max\left\{0,\: \Lambda\left(1 - \frac{u+p}{m \left(R_i + (1-R_i) A_i^\alpha \right)} \right)\right\}.
\end{align*}

Notice that while the aware firm observes the starting reputation ($R_1$) and plans the optimal advertisement and capacity policies of all periods accounting for WoM, the naive firm observes $R_i$ at the beginning of each period $i$, and makes the optimal advertisement and capacity decisions of each period individually.

\subsection{Constant Capacity}
\label{subsec: constCap}

In this section, we model the system under a fixed service capacity that is predetermined as $S$. In other words, we consider the service capacity as a problem parameter rather than a decision variable. We assume $S <= \Lambda$, as it is not meaningful for the firm to invest in a capacity level greater than the market size. In this new setting, the firm only makes the advertisement decisions, choosing between the high ($H$) and low ($L$) advertisement levels for each period.

\textbf{Aware Firm Under Constant Capacity: }
The aware firm tries to maximize profit considering the endogenous demand through the advertising decision. The problem can be modeled as a dynamic program, where stages correspond to time periods, the state corresponds to the reputation $R_i$, and the actions are $A_i$. The formulation is the same as in the variable capacity scenario except the exclusion of service capacity from the action space. The dynamic program is given in Appendix \ref{model:const_cap}.

\textbf{Naive Firm Under Constant Capacity: }
The naive firm tries to maximize profit by choosing the advertisement level in each period $i$ based on the observed reputation $R_i$. Having fewer actions than the variable capacity problem and less information than the aware firm, this is the most restrictive scenario we examine. The problem the company solves in each period is given in Appendix \ref{model:const_cap}.

%% file: endogenous_demand.tex
\begin{figure}[H]
    \centering
    \resizebox{.8\hsize}{!}{
    \tikzset{every picture/.style={line width=0.75pt}} %set default line width to 0.75pt

    \begin{tikzpicture}[x=0.75pt,y=0.75pt,yscale=-1,xscale=1]
    
    % % Text Node
    \draw (-35,92) node [anchor=north west][inner sep=0.75pt]   [align=left] {$\displaystyle R_{i}$};
    %Straight Lines [id:da9255911732654596] 
    \draw (-10,101) -- (100,101) ;
    \draw [shift={(100,101)}, rotate = 180] [line width=0.75]    (10.93,-3.29) .. controls (6.95,-1.4) and (3.31,-0.3) .. (0,0) .. controls (3.31,0.3) and (6.95,1.4) .. (10.93,3.29)   ;

    % % Text Node
    \draw (20,32) node [anchor=north west][inner sep=0.75pt]   [align=left] {$\displaystyle A_{i}$};
    %Straight Lines [id:da8618825922990576] 
    \draw    (30,56.67) -- (30,101) ;
    \draw (0,110) node [anchor=north west][inner sep=0.75pt]   [align=left] {\color{blue}\scriptsize reputation is};
    \draw (0,122) node [anchor=north west][inner sep=0.75pt]   [align=left] {\color{blue}\scriptsize updated by};
    \draw (0,134) node [anchor=north west][inner sep=0.75pt]   [align=left] {\color{blue}\scriptsize advertisements};

    % % Text Node
    \draw (55,32) node [anchor=north west][inner sep=0.75pt]   [align=left] {$\displaystyle S_{i}$};
    %Straight Lines [id:da6261558199591686] 
    \draw    (80,41) -- (380,41) ;
    %Straight Lines [id:da02668457926064649] 
    \draw    (380,41) -- (380,84.33) ;
    \draw [shift={(380,86.33)}, rotate = 270] [line width=0.75]    (10.93,-3.29) .. controls (6.95,-1.4) and (3.31,-0.3) .. (0,0) .. controls (3.31,0.3) and (6.95,1.4) .. (10.93,3.29)   ;

    % % Text Node
    \draw (110,88) node [anchor=north west][inner sep=0.75pt]   [align=left] {$\displaystyle R_{i}^{'}$};
    % r' to d:
    \draw    (135,101) -- (229,101) ;
    \draw [shift={(231,101)}, rotate = 180] [line width=0.75]    (10.93,-3.29) .. controls (6.95,-1.4) and (3.31,-0.3) .. (0,0) .. controls (3.31,0.3) and (6.95,1.4) .. (10.93,3.29)   ;
    \draw (150,110) node [anchor=north west][inner sep=0.75pt]   [align=left] {\color{blue}\scriptsize demand is};
    \draw (150,122) node [anchor=north west][inner sep=0.75pt]   [align=left] {\color{blue}\scriptsize realized};
    % r' to next r:
    % down 
    \draw    (125,118) -- (125,200) ;
    % right
    \draw    (125,200) -- (500,200) ;
    % up
    \draw    (500,200) -- (500,120) ;
    \draw [shift={(500,120)}, rotate = 450] [line width=0.75]    (10.93,-3.29) .. controls (6.95,-1.4) and (3.31,-0.3) .. (0,0) .. controls (3.31,0.3) and (6.95,1.4) .. (10.93,3.29)   ;

    % % Text Node
    \draw (236,92) node [anchor=north west][inner sep=0.75pt]   [align=left] {$\displaystyle D_{i}$};
    % d to f: 
    \draw    (264,100) -- (350,100) ;
    \draw [shift={(350,100)}, rotate = 180] [line width=0.75]    (10.93,-3.29) .. controls (6.95,-1.4) and (3.31,-0.3) .. (0,0) .. controls (3.31,0.3) and (6.95,1.4) .. (10.93,3.29)   ;
    \draw (270,60) node [anchor=north west][inner sep=0.75pt]   [align=left] {\color{blue}\scriptsize firm operates, and};
    \draw (270,75) node [anchor=north west][inner sep=0.75pt]   [align=left] {\color{blue}\scriptsize fill rate is realized};
    
    % % Text Node
    \draw (370,92) node [anchor=north west][inner sep=0.75pt]   [align=left] {$\displaystyle F_{i}$};
    % f to next r:
    \draw (400,101) -- (480,101) ;
    \draw [shift={(480,101)}, rotate = 180] [line width=0.75]    (10.93,-3.29) .. controls (6.95,-1.4) and (3.31,-0.3) .. (0,0) .. controls (3.31,0.3) and (6.95,1.4) .. (10.93,3.29)   ;
    \draw (393,125) node [anchor=north west][inner sep=0.75pt]   [align=left] {\color{blue}\scriptsize reputation is};
    \draw (393,140) node [anchor=north west][inner sep=0.75pt]   [align=left] {\color{blue}\scriptsize updated by WoM};
    \draw (393,155) node [anchor=north west][inner sep=0.75pt]   [align=left] {\color{blue}\scriptsize communication on};
    \draw (393,170) node [anchor=north west][inner sep=0.75pt]   [align=left] {\color{blue}\scriptsize the actual fill rate};

    \draw (490,92) node [anchor=north west][inner sep=0.75pt]   [align=left] {$\displaystyle R_{i+1}$};

    \end{tikzpicture}
    }
    \caption{Endogenous Demand Structure}
    \label{fig:endo_demand}
\end{figure}
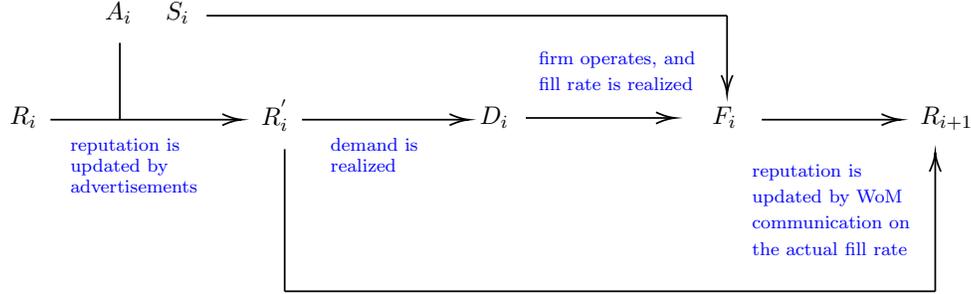

%% file: sections/structural_results.tex
\section{Analysis}
\label{sec:3}
In this section, we present our structural results that partially characterize the operational decisions under WoM communication and provide insights. We also discuss the \emph{value of information (VoI)} in terms of the profit difference resulting from integrating the WoM communication in decision making or not.
\subsection{Variable Capacity Firms}
\label{var.cap.firms.analysis}
Our first result is regarding the optimal capacity allocation, and is applicable for both aware and naive firms.

\newtheorem{Proposition}{Proposition}
\begin{Proposition}
In the variable-capacity problem, in each period, the optimal capacity is equal to the demand to be realized for both aware and naive firms.
\label{Prop:SiDi}
\end{Proposition}
Proof of the proposition is given in Appendix \ref{Proof:Prop:SiDi}.

Proposition \ref{Prop:SiDi} simplifies the optimal capacity decision in the problem to $S^*_i=D^*_i$ for all $i$ and for both aware and naive firms. Therefore, the problem reduces to finding the optimal advertisement policies. %Since the firm has to choose between discrete advertisement levels at each period, the number of solutions is finite, which implies complete enumeration as an obvious solution method. However, the number of iterations required grows exponentially with the number of periods $N$. The problem complexity is $O(2^N)$. Having said that, let us pause the discussion on the solution method, and continue with the structural analysis and the resulting implications.

An important implication of Proposition \ref{Prop:SiDi} is that we have $F_i=1$ for all periods. Therefore, as a result of WoM communication, the reputation of the company increases throughout the horizon. Hence, regarding the relationship between the reputation and demand between consecutive periods, we can write the following: $R^*_{i+1} \geq R^*_i$, $D^*_{i+1} \geq D^*_i$.

Accordingly, we have the following corollary.

\newtheorem{Corollary}{Corollary}
\begin{corollary}
In the variable-capacity problem, total profit is a non-decreasing function of the starting reputation $R_1$ for both aware and naive firms.
\label{start w big R}
\end{corollary}

Corollary \ref{start w big R} results from the fact that we have $R^*_{i} \geq R^*_{j}$ for all $i>j$; higher reputation values lead to higher demand and greater profit as a result.

Our next result is regarding the optimal advertisement policy of the aware firm.

\begin{Proposition}
\label{Prop:SwitchingPoint}
Under variable capacity, for the aware firm, if $A^*_{i}=L$, then $A^{*}_j=L$ for all $j>i$.
\end{Proposition}
The proof is in Appendix \ref{Proof:Prop:SwitchingPoint}. Proposition \ref{Prop:SwitchingPoint} states that under variable capacity, in an optimal advertising sequence, firm never makes high level advertisement after making low level advertisement for the first time. Therefore, an optimal advertising sequence should consist of a number of high level advertisements in the earlier periods and continue with low level advertisements after a point. Thus the problem becomes when to make the first low level advertisement, which reduces the order of complexity from $O(2^N)$ to $O(N)$.

We call the first period we observe a low level advertisement, \emph{the switching point}. For example, in a problem with $N=5$ and $\{H, H, L, L, L\} $ as the advertising sequence, the switching point would be 3.

%In Proposition \ref{Prop:SiDi}, we have seen that it is always more profitable for the firm to fully meet the demand and exploit positive WoM. When Proposition \ref{Prop:SwitchingPoint} is considered, it is seen that firm prefers to boost consumer expectations through advertisement when the returns of advertisement are high, that is at the earlier periods when $R_i$ is the lowest. Then, in the later periods, the firm does not invest in advertising efforts as much, since advertisement's marginal return on reputation is lower at high $R_i$. Moreover, extra effort put into advertisement at any period $i$ will contribute to $R_j$ for every $j > i$.

Propositions \ref{Prop:SiDi} and \ref{Prop:SwitchingPoint} together imply that under variable capacity, the optimal strategy suggests aggressive growth. Firm builds up expectations by advertisement and then builds them up even further by fully meeting them.

Next, we analyze the optimal strategy of the naive firm under variable capacity, which is rather straight-forward, as the firm picks the advertisement level that results in a higher $\pi_i$ than the alternative in each period.

The optimal decision of the naive firm under variable capacity, for any period $i$, can be summarized as below, where $\iota$, $\rho$, $\kappa$, $\tau$ and $\nu$ are given by Appendix \ref{thresholds:naive}.
\begin{equation}
	 \small{
	 	 A_i = 
		 \begin{cases} 
			H, & if \:\: \tau > 0 \:\: and \:\: \iota < \rho \:\: and \:\: R_i \: \in \: \left(\:min \: \{ max\{ \iota, \: \kappa \}, \: \rho \}, \: max\{\rho, \: \nu\} \: \right) \\
			L, & otherwise \\
		\end{cases}
	}
\label{varCapacityIncDecision}
\end{equation}
Note that since we always have $R_{i+1}\geq R_i$, as implied by Proposition \ref{Prop:SiDi}:
\begin{itemize}
    \item [(i)] firm will make low level advertisement when $R_i \leq min \: \{ max\{ \iota, \: \kappa \}, \: \rho \}$,
    \item [(ii)] then, as reputation increases between consecutive periods, $min \: \{ max\{ \iota, \: \kappa \}, \: \rho \} < R_i < max\{\rho, \: \nu\}$ will hold at some later period, and the firm will start making high level advertisements,
    \item [(iii)] finally, reputation will exceed $max\{\rho, \: \nu\}$, and firm will continue with low level advertisements until the end of the last period.
\end{itemize}
Thus, firm will adopt an advertisement pattern such as: $A = \{L,\: L,\: L,\: H, \: H, \: H, \: L, \: L\}$

We say that such solutions have \emph{two switching points}. In the above example, where $N=8$, the first switching point is at $j=4$ and the second is at $j=7$. Depending on the starting reputation, $R_1$, firm may adopt a solution with one or two switching points. While the structure of the optimal advertisement policies are different for aware and naive firms, there is a relation between the optimal advertisement strategies of the two. The number of high level advertisement periods of the naive firm's problem creates a bound for that of the aware firm's problem, which leads to the following.

\begin{Proposition}
	In the variable capacity problem, the naive firm never makes more high level advertisements than its aware equivalent.
	\label{Prop:LBForH}
\end{Proposition}

The proof is presented in Appendix \ref{Proof:Prop:LBForH}. Proposition \ref{Prop:LBForH}, implies that the aware firm may prefer high advertisement costs as an investment on the reputation, and be rewarded in the upcoming periods. However, the naive firm makes low level advertisement if the returns on advertisement investment will not be realized immediately.This is due to the \emph{myopic behavior} of the naive firm.

%Using Proposition \ref{Prop:LBForH}, we can reduce the number of iterations needed to find the optimal advertising sequence of the aware firm even further from $N+1$ to $N - H_{R_1}^I + 1$.

We have observed that company can make more informed decisions under complete information on WoM communication and can interfere in customer demand. As a result, the aware firm is able to gain reputation faster and gain higher profits.

\subsection{Constant Capacity Firms}
Constant capacity problem introduces new challenges into the advertisement decision. Since the firm cannot change service capacity, advertisement is the only tool to manage demand. Due to the discrete advertisement choices, it is possible that the resulting reputation may attract demand that is larger or smaller than what the firm can serve. By the advertising actions, the firm may end up increasing the reputation more than necessary, and suffer negative WoM due to low fill rate. When reputation becomes low enough, the firm will start enjoying the positive WoM as a result of low demand. Therefore, there may be fluctuations in the reputation, and the single or two switching point solutions are often not optimal for the constant capacity firm.

Our analysis is based on the conditions of the firm at the beginning of each period. We define a set of states that depend on firm's reputation before advertisement.  At the beginning of a period $i$, the firm can be in one of the following states:

\textbf{State A:} Demand exceeds capacity regardless of the advertisement decision ($D_i \geq S$).
        
\textbf{State B:} Demand exceeds capacity under high level advertisement, and is positive but below capacity under low level advertisement ($D_i \geq S$ when $A_i = H$, and $S > D_i > 0$ when $A_i = L$).

\textbf{State C:} Demand is positive and less than the capacity regardless of the advertisement decision ($S > D_i > 0$).

\textbf{State D:} Demand exceeds capacity under high level advertisement and is zero under low level advertisement ($D_i \geq S$ when $A_i = H$, and $D_i = 0$ when $A_i = L$).

\textbf{State E:} Demand is positive but less than the capacity under high level advertisement, and is zero under low level advertisement ($S > D_i > 0$ when $A_i = H$, and $D_i = 0$ when $A_i = L$).

 \textbf{State F:} Demand is zero regardless of the advertisement decision ($D_i = 0$).

A summary of these states is presented in Table \ref{Tab:States}.

\begin{table}[H]
\caption{Demand in each state based on the advertisement level.}
    \centering
    \begin{tabular}{c|c|c}
    \hline \hline
    State & Low Advertisement & High Advertisement \\ \hline
    A & $D_i\geq S$ &  $D_i\geq S$ \\\hline
    B & $ S>D_i > 0$ &  $D_i\geq S$ \\ \hline
    C & $ S>D_i > 0$ &  $ S>D_i > 0$ \\ \hline 
    D & 0 &  $D_i\geq S$ \\ \hline 
    E & 0 &  $S > D_i>0 $ \\ \hline
    F & 0 & 0 \\\hline \hline
    \end{tabular}
    \label{Tab:States}
\end{table}

We denote the set of possible states a firm can observe throughout the planning horizon by $\Omega$. Note that $\Omega$ depends on the problem parameters. We discuss how to determine the set $\Omega$ based on the parameters in Appendix \ref{App:Omega}.

Note that states with higher reputation values, such as A, yield lower fill rates than the others. For states C, E and F, we have $R_{i+1} > R_i$ for period $i$. Same statement holds for states B and D if $A_i = L$. Then we can make interpretations such as: the state in period $i+1$ cannot be $E$ or $F$ when state in period $i$ is C. When we are at state A or when $A_i = H$ in states B and D; the relationship between $R_i$ and $R_{i+1}$ may differ depending on the advertisement choice and the fill rate. Based on these state definitions, we will now present some structural results for the optimal advertisement policies for the firms with constant capacity.

Under constant capacity, the aware firm's problem is not tractable. Therefore, we partially characterize the aware firm's advertisement policy. Our first result is in Proposition \ref{Prop:ConstCapCEF}, and proven by Appendix \ref{Proof:Prop:ConstCapCEF}.

\begin{Proposition}
	Under constant capacity, if $\Lambda \left(1-\frac{u+p}{m}\right) \leq S$; then $A^*_{j} = L$ for all $j > i$ if $A^*_i = L$ for the aware firm.
	\label{Prop:ConstCapCEF}
\end{Proposition}

By Proposition \ref{Prop:ConstCapCEF}, if demand can never exceed the service capacity, the aware firm's problem has a single switching point from $H$ to $L$ as in the variable capacity problem.

%When this condition holds, we say that demand can never reach or exceed capacity, and we have $\Omega \subseteq \{C, \: E, \: F\}$. This allows us to decrease complexity from $O(2^N)$ to $O(N)$.

The next result is regarding the states with high demand potential, such as A, B and D. By Lemma \ref{Lemma:ConstCapABD}, when the set of states a firm can be in (i.e., $\Omega$) is a subset of $\{A,B,D\}$, we can calculate an upper bound for the number of high level advertisements in the optimal solution for the aware firm. The proof of the lemma is provided in Appendix \ref{Proof:Lemma:ConstCapABD}.

\begin{lemma}
	Under constant capacity, if the firm can always achieve $D_i\ge S$ by high level advertisement, the optimal number of high level advertisements in the planning horizon for the aware firm cannot be greater than the high level advertisements determined by the following rule:
	%if we have $D_i \geq S$ whenever $A_i=H$, the optimal number of high level advertisements of the aware firm cannot be greater than the one obtained by the following rule:
	choose low level advertisement in state A, and high level advertisement otherwise.
	\label{Lemma:ConstCapABD}
\end{lemma}

%Lemma \ref{Lemma:ConstCapABD} implies that when the firm is able to seize enough demand to utilize service capacity just by a single advertisement action at each period, making high level advertisement only when demand needs a boost to exceed service capacity provides a lower bound to the optimal amount of profit. Having complete information on the endogenous demand structure, firm can be able to come up with less costly advertising policies to keep demand high.

%Lemma \ref{Lemma:ConstCapABD} reduces the number of iterations required from $2^N$ (or $\sum_{k=0}^{N} \binom{N}{k}$) to $\sum_{k=0}^{h^*} \binom{N}{k}$.

For our final special case, we present Proposition \ref{Prop:AbsorbingA}, which is proven by Appendix \ref{Proof:Prop:AbsorbingA}.

\begin{Proposition}
	Suppose, under constant capacity, a firm is at state A in some period $i$. Then, if the intervals $[X_1, \: 1)$, where $X_1=\frac{\frac{(u+p) \Lambda}{m (\Lambda -S)}-L^\alpha}{1-L^\alpha}$, and the one defined by the below expression, where $A_i=L$, do not have an intersection, the aware firm makes low level advertisement in all remaining periods.
    \begin{equation}
    	\frac{\left(\frac{- \left(\frac{w S}{\Lambda} - X_1 - (1-w) \frac{u+p}{m}\right) \pm \sqrt{\left(\frac{w S}{\Lambda} - X_1 - (1-w) \frac{u+p}{m}\right)^2 - 4 (1-w) X_1 \frac{u+p}{m}}}{2 (1-w)}\right) - A_i^\alpha}{1 - A_i^\alpha}
    	\label{roots for R to drop from a}
    \end{equation}
	\label{Prop:AbsorbingA}
\end{Proposition}

Note that Proposition \ref{Prop:AbsorbingA} is relevant under the condition that $A \in \Omega$ (state $A$ can be observed) and $\Omega \neq \{A\}$ (at least one other state can be observed). Therefore, we need to have the reputation level yielding $D_i \geq S$ regardless of the advertisement level ($X_1$) between $0$ and $1$, which translates as follows.
\begin{equation}
	1 > \frac{u+p}{m} \frac{\Lambda}{\Lambda - S} > L^\alpha
	\label{cond_for_prop_all_l_after_a}
\end{equation} 

When the roots of Equation \ref{roots for R to drop from a} are complex, negative or define an interval that simply does not coincide with $[X_1,\:1)$, Proposition \ref{Prop:AbsorbingA} is applicable. At low $w$, the roots of Equation \ref{roots for R to drop from a} define a large interval but they are negative. At slightly higher $w$, the roots are positive but the interval between the roots is narrower, so it does not intersect $[X_1,\:1)$. At larger $w$, roots are complex numbers. Therefore, up to this point, the proposition is applicable. Finally, when $w$ increases further, roots become positive real numbers, and the interval between them starts to get wider, making the proposition not applicable. The same behavior is observed for $\alpha$ and $\frac{u+p}{m}$, while $L$ behaves in the exact opposite manner. To sum up, when Equation \ref{cond_for_prop_all_l_after_a} holds, and $w$, $\alpha$ and $\frac{u+p}{m}$ are low while $L$ is high; it is easy for the firm to keep demand high once it is raised above a certain level. Consumers' tend to stick with their existing opinions rather than relying on WoM. This may occur for services, for which WoM communication has not become a habit yet. Alternatively, highly needed or beneficial services may invite such situations, as well as the ones without substitutes.

Next, we analyze the optimal strategy of the naive firm under constant capacity, which is again very straight-forward, as the firm simply maximizes $\pi_i$. The optimal decision for any period $i$ can be summarized as below, where the values of $\gamma$, $\omega$, $\phi$ and $\beta$ are given by Appendix \ref{thresholds:naive}.
\begin{equation}
	A_i = 
	\begin{cases} 
		H, & \text{if} \:\: \:
		\begin{split}
			& R_i \: \in \: \left( \gamma, \: min\{X_2, \: \beta\}\right) \:\: \text{or} \\
			& R_i \: \in \: \left(max\{Y_1, \: X_2\}, \: \omega\right) \:\: \text{or}\\
			& \phi > 0 \:\: \text{and} \:\: X_2 < R_i < Y_1\\
		\end{split} \\
		
		L,  & \text{otherwise} \\
	\end{cases}
	\label{constCapacityIncDecision}
\end{equation}

Notice that the naive firm's optimal advertisement decision under constant capacity is more complex than its variable capacity equivalent. In the variable capacity scenario, the firm picks high level advertisement only if $R_i$ falls within the interval defined in Equation \ref{varCapacityIncDecision}. Under constant capacity, however, by Equation \ref{constCapacityIncDecision}, we define a more complex decision mechanism. This results from the fact that under constant capacity, we can observe $S > D > 0$ and $D>S$, which is not the case in the variable capacity scenario. Moreover, unlike the variable capacity scenario, we can observe instances, where the naive firm makes more frequent high level advertisements than the aware firm. This highlights the importance of having complete information on endogenous demand structure; since, by having complete information, firms can make more profit even if they place less effort in advertising.

A final note on the naive firm under constant capacity is that when $\Omega \subseteq \{C,\:E,\:F\}$ (i.e., when the firm cannot possibly generate enough demand to fully utilize its capacity), the naive firms' decision mechanism (given by Equation \ref{constCapacityIncDecision}) becomes equivalent to the one given by Equation \ref{varCapacityIncDecision} under variable capacity.

It is clear that aware firms are able to generate more profit than their naive equivalents both under variable and constant capacity scenarios. Depending on the parameters, value of information (VoI) may be larger, which indicates how advantageous the aware firm is. Having said that, access to reliable information on customer behavior (e.g., sensitivity to advertisement, reliance on past information) does not come at no cost. Therefore, it is important to understand the factors affecting the VoI. To explore the relationship between the parameters and VoI, we perform extensive computational studies presented in Section \ref{sec:4}.

%% file: sections/computational_results.tex
\section{Computational Results and Discussion}
\label{sec:4}
We perform computational experiments to observe the impact of problem parameters on the VoI, which are the starting reputation ($R_1$), low advertisement level ($L$), high advertisement level ($H$), customers' advertisement resistance ($\alpha$), weight of WoM ($w$), market size ($\Lambda$), capacity cost ($c$), reputation's monetary value ($m$), reservation utility ($u$), price ($p$), cost of advertisement ($c_A$), number of periods ($N$) and service capacity ($S$).

We use the parameter distributions in Table \ref{tab:prm_distr_var} when generating the test instances. Using these distributions, we first generate 100 random scenarios for each problem (i.e., for variable capacity problem and for constant capacity problem). To test the impact of a parameter, for each parameter set, we generate 50 of random values for the parameter we are testing using the original distribution while keeping others constant. In total, we consider 60,000 instances for each problem as we test 12 parameters ($S$ is excluded for the variable capacity problem and $c$ is excluded for the constant capacity problem). We solve these test instances for both aware and naive firms, and observe how VoI changes with each parameter. We perform all computations using Python Programming Language version 3.8.5 on an Intel(R) Core(TM) i7-8700 CPU @ 3.20 GHz, 16GB RAM Windows 10 PC.

The parameter distributions are determined according to a preliminary analysis to avoid trivial cases. Note that the distribution of the problem horizon, $N$, is different for variable and constant capacity scenarios since the computation time increases significantly for the constant capacity problem.

\input{distributions/prm_distr_var}

Finally, to make a healthy comparison regarding VoI, we need to eliminate the differences resulting from the size of the problem. Therefore, we base our comparison on \% VoI, which is calculated by:
\begin{equation}
    \% \: VoI = \begin{dcases}
    \frac{\Pi^*_{aware} - \Pi^*_{naive}}{\Pi^*_{aware}} \: 100, & if \:\: \Pi^*_{aware}>0\\
    0, & if \:\: \Pi^*_{aware}=0\\
    %100, & if \:\: \Pi^*_{naive}=0 \text{ and } \Pi^*_{aware}>0
    \end{dcases}
\end{equation}
We assume that the profit is zero for the cases where the total profit in the optimal solutions of the two types of firms (aware and naive) is negative since a firm would not operate in such cases.

\subsection{Value of Information Under Variable Capacity}
In this section, we compare the aware and naive firms under variable capacity, and discuss the results for each parameter. A key indicator of how \% VoI changes is the location of the switching points, introduced in Section \ref{var.cap.firms.analysis}. For the sake of brevity, we will refer to those as follows.
\begin{itemize}
    \item $SP_0$: the period number in which the aware firm switches from $H$ to $L$.
    \item $SP_1$: the period number in which the naive firm switches from $L$ to $H$.
    \item $SP_2$: the period number in which the naive firm switches from $H$ to $L$.
\end{itemize}

In the following part, we present and discuss the VoI results for the variable capacity problem. Note that under variable capacity, the reputation is increasing for both firms, but at a faster rate for the aware firm. Therefore, the contribution of the early periods to the VoI is larger. The difference between two firms is greater if the naive firm is discouraged from advertisement in the problem setting (i.e., by high advertisement resistance, high reservation utility and similar factors). 

\paragraph{Impact of Advertisement Resistance ($\alpha$) on VoI:} We present a representative set of computational results for different advertisement resistance parameters ($\alpha$) in Figure \ref{fig:var_cap_alpha_prft} displaying the profits obtained by both types of firms and the VoI percentages. Corresponding optimal advertisement sequences are given in Appendix \ref{exp:varcap:alpha} along with the parameter values. As discussed earlier, when $\alpha$ increases, consumers become more resistant towards advertisement. Notice that the VoI is zero for very low and very high advertisement resistance. The reason behind the former is that when the consumers are easily impressionable, advertisement decision is less impactful, and the naive firm is not disadvantaged by the lack of WoM knowledge. For very large $\alpha$ values, it becomes rather difficult to attract more customers through advertisements for both firms, bringing the impact of WoM knowledge gap to zero. On the other hand, for intermediate advertisement resistance levels, the VoI is significant. In this region, as $\alpha$ increases, the aware firm tries to penetrate this resistance through making high level advertisements for a longer period of time, delaying $SP_0$ in the planning horizon.
Meanwhile, the naive firm is reluctant to invest in advertisement since it cannot foresee the benefits of it due to myopic behavior. Therefore, the naive firm operates longer under low advertisement policy (i.e., $SP_1$ moves forward in the planning horizon) as $\alpha$ increases. As a result, it takes longer to build on reputation, and $SP_2$ is delayed. While the profits of both firms decline with the advertisement resistance, the naive firm's profit decreases faster due to the lack of foresight, resulting in higher VoI before dropping to zero.
\input{charts/var_cap_alpha_profits}

Among other problem parameters, the reservation utility $u$ affects the optimal advertisement decisions and the VoI in a similar manner to the advertisement resistance, while the monetary value of reputation $m$ works in the opposite direction. For brevity, we do not discuss them in detail.

%%%%%%%%%%%%%%%%%%%%%%%%%%%%%%%
\paragraph{Impact of the Starting Reputation ($R_1$) on VoI:}
We display the results of a representative computational setting for the impact of initial company reputation in  Figure \ref{fig:var_cap_r1_prft_sample1}. Corresponding optimal advertisement sequences and the parameter set are given in Appendix \ref{exp:varcap:r1_sample1}. We observe that for low initial reputation levels (e.g., $R_1\le 0.4$ in this example), none of the firms can generate profit. This is the region where even the aware firm cannot recover from this low reputation. As the initial reputation increases, we see that the aware firm sees the potential to build up using high level advertisements. Meanwhile the naive firm cannot capture this opportunity due to myopic behavior, and the VoI jumps to 100\%. With further increase in the initial reputation, it also becomes profitable for the naive firm to operate, and hence the VoI starts to decrease until it reaches zero at a point where the optimal policy of the naive firm coincides with that of the aware firm. For larger company reputations, the reputation itself acts as the main driving force to attract more customers decreasing the impact of the advertisement. Therefore, the gap between the aware firm and the naive firm decreases, and the VoI converges to zero as the initial reputation increases.

%At significantly low $R_1$, the firms may fail to realize profit (depending on the other parameters, this region may not be visited); therefore \% VoI is 0. As we increase $R_1$ and repeat the experiment, once some $R_1$ is reached, it becomes optimal for the aware firm to enter the market by making a number of high level advertisements; therefore, $SP_0$ jumps to some time period. At that point we have $\% VoI = 100\%$. As $R_1$ is further increased, the need for advertising drops, and $SP_0$ moves to earlier periods. Meanwhile, the naive firm starts to realize profit at some $R_1$, and \% VoI starts to drop from 100\%. It is observed that as $R_1$ goes up, the gap between the two firms' reputation curves becomes narrower (see Figure \ref{fig:r_i_var_cap}), leading to similar demand realization and similar profits. As $R_1$ goes up, reputation grows faster, and it becomes easier for the switching conditions to be satisfied; therefore, both $SP_1$ and $SP_2$ move backwards. At a large enough $R_1$, the necessity for advertisement is gone for the naive firm, and for an even larger $R_1$ the same happens for the aware firm, leading to zero \% VoI. 

\input{charts/var_cap_r1_profits}

%We present the results of another experiment in Appendix \ref{exp:varcap:r1_sample2} that captures the motion of $SP_1$ and $SP_2$.

\paragraph{Impact of Price ($p$) on VoI:} We present the computational results for a representative case for the changes in price in Figure \ref{fig:var_cap_p_prft}. The corresponding parameter set and the optimal advertising sequences are given in Appendix \ref{exp:varcap:p_sample1}. %The profit function under variable capacity is concave in price in positive demand region due to the utility-based demand. 
The availability of the WoM knowledge creates a significant difference between the aware and naive firms as price changes. Increasing prices first encourage both firms by increasing the revenues but can quickly become discouraging for the naive firm as it becomes more difficult to attract demand. We observe that the naive firm stops operating during mid-price ranges. On the other hand, there is a potential for higher profit for larger prices as long as demand can be maintained. This opportunity is available to the aware firm through WoM knowledge. As the price increases, the aware firm's advertisement strategy becomes more aggressive with more frequent high level advertisements as it can be seen in Appendix \ref{exp:varcap:p_sample1}. Therefore, the aware firm can operate in a much larger price region. Indeed, in some test instances, the aware firm can achieve its highest profit at a price value that is too high for the naive firm to operate as in Figure \ref{fig:var_cap_p_prft}.

%When price is significantly low, it is not profitable to operate for neither of the firms, leading to zero \% VoI. As we increase $p$ between test instances, at some point the aware firm is able to generate profit and, hence, starts to operate making $\% VoI=100\%$. As price goes up further, the naive firm starts to operate as well, lowering \% VoI. For a while, due to price's positive impact on revenue, even though demand is lower, both firms' profits go up. However, the increase in profit slows down after a while for both firms starting with the naive firm, and \% VoI starts to go up. Eventually, profits start to decline, and finally, the negative impact of price on demand outweighs its improvement on revenue so much that the naive firm is unable to operate, making $\% VoI=100\%$ again. In the end, the aware firm stops operating as well, and \% VoI takes a final value of $\% 0$.  

%Another sample is given by Appendix \ref{exp:varcap:p_sample2} that captures the motion of $SP_1$ and $SP_2$.
\input{charts/var_cap_p_profits}

\paragraph{Discussion on the Other Parameters:}

In the computational results, we also observe that the weight customers assign to the fill rate information obtained by WoM communication $w$ is not among the key parameters affecting the VoI. Considering that the main difference between the aware and the naive firms is the knowledge on WoM, this observation is counter-intuitive. This is due to the nature of the variable capacity. Under variable capacity, both firms choose to meet the entire demand, and achieve the fill rate of 1. Therefore, both firms' customers share the same information through WoM, leaving little room for the coefficient $w$ to act.

%Under variable capacity, it is optimal for both firms to fulfill the entire demand; therefore, we only observe positive WoM. $w$ acts as a means of advertisement once demand is realized for the first time, and both firms' reputation curves converge faster under high $w$. In the end, \% VoI decreases rather insignificantly with $w$. % (We may observe a rare special case with very low advertisement levels and very high $\alpha$, where the firms would depend solely on $w$ to increase reputation. In such a case, the firms may not operate at low $w$, and \% VoI may be $0\%$, $100\%$ and then start decreasing.).

Increasing capacity cost $c$ and advertisement cost $c_A$ decrease the incentive for advertisement to attract more demand. As these cost parameters increase, the naive firm avoids the high level advertisements ($SP_1$ is delayed in the planning horizon), enlarging the reputation gap between the naive firm and the aware firm. However, this does not translate into an increasing profit gap due to the higher operational costs for the aware firm. Therefore, the impact of the capacity cost and the advertisement cost on the VoI is not significant while it is for the advertisement policy.

It is seen that the aware firm is inclined to increase advertisement in longer problem horizons, since there is simply more time to collect the returns; therefore, $SP_0$ moves forward. The naive firm is not able to develop long-term advertising strategies; therefore, its switching points do not move with $N$. Yet, as the planning horizon gets longer, the advantage gained by the aware firm in the earlier periods become less significant as the $\pi_i$'s of the two firms collected in late periods are more similar. As a result, $N$ reduces the VoI.

In terms of the advertisement, as the high advertisement level $H$ increases, any difference on the advertisement policy between the naive and aware firms becomes more important. However, the aware firm avoids high level advertisements for larger $H$ due increasing cost. Therefore, increasing $H$ does not have a significant or consistent impact on the VoI. For the low advertisement level $L$, increasing $L$ reduces the gap between the aware and the naive firms for the WoM knowledge to play a role. Therefore, VoI decreases with $L$.

Finally, as the market size $\Lambda$ increases, the naive firm follows a more aggressive advertisement policy. Although this decision is still due to myopic approach, it boosts the reputation growth decreasing the VoI in larger markets.

%when the market size, $\Lambda$ gets larger, the naive firm is more willing to invest in high level advertisement in the early periods, boosting the reputation growth. This indicates that VoI decreases in $\Lambda$.

\subsection{Value of Information Under Constant Capacity}

In this section, we present and discuss the computational results under constant capacity for aware and naive firms. How the VoI on the WoM communication changes with respect to a parameter depends highly on the set of states the firm can be in as a function of the parameters ($\Omega$). Our discussion on the VoI revolves around this concept. Note that when $\Omega \subseteq \{C,\:E,\:F\}$, the parameters impact the VoI in the same way as they do for the variable capacity problem due to Proposition \ref{Prop:ConstCapCEF}. Therefore, we focus on other cases in this section.

% Moreover, if the content of $\Omega$ changes with the parameter of concern, the relationship between the parameter and VoI differs for different intervals of the parameter. Consequently, another thing to note is that when $\Omega \subseteq \{C,\:E,\:F\}$, all parameters impact the VoI in the same way as they do for the variable capacity problem as a result of Proposition \ref{Prop:ConstCapCEF}. Therefore, this section focuses more on the other $\Omega$ types for the sake of brevity.

%Similar to before, our discussion revolves around the value of information on the WoM communication. 
%In this section, we compare the aware and naive firms under constant capacity, and report how the parameters impact VoI.
%changing a parameter impacts $\Omega$. 
%For instance, suppose $\Omega = \{A,\:B,\:D\}$. For the same parameter set, if the reputation's monetary value $m$ were smaller, it would be harder to realize demand, and $\Omega$ would become $\{A,\:B,\:D,\:E\}$ or $\{A,\:B,\:D,\:E,\:F\}$. These indicate more complex problems and may require more elaborate advertisement policies, making it harder for the naive firm to approach the aware firm's profit with its simpler decision mechanism (given by Equation \ref{constCapacityIncDecision}).
%

Under constant capacity, one of the key observations is that the reputation of the naive firm may have a cyclic behavior when customers are highly sensitive to the information obtained by the WoM communication. In such cases, the naive firm cycles between attracting a small demand (resulting in high fill rate and improvement in reputation) in one period, and attracting a large demand (resulting in low fill rate and a decrease in reputation) in the next period. A representative example is presented in Figure \ref{r_const_cap_fluct} for the parameters $R_1=0.18$, $L=0.68$, $H=0.82$, $\alpha=6.04$, $w=0.99$, $\Lambda=653$, $c=0.70$, $m=39.22$, $u=2.70$, $p=16.13$, $c_A=731.85$, $N=18$, $S=88$. Note that in this particular example, the two firms' advertisement costs are the same, yet the policies are different. The aware firm is able to stay in high demand states longer while naive firm's investments on reputation (through advertisement) is immediately overridden by negative WoM.

%We represent a problem instance below to illustrate a typical situation observed under constant capacity. Here, the naive firm is greatly impacted by the WoM and the improvements on reputation get overridden every period while the aware firm's reputation does not fluctuate as much, yielding higher demand. In this particular example, the two firms' advertisement costs are the same, yet the timing is different. With the same total cost, the aware firm visits the state $A$ 11 times, $B$ 3 times, $D$ twice, $E$ and $F$ once each while the naive firm visits state $A$ 8 times, $E$ 8 times, and $F$ twice (the others are not visited) yielding less revenue. The corresponding parameter values are: $R_1=0.18$, $L=0.68$, $H=0.82$, $\alpha=6.04$, $w=0.99$, $\Lambda=653$, $c=0.70$, $m=39.22$, $u=2.70$, $p=16.13$, $c_A=731.85$, $N=18$, $S=88$.

%\begin{minipage}[b]{0.45\textwidth}
\input{charts/con_cap_fluct}
%\end{minipage}
%\begin{minipage}[b]{0.10\textwidth}

%\end{minipage}
%\begin{minipage}[b]{0.45\textwidth}
%\input{charts/con_cap_pi.tex}
%\end{minipage}

%Finally, let us note that since capacity is constant, the impact of capacity cost $c$ will not be studied.

\paragraph{Impact of the Service Capacity ($S$) on VoI:}

 Under constant capacity, the impact of the capacity level on the VoI is not trivial, and it depends on the other parameters. In general, the extreme capacity regions (very small or very large) are not the most profitable. When the capacity is very low, the firms cannot recover from the negative WoM; when the capacity is very large, operational costs become dominant, and both firms stop operating. In moderate capacity regions, when customers are highly sensitive to the information obtained by WoM  (i.e., $w$ is large), the profits of both firms have an increasing trend as the capacity increases. Since customers are less susceptible to advertisement under high $w$, the naive firm's myopic advertisement policy does not hurt the overall profits, and the VoI reduces as capacity increases. An example case is presented in Figure \ref{fig:con_cap_s1_prft}. On the other hand, when $w$ is small, customer perception relies more on the information obtained through advertisement and promotion activities. Especially for larger capacity levels, the aware firm has a greater potential due to being able manage the endogenous demand better. Therefore, the VoI is significant and increases with capacity as shown in Figure \ref{fig:con_cap_s2_prft}. For increasing capacity, the VoI on the endogenous demand structure is also significant for cases where the unit capacity cost is large. An example case is illustrated in Figure \ref{fig:con_cap_s3_prft}. In such cases, the naive firm is immediately discouraged from operating due to the large unit capacity cost as it cannot foresee the potential profit that could be received in the long run. The aware firm is the only one that provides service with the entire profit generated based on the knowledge of endogenous demand structure. The corresponding set of parameters and the optimal advertisement policies for the experiments in Figure \ref{fig:ConstSMain} are in Appendix \ref{exp:concap:s1_sample}. All in all, the profits of both firms follow concave patterns with $S$ while VoI first declines, then goes up with $S$ within the region both firms operate.

\newcommand\resize{0.65} %Adjusts size of figures
\begin{figure}[ht]
    \centering
    \begin{subfigure}{.32\textwidth}
        \input{charts/con_cap_s1}
        \caption{Large $w$}
        \label{fig:con_cap_s1_prft}
    \end{subfigure}\hfill
    \begin{subfigure}{.32\textwidth}
        \input{charts/con_cap_s2}
        \caption{Small $w$}
        \label{fig:con_cap_s2_prft}
    \end{subfigure}\hfill
    \begin{subfigure}{.32\textwidth}
        \input{charts/con_cap_s3}
        \caption{Large $c$}
        \label{fig:con_cap_s3_prft}
    \end{subfigure}
    \caption{Profits and the \%VoI with respect to service capacity $S$ under constant capacity.}
    \label{fig:ConstSMain}
\end{figure}

% \begin{minipage}{0.45\textwidth}
% \input{charts/con_cap_s1}
% \end{minipage}\hfill
% \begin{minipage}{0.45\textwidth}
% \input{charts/con_cap_s2}
% \end{minipage}\hfill \\

%It is observed that when the service capacity is significantly low, a strict upper bound is set for the revenue and low $F_i$'s are observed frequently. Both firms fail to realize profit and stay out of the market. As $S$ gets higher, the aware firm enters the market yielding VoI $ = 100\%$, to be followed by the naive firm. VoI starts decreasing. This behavior is shown by Figure \ref{fig:con_cap_s1_prft}, whose parameters, $\Omega$s and optimal advertising sequences are given by Appendix \ref{exp:concap:s1_sample}. Note that as $S$ goes up, it becomes harder to have demand exceed capacity; therefore, the reputation range covered by states A, B and D shrink. This is reflected in $\Omega$. However, when $S$ is too large; the firms leave the market due to high capacity costs (first the naive firm, then the aware firm). This leads to an increase in VoI up to 100\% and a sudden drop to zero when the aware firm stops operating. This behavior is presented by Figures \ref{fig:con_cap_s2_prft} and \ref{fig:con_cap_s3_prft} with the respective Appendices  \ref{exp:concap:s2_sample} and \ref{exp:concap:s3_sample}. All in all, the profits of both firms follow concave patterns with $S$.

\paragraph{Impact of the Weight of WoM Communication ($w$) on VoI:}
The weight customers assign to the recent information obtained by WoM communication has significant impact on VoI. In the markets with highly variable demand potential, various states are possible for the firms to visit (i.e., high cardinality $\Omega$). Accordingly, the impact of WoM in a period can be negative or positive, depending on the fill rate. In such cases, a customer group that is more sensitive to the operational performance increases the need to better understand the endogenous demand structure, and VoI as a result. The aware firm's response to increasing $w$ depends on the other parameters. For instance, under low advertisement cost, the aware firm chooses a more aggressive advertisement strategy to respond to the negative WoM due to high demand. However, for higher advertisement costs, this action becomes less profitable, and the firm switches to a low advertisement strategy and maintains the high reputation through high fill rate. To show the impact of $w$, we present representative computational results in Figure \ref{fig:con_cap_w_prft}. The corresponding parameter set and the optimal advertisement policies of both firms are given in Appendix \ref{exp:concap:w_sample}. Notice that the VoI increases as $w$ increases. When $w$ becomes too large (around 0.72 in this example), it becomes too difficult to manage demand, and both companies stop operating.
\input{charts/con_cap_w}
%Weight of WoM, $w$, does not affect $\Omega$; yet, its impact on \% VoI is still significant when $\Omega \not \subseteq \{C,\:E,\:F\}$. This is due to the fact that with these $\Omega$ types, we can have $F_i<1$, $F_i=1$ and the case where fill rate is not applicable due to zero demand.
%Therefore, we can observe WoM's effect on reputation in both directions, by positive and negative WoM. In addition to this direction variety, when $w$ is high, fluctuations in reputation are stronger and more frequent, requiring elaborate planning for the advertisement efforts. In some instances, we observe that when $w$ is high, firm makes more high level advertisements to recover the impacts of negative WoM (typically under low $c_A$) while in other instances, firm reduces the number of high level advertisements to keep the fill rate high. The second behavior sometimes results from the fact that advertisement's effect lasts so short that it is not worth it. Also, sometimes firm tries to make use of WoM as a tool for advertising. In the end, as $w$ goes up, we observe higher VoI. A sample is given by Figure \ref{fig:con_cap_w_prft}, whose parameter set, $\Omega$ and optimal advertising sequences are given in Appendix \ref{exp:concap:w_sample}.
%\input{charts/con_cap_w}
\paragraph{Impact of Advertisement Resistance ($\alpha$) on VoI:} We present the computational results for a representative case for increasing $\alpha$ values in Figure \ref{fig:con_cap_alpha_prft}. The parameter set, the set of possible states that the firm can visit and optimal advertising sequences are given in Appendix \ref{exp:concap:alpha_sample}.
As $\alpha$ increases, it becomes more difficult to change the customer perception through advertisement. We observe from our computational results that for very small $\alpha$'s where customers are highly sensitive to the information obtained by advertisement, both aware and the naive firms can achieve the same profit. As $\alpha$ increases, the impact of the advertisement decreases, and the WoM communication becomes the main factor affecting customer perception on operational performance. Being fully informed on the endogenous demand structure, the aware firm responds to the increase in $\alpha$ by advertising more aggressively compared to the naive firm. Therefore, we see that for increasing $\alpha$ the VoI increases. For very large $\alpha$ values, customers rely almost entirely on the fill rate information provided by the WoM; and depending on the capacity, no profit can be obtained in that region. Therefore, both firms stop operation.

\input{charts/con_cap_alpha}

\paragraph{Discussion on the Other Parameters:} Remaining problem parameters do not have significant or consistent impacts on the VoI. Company's initial reputation $R_1$ does not affect the set of states $\Omega$. When the firm's demand attraction power is low (i.e., $\Omega \subseteq \{C,\:E,\:F\}$), the value of information is larger for companies with low initial reputation, as it is the case under variable capacity. When $\Omega \not \subseteq \{C,\:E,\:F\}$, we observe that $R_1$ has no regular impact on VoI, and its effect highly depends on other parameters. When WoM has notable weight (high $w$), $R_i$ fluctuates, which reduces the importance of the starting reputation. On the other hand, when the planning horizon $N$ is short, $R_1$ becomes important, since the firm can employ an advertising policy that is just enough to keep $R_i$ high until the end of the $N^{th}$ period.

For the firms with similar attributes and $\Omega \subseteq \{C,\:E,\:F\}$, the VoI decreases as the length of the planning horizon $N$ increases. In other cases, the impact of $N$ is not clear, and it is outweighed by the other parameters. When $N$ is sufficiently large, the policies and states of the firms start to repeat themselves keeping the VoI almost stable.

%To begin with the starting reputation, $R_1$, which does not affect $\Omega$, we observe that its impact on VoI is not regular or significant. When $\Omega \subseteq \{C,\:E,\:F\}$, the starting reputation decreases VoI, as it is the case under variable capacity. When $\Omega \not \subseteq \{C,\:E,\:F\}$, we observe that $R_1$ has no regular impact on VoI, and its effect highly depends on other parameters. When WoM has notable weight (high $w$), $R_i$ fluctuates, which reduces the importance of the starting reputation. On the other hand, when the planning horizon $N$ is short, $R_1$ becomes important, since the firm can employ an advertising policy that is just enough to keep $R_i$ high until the end of the $N^{th}$ period.

%Length of the planning horizon, $N$, does not have a regular impact on VoI either, unless we have $\Omega \subseteq \{C,\:E,\:F\}$. When $\Omega \not \subseteq \{C,\:E,\:F\}$, we observe that \% VoI makes small fluctuations as we increase $N$. When $N$ is sufficiently large, the advertisement patterns of both firms start to repeat themselves. This means that the profits of the two firms increase in an almost linear fashion, keeping the \% VoI almost stable.

The rest of the problem parameters (price, reservation utility, monetary value of reputation, advertisement levels, cost of advertisement and market size) affect the VoI in a similar manner to the variable capacity problem. Therefore, we do not discuss them in this section.

\input{sections/summary}

%% file: distributions/prm_distr_var.tex
\begin{table}[H]
    \centering
    \caption{Distributions for the parameter selection in the computational experiments.}
    \begin{tabular}{cc}
    \hline\hline
        Parameter & Distribution \\ \hline 
         $R_1    $ & $  U(0, 1)$\\
         $w      $ & $  U(0, 1)$\\
         $L      $ & $  U(0, 1)$\\
         $H      $ & $  U(L, 1)$\\
         $N_{variable}      $ & $  U(5, 50)$\\
         $N_{constant}      $ & $  U(5, 25)$ \\
         $ c_A   $ & $  U(1, 1000)$\\
         $\alpha $ & $  U(0, 10)$\\
         $ m     $ & $  U(5, 40)$\\ 
         $\Lambda$ & $  U(500, 2500)$\\
         $ u     $ & $ U(0, 5)$\\
         $ c     $ & $ U(0, 10)$\\
         $ p     $ & $ U(c, 25)$\\
         $S_{constant}      $ & $ U(1, \Lambda)$\\
    \hline\hline
    \end{tabular}
    \label{tab:prm_distr_var}
\end{table}

%% file: charts/var_cap_alpha_profits.tex
\definecolor{ggreen}{HTML}{6bc976}

\begin{figure}[H]
\centering
\begin{tikzpicture}
\pgfplotsset{every axis y label/.append style={yshift=-0.5cm, scale=0.75}, every x tick label/.append style={font=\tiny}, every y tick label/.append style={font=\tiny}}
\begin{axis}[
    axis y line*=left,
    xlabel = \(\alpha\),
    ylabel = {\( \Pi \)},
    height=5cm,
    width=10cm,
    x label style={at={(axis description cs:0.5,0.05)},anchor=north}
]
\addplot[ultra thick,mark=x,black]
  coordinates{
(0.52,	334935.2389 )
(0.69,	333838.6592 )
(1.16,	331259.7812 )
(1.27,	330727.8932 )
(1.39,	330191.4657 )
(1.63,	329272.0003 )
(2.06,	327977.5642 )
(2.08,	327921.8869 )
(2.30,	327434.2775 )
(2.43,	327181.6337 )
(2.52,	327033.2916 )
(2.65,	326828.1618 )
(2.89,	326508.0375 )
(3.12,	326266.0939 )
(3.24,	326155.7441 )
(3.44,	326004.2525 )
(3.67,	325865.9753 )
(3.67,	325864.1025 )
(3.68,	325857.0953 )
(4.09,	325672.8067 )
(4.12,	325663.5105 )
(4.39,	325580.2013 )
(5.00,	301832.3077 )
(5.09,	301813.2331 )
(5.10,	301809.9604 )
(5.12,	301805.8799 )
(5.20,	289971.2225 )
(5.33,	278118.6971 )
(5.44,	266268.1066 )
(5.55,	254415.3836 )
(6.11,	159624.6662 )
(6.31,	112230.3614 )
(6.59,	17469.21531 )
(6.61,	17466.44104 )
(7.15,	0           )
(7.23,	0           )
(7.32,	0           )
(7.39,	0           )
(7.53,	0           )
(7.78,	0           )
(7.94,	0           )
(8.11,	0           )
(8.47,	0           )
(9.24,   0          )
(9.25,   0          )
(9.42,   0          )
(9.71,   0          )
};

\addplot[ultra thick,mark=o,red, mark options={scale=0.5}]
  coordinates{
(0.52,	334935.2389)
(0.69,	333838.6592)
(1.16,	331259.7812)
(1.27,	330727.8932)
(1.39,	330191.4657)
(1.63,	329272.0003)
(2.06,	327983.3165)
(2.08,	327927.0147)
(2.30,	327434.2775)
(2.43,	327181.6337)
(2.52,	327033.2916)
(2.65,	326828.1618)
(2.89,	326508.0375)
(3.12,	326266.0939)
(3.24,	326155.7441)
(3.44,	326004.2525)
(3.67,	325865.9753)
(3.67,	325864.1025)
(3.68,	325857.0953)
(4.09,	325672.8067)
(4.12,	325663.5105)
(4.39,	325580.2013)
(5.00,	325456.0969)
(5.09,	325444.3265)
(5.10,	325442.2977)
(5.12,	325439.7642)
(5.20,	313644.0294)
(5.33,	313614.6792)
(5.44,	313593.4336)
(5.55,	313574.7668)
(6.11,	301646.9317)
(6.31,	289758.1011)
(6.59,	277846.7849)
(6.61,	277844.3744)
(7.15,	242128.4995)
(7.23,	230232.5285)
(7.32,	218332.2775)
(7.39,	218322.7585)
(7.53,	194525.5855)
(7.78,	158816.4050)
(7.94,	123161.2973)
(8.11,	87406.12324)
(8.47,	0	0      )
(9.24,	0	0      )
(9.25,	0	0      )
(9.42,	0	0      )
(9.71,	0	0      )
};

\end{axis}
\pgfplotsset{every axis y label/.append style={rotate=180, yshift=13.3cm}}
\begin{axis}[
    axis y line*=right,
    axis x line=none,
    ylabel = {\(\% VoI\)},
    scaled y ticks={real:0.01},
    every y tick scale label/.style={
            color=white
        },
    height=5cm,
    width=10cm,
    legend style={at={(0.3,0.5)},
	anchor=north east, nodes={scale=0.7, transform shape}},
]
\addlegendimage{smooth,mark=x,red}\addlegendentry{$\Pi_{aware}$}
\addlegendimage{smooth,mark=o,black, mark options={scale=0.5}}\addlegendentry{$\Pi_{naive}$}

\addplot[ultra thick,mark=*,ggreen, mark options={scale=0.5}]
  coordinates{
(0.52,	0          )
(0.69,	0          )
(1.16,	0          )
(1.27,	0          )
(1.39,	0          )
(1.63,	0          )
(2.06,	1.75384E-05)
(2.08,	1.5637E-05 )
(2.30,	0          )
(2.43,	0          )
(2.52,	0          )
(2.65,	0          )
(2.89,	0          )
(3.12,	0          )
(3.24,	0          )
(3.44,	0          )
(3.67,	0          )
(3.67,	0          )
(3.68,	0          )
(4.09,	0          )
(4.12,	0          )
(4.39,	0          )
(5.00,	0.072586716)
(5.09,	0.072611785)
(5.10,	0.07261606 )
(5.12,	0.072621379)
(5.20,	0.07547667 )
(5.33,	0.113183421)
(5.44,	0.150913004)
(5.55,	0.188661173)
(6.11,	0.470822842)
(6.31,	0.612675673)
(6.59,	0.937126444)
(6.61,	0.937135884)
(7.15,	1          )
(7.23,	1          )
(7.32,	1          )
(7.39,	1          )
(7.53,	1          )
(7.78,	1          )
(7.94,	1          )
(8.11,	1          )
(8.47,	0          )
(9.24,	0          )
(9.25,	0          )
(9.42,	0          )
(9.71,	0          )
(9.73,	0          )
(9.85,	0          )
(9.89,	0          )
}; \addlegendentry{\% VoI}
\end{axis}
\end{tikzpicture}
    \caption{Profits and the \%VoI with respect to advertisement resistance $\alpha$ under variable capacity.}
    \label{fig:var_cap_alpha_prft}
\end{figure}

%% file: charts/var_cap_r1_profits.tex
\definecolor{ggreen}{HTML}{6bc976}

\begin{figure}[H]
\centering
\begin{tikzpicture}
\pgfplotsset{every axis y label/.append style={yshift=-0.5cm, scale=0.75}, every x tick label/.append style={font=\tiny}, every y tick label/.append style={font=\tiny}}
\begin{axis}[
    axis y line*=left,
    xlabel = \(R_1\),
    ylabel = {\( \Pi \)},
    height=5cm,
    width=10cm,
    x label style={at={(axis description cs:0.5,0.05)},anchor=north}
]
\addplot[ultra thick,mark=x,red]
  coordinates{
(0.4098,	0          )
(0.4352,	0          )
(0.4504,	0          )
(0.4531,	13.922501  )
(0.4594,	219.139617 )
(0.4710,	1041.944411)
(0.5063,	2959.915426)
(0.5383,	4350.765321)
(0.5648,	5466.334057)
(0.5860,	6305.475572)
(0.5875,	6360.400196)
(0.5964,	6697.896047)
(0.5993,	6807.650428)
(0.6283,	7889.786917)
(0.6391,	8283.522247)
(0.6983,	10337.92513)
(0.7531,	12162.44367)
(0.7655,	12553.00607)
(0.8175,	14059.6502 )
(0.8181,	14076.20477)
(0.8354,	14539.28728)
(0.8421,	14714.68472)
(0.8462,	14820.3636 )
(0.8629,	15238.49697)
(0.8795,	15642.55073)
(0.9304,	16791.55511)
(0.9441,	17080.44925)
(0.9955,	18101.05112)
};

\addplot[ultra thick,mark=o,black, mark options={scale=0.5}]
  coordinates{
(0.4098,	0          )
(0.4352,	0          )
(0.4504,	0          )
(0.4531,	0          )
(0.4594,	0          )
(0.4710,	0          )
(0.5063,	0          )
(0.5383,	2863.851113)
(0.5648,	4358.735182)
(0.5860,	5467.926442)
(0.5875,	5540.074482)
(0.5964,	5982.191746)
(0.5993,	6125.521458)
(0.6283,	7468.926342)
(0.6391,	7943.208295)
(0.6983,	10293.19809)
(0.7531,	12162.44367)
(0.7655,	12553.00607)
(0.8175,	14059.6502 )
(0.8181,	14076.20477)
(0.8354,	14539.28728)
(0.8421,	14714.68472)
(0.8462,	14820.3636 )
(0.8629,	15238.49697)
(0.8795,	15642.55073)
(0.9304,	16791.55511)
(0.9441,	17080.44925)
(0.9955,	18101.05112)
};

\end{axis}
\pgfplotsset{every axis y label/.append style={rotate=180, yshift=13.3cm}}
\begin{axis}[
    axis y line*=right,
    axis x line=none,
    ylabel = {\(\% VoI\)},
    scaled y ticks={real:0.01},
    every y tick scale label/.style={
            color=white
        },
    height=5cm,
    width=10cm,
    legend style={at={(0.9,0.5)},
	anchor=north east, nodes={scale=0.7, transform shape}},
]
\addlegendimage{smooth,mark=x,red}\addlegendentry{$\Pi_{aware}$}
\addlegendimage{smooth,mark=o,black, mark options={scale=0.5}}\addlegendentry{$\Pi_{naive}$}

\addplot[ultra thick,mark=*,ggreen, mark options={scale=0.5}]
  coordinates{
(0.4098,	0           )
(0.4352,	0           )
(0.4504,	0           )
(0.4531,	1           )
(0.4594,	1           )
(0.4710,	1           )
(0.5063,	1           )
(0.5383,	0.34175923  )
(0.5648,	0.202621878 )
(0.5860,	0.13282886  )
(0.5875,	0.128973915 )
(0.5964,	0.106855092 )
(0.5993,	0.100200352 )
(0.6283,	0.053342451 )
(0.6391,	0.041083242 )
(0.6983,	0.0043265   )
(0.7531,	0           )
(0.7655,	0           )
(0.8175,	0           )
(0.8181,	0           )
(0.8354,	0           )
(0.8421,	0           )
(0.8462,	0           )
(0.8629,	0           )
(0.8795,	0           )
(0.9304,	0           )
(0.9441,	0           )
(0.9955,	0           )
}; \addlegendentry{\% VoI}
\end{axis}
\end{tikzpicture}
    \caption{Profits and the \% VoI with respect to the initial reputation $R_1$ under variable capacity.}
    \label{fig:var_cap_r1_prft_sample1}
\end{figure}

%% file: charts/var_cap_p_profits.tex
\definecolor{ggreen}{HTML}{6bc976}

\begin{figure}[H]
\centering
\begin{tikzpicture}
\pgfplotsset{every axis y label/.append style={yshift=-0.5cm, scale=0.75}, every x tick label/.append style={font=\tiny}, every y tick label/.append style={font=\tiny}}
\begin{axis}[
    axis y line*=left,
    xlabel = \(p\),
    ylabel = {\( \Pi \)},
    height=5cm,
    width=10cm,
    x label style={at={(axis description cs:0.5,0.05)},anchor=north}
]
\addplot[ultra thick,mark=x,red]
  coordinates{
(9.11	    ,    0)
(9.23	    ,    0)
(9.344990039,	1292.935363)
(9.386497881,	2429.478096)
(9.440295385,	3866.691795)
(10.0333725	,   17897.14718)
(10.12339138,	20001.7186 )
(10.66502263,	29914.41331)
(11.09858279,	36995.17563)
(11.18861406,	37261.75848)
(11.78443128,	43639.07399)
(11.88940889,	45122.55525)
(13.11541117,	48917.04339)
(13.90147873,	47323.42343)
(14.02112145,	47964.87598)
(14.30225707,	46761.35501)
(14.42442471,	44659.97259)
(14.63178894,	42741.36887)
(15.04019745,	38360.43631)
(15.37813373,	36110.26181)
(15.55878329,	33463.30692)
(15.59926373,	33491.82277)
(15.72304701,	30711.10452)
(15.99936027,	27873.70589)
(16.37023988,	21936.45479)
(16.54091655,	18908.59319)
(16.69206676,	15885.1058 )
(16.92483118,	12749.27012)
(16.99309642,	12645.28378)
(17.1540663	,   9622.041829)
(17.1664568	,   9601.250989)
(17.59803683,	1077.951511)
(17.69970841,	900.016031 )
(17.73211874,	0          )
(17.84123322,	0          )
(18.02423551,	0          )
(18.68942424,	0          )
%(19.26784249,	0          )
%(19.36163043,	0          )
%(19.53199252,	0          )
%(20.00216878,	0          )
%(20.46670049,	0          )
%(20.64435679,	0          )
%(20.74625511,	0          )
%(20.94800721,	0          )
%(21.79363065,	0          )
%(22.40312019,	0          )
%(22.65787362,	0          )
%(23.64332286,	0          )
%(24.22902297,	0          )
%(24.28630851,	0          )
%(24.65882216,	0          )
};

\addplot[ultra thick,mark=o,black, mark options={scale=0.5}]
  coordinates{
(9.11	    ,      0)
(9.23	    ,      0)
(9.344990039,	   0)
(9.386497881,	   699.947776)
(9.440295385,	   1565.889637)
(10.0333725	,      10398.01424)
(10.12339138,	   11080.82055)
(10.66502263,	   14517.79307)
(11.09858279,	   15171.85574)
(11.18861406,	   14557.6207)
(11.78443128,	   10514.61079)
(11.88940889,	   9257.065873)
(13.11541117,	   0)
(13.90147873,	   0)
(14.02112145,	   0)
(14.30225707,	   0)
(14.42442471,	   0)
(14.63178894,	   0)
(15.04019745,	   0)
(15.37813373,	   0)
(15.55878329,	   0)
(15.59926373,	   0)
(15.72304701,	   0)
(15.99936027,	   0)
(16.37023988,	   0)
(16.54091655,	   0)
(16.69206676,	   0)
(16.92483118,	   0)
(16.99309642,	   0)
(17.1540663	,      0)
(17.1664568	,      0)
(17.59803683,	   0)
(17.69970841,	   0)
(17.73211874,	   0)
(17.84123322,	   0)
(18.02423551,	   0)
(18.68942424,	   0)
%(19.26784249,	   0)
%(19.36163043,	   0)
%(19.53199252,	   0)
%(20.00216878,	   0)
%(20.46670049,	   0)
%(20.64435679,	   0)
%(20.74625511,	   0)
%(20.94800721,	   0)
%(21.79363065,	   0)
%(22.40312019,	   0)
%(22.65787362,	   0)
%(23.64332286,	   0)
%(24.22902297,	   0)
%(24.28630851,	   0)
%(24.65882216,	   0)
};

\end{axis}
\pgfplotsset{every axis y label/.append style={rotate=180, yshift=13.3cm}}
\begin{axis}[
    axis y line*=right,
    axis x line=none,
    ylabel = {\(\% VoI\)},
    scaled y ticks={real:0.01},
    every y tick scale label/.style={
            color=white
        },
    height=5cm,
    width=10cm,
    legend style={at={(0.9,0.72)},
	anchor=north east, nodes={scale=0.7, transform shape}},
]
\addlegendimage{smooth,mark=x,red}\addlegendentry{$\Pi_{aware}$}
\addlegendimage{smooth,mark=o,black, mark options={scale=0.5}}\addlegendentry{$\Pi_{naive}$}

\addplot[ultra thick,mark=*,ggreen, mark options={scale=0.5}]
  coordinates{
(9.11	    ,   0          )
(9.23	    ,   0          )
(9.344990039,	1          )
(9.386497881,	0.71189377 )
(9.440295385,	0.595031174)
(10.0333725	,   0.419012755)
(10.12339138,	0.446006577)
(10.66502263,	0.514689025)
(11.09858279,	0.589896372)
(11.18861406,	0.609314716)
(11.78443128,	0.759055135)
(11.88940889,	0.794846151)
(13.11541117,	1          )
(13.90147873,	1          )
(14.02112145,	1          )
(14.30225707,	1          )
(14.42442471,	1          )
(14.63178894,	1          )
(15.04019745,	1          )
(15.37813373,	1          )
(15.55878329,	1          )
(15.59926373,	1          )
(15.72304701,	1          )
(15.99936027,	1          )
(16.37023988,	1          )
(16.54091655,	1          )
(16.69206676,	1          )
(16.92483118,	1          )
(16.99309642,	1          )
(17.1540663	,   1          )
(17.1664568	,   1          )
(17.59803683,	1          )
(17.69970841,	1          )
(17.73211874,	0          )
(17.84123322,	0          )
(18.02423551,	0          )
(18.68942424,	0          )
%(19.26784249,	0          )
%(19.36163043,	0          )
%(19.53199252,	0          )
%(20.00216878,	0          )
%(20.46670049,	0          )
%(20.64435679,	0          )
%(20.74625511,	0          )
%(20.94800721,	0          )
%(21.79363065,	0          )
%(22.40312019,	0          )
%(22.65787362,	0          )
%(23.64332286,	0          )
%(24.22902297,	0          )
%(24.28630851,	0          )
%(24.65882216,	0          )
}; \addlegendentry{\% VoI}
\end{axis}
\end{tikzpicture}
    \caption{Profits and the VoI with respect to price $p$ under variable capacity.}
    \label{fig:var_cap_p_prft}
\end{figure}

%% file: charts/con_cap_fluct.tex
\definecolor{bblue}{HTML}{BACEDD}
\definecolor{ggreen}{HTML}{6bc976}

\begin{figure}[H]
    \centering
\begin{tikzpicture}
\pgfplotsset{every axis y label/.append style={yshift=-0.2cm}, every x tick label/.append style={font=\tiny}, every y tick label/.append style={font=\tiny}}
\begin{axis}[
    axis lines = left,
    xlabel = \( time \:\: period \: (i)\),
    ylabel = {\( R_i\)},
    extra y tick style={yticklabel=\pgfmathprintnumber{\tick}},
    ymin=0.0, ymax=1,
    height=5cm, width=8cm,
    legend columns=6,
    legend style={at={(1.2,1.3), font=\scriptsize, nodes={scale=0.3, transform shape}}}
]
% X_1

% R_i
\addplot [
    domain=0:20, 
    samples=100, 
    color=red,
    mark=x,
    ultra thick,
    smooth
]
    coordinates {
(1,	0.177253582)
(2,	0.422440834)
(3,	0.699027679)
(4,	0.399609882)
(5,	0.456640914)
(6,	0.602254947)
(7,	0.540389069)
(8,	0.465013003)
(9,	0.583483082)
(10,	0.587804422)
(11,	0.575913557)
(12,	0.610285935)
(13,	0.523057829)
(14,	0.486227956)
(15,	0.542140141)
(16,	0.746273741)
(17,	0.361940645)
(18,	0.995367579)
    };
% naive
\addplot [
    domain=0:20, 
    samples=100, 
    color=black,
    mark=x,
    ultra thick,
    smooth
]
    coordinates {
(1,	0.177253582)
(2,	0.255406229)
(3,	0.326135159)
(4,	0.995107625)
(5,	0.267890167)
(6,	0.994684757)
(7,	0.267978508)
(8,	0.994685398)
(9,	0.267978374)
(10,	0.994685397)
(11,	0.267978374)
(12,	0.994685397)
(13,	0.267978374)
(14,	0.994685397)
(15,	0.267978374)
(16,	0.994685397)
(17,	0.267978374)
(18,	0.994685397)
    };
\legend {aware, naive}
\end{axis}
\end{tikzpicture}
    \caption{A case of a cycling reputation under constant capacity.}
    \label{r_const_cap_fluct}
\end{figure}

%% file: charts/con_cap_s1.tex
\definecolor{ggreen}{HTML}{6bc976}

% \begin{figure}[H]
\centering
\pgfkeys{/pgf/number format/.cd,1000 sep={}}
\begin{tikzpicture}[scale=\resize]
\pgfplotsset{every axis y label/.append style={yshift=-0.5cm, scale=0.75}, every axis x label/.append style={yshift=0.5cm}, every x tick label/.append style={font=\tiny}, every y tick label/.append style={font=\tiny}}
\begin{axis}[
    axis y line*=left,
    xlabel = \(S\),
    ylabel = {\( \Pi \)},
    height=5cm,
    width=8cm,
    x label style={at={(axis description cs:0.5,-0.1)},anchor=north}
]
\addplot[ultra thick,mark=x,red]
  coordinates{
(14.00, 0               )
(20.00, 0               )
(26.00, 0               )
(35.00, 0               )
(59.00, 4423.083339     )
(81.00, 11129.81981     )
(88.00, 11804.74606     )
(93.00, 12975.0781      )
(100.0, 	14693.2629  )
(142.0, 	25084.1302  )
(228.0, 	48552.28216 )
(307.0, 	68253.01351 )
(314.0, 	69971.19831 )
(316.0, 	70462.10826 )
(318.0, 	70953.0182  )
(324.0, 	72425.74803 )
(335.0, 	75125.75272 )
(336.0, 	75371.20769 )
(380.0, 	75693.94503 )
(400.0, 	75440.66582 )
(420.0, 	75187.38662 )
(426.0, 	75111.40285 )
(449.0, 	74820.13177 )
(454.0, 	74756.81196 )
(472.0, 	74528.86068 )
(474.0, 	74503.53276 )
(480.0, 	74427.54899 )
(488.0, 	74326.23731 )
(494.0, 	74250.25355 )
(495.0, 	74237.58959 )
(509.0, 	74060.29415 )
(541.0, 	73655.04741 )
(543.0, 	73629.71949 )
(545.0, 	73604.39157 )
(556.0, 	73465.08801 )
(585.0, 	73097.83316 )
(588.0, 	73059.84128 )
(600.0, 	72907.87375 )
(627.0, 	72565.94683 )
(639.0, 	72413.9793  )
(652.0, 	72323.43129 )
};

\addplot[ultra thick,mark=o,black, mark options={scale=0.5}]
  coordinates{
(14.00, 0           )
(20.00, 0           )
(26.00, 0           )
(35.00, 0           )
(59.00, 0           )
(81.00, 1116.278306 )
(88.00, 2340.979145 )
(93.00, 4428.817343 )
(100.0, 7249.759043 )
(142.0, 22141.20423 )
(228.0, 43885.97748 )
(307.0, 62312.24661 )
(314.0, 63917.50438 )
(316.0, 64376.14945 )
(318.0, 64834.79453 )
(324.0, 66210.72976 )
(335.0, 68733.27769 )
(336.0, 68962.60022 )
(380.0, 69225.86791 )
(400.0, 68972.5887  )
(420.0, 68719.3095  )
(426.0, 68643.32574 )
(449.0, 68352.05465 )
(454.0, 68288.73485 )
(472.0, 68060.78356 )
(474.0, 68035.45564 )
(480.0, 67959.47188 )
(488.0, 67858.1602  )
(494.0, 67782.17643 )
(495.0, 67769.51247 )
(509.0, 67592.21703 )
(541.0, 67186.9703  )
(543.0, 67161.64238 )
(545.0, 67136.31446 )
(556.0, 66997.01089 )
(585.0, 66629.75604 )
(588.0, 66591.76416 )
(600.0, 66439.79664 )
(627.0, 66097.86971 )
(639.0, 65945.90218 )
(652.0, 65893.86736 )
};

\end{axis}
\pgfplotsset{every axis y label/.append style={rotate=180, yshift=10.5cm}}
\begin{axis}[
    axis y line*=right,
    axis x line=none,
    ylabel = {\(\% VoI\)},
    scaled y ticks={real:0.01},
    every y tick scale label/.style={
            color=white
        },
    height=5cm,
    width=8cm,
    legend style={at={(0.9,0.72)},
	anchor=north east, nodes={scale=0.7, transform shape}},
]
\addlegendimage{smooth,mark=x,red}\addlegendentry{$\Pi_{aware}$}
\addlegendimage{smooth,mark=o,black, mark options={scale=0.5}}\addlegendentry{$\Pi_{naive}$}

\addplot[ultra thick,mark=*,ggreen, mark options={scale=0.5}]
  coordinates{
(14.00, 0          )
(20.00, 0          )
(26.00, 0          )
(35.00, 0          )
(59.00, 1          )
(81.00, 0.899703829)
(88.00, 0.801691698)
(93.00, 0.658667385)
(100.0, 0.506592981)
(142.0, 0.117322225)
(228.0, 0.096108864)
(307.0, 0.087040361)
(314.0, 0.08651694 )
(316.0, 0.086372079)
(318.0, 0.086229224)
(324.0, 0.085812276)
(335.0, 0.085090329)
(336.0, 0.085027263)
(380.0, 0.08545039 )
(400.0, 0.085737275)
(420.0, 0.086026093)
(426.0, 0.086113118)
(449.0, 0.086448352)
(454.0, 0.086521575)
(472.0, 0.086786207)
(474.0, 0.08681571 )
(480.0, 0.086904341)
(488.0, 0.087022798)
(494.0, 0.087111852)
(495.0, 0.087126712)
(509.0, 0.087335288)
(541.0, 0.087815803)
(543.0, 0.087846011)
(545.0, 0.087876239)
(556.0, 0.088042869)
(585.0, 0.08848521 )
(588.0, 0.088531223)
(600.0, 0.088715756)
(627.0, 0.08913378 )
(639.0, 0.089320835)
(652.0, 0.088900151)
}; \addlegendentry{\% VoI}
\end{axis}
\end{tikzpicture}
% \end{figure}

%% file: charts/con_cap_s2.tex
\definecolor{ggreen}{HTML}{6bc976}

% \begin{figure}[H]
\centering
\pgfkeys{/pgf/number format/.cd,1000 sep={}}
\begin{tikzpicture}[scale=\resize]
\pgfplotsset{every axis y label/.append style={yshift=-0.5cm, scale=0.75}, every axis x label/.append style={yshift=0.5cm}, every x tick label/.append style={font=\tiny}, every y tick label/.append style={font=\tiny}}
\begin{axis}[
    axis y line*=left,
    xlabel = \(S\),
    ylabel = {\( \Pi \)},
    height=5cm,
    width=8cm,
    x label style={at={(axis description cs:0.5,-0.1)},anchor=north}
]
\addplot[ultra thick,mark=x,red]
  coordinates{
(511.0000,	21418.67803)
(647.0000,	26733.06259)
(797.0000,	31967.81806)
(868.0000,	33892.18443)
(976.0000,	36452.02477)
(1059.000,	37720.81746)
(1072.000,	37808.95328)
(1121.000,	37903.39662)
(1525.000,	29393.04324)
(1548.000,	28874.53221)
(1570.000,	28378.56513)
(1651.000,	26552.50454)
(1676.000,	25988.90559)
(1736.000,	24636.26812)
%(1817.000,	22810.20753)
%(1896.000,	21029.23486)
%(1922.000,	20443.09195)
%(1928.000,	20307.8282 )
%(2004.000,	18594.4874 )
%(2107.000,	16272.45974)
};

\addplot[ultra thick,mark=o,black, mark options={scale=0.5}]
  coordinates{
(511.0000,  5664.056562)
(647.0000,  6503.577022)
(797.0000,  7052.533954)
(868.0000,  6981.314919)
(976.0000,  6225.394322)
(1059.000,  4886.818903)
(1072.000,  4593.747451)
(1121.000,  3489.093514)
(1525.000,  0          )
(1548.000,  0          )
(1570.000,  0          )
(1651.000,  0          )
(1676.000,  0          )
(1736.000,  0          )
%(1817.000,  0          )
%(1896.000,  0          )
%(1922.000,  0          )
%(1928.000,  0          )
%(2004.000,  0          )
%(2107.000,  0          )
};

\end{axis}
\pgfplotsset{every axis y label/.append style={rotate=180, yshift=10.5cm}}
\begin{axis}[
    axis y line*=right,
    axis x line=none,
    ylabel = {\(\% VoI\)},
    scaled y ticks={real:0.01},
    every y tick scale label/.style={
            color=white
        },
    height=5cm,
    width=8cm,
    legend style={at={(0.9,0.6)},
	anchor=north east, nodes={scale=0.7, transform shape}},
]
\addlegendimage{smooth,mark=x,red}\addlegendentry{$\Pi_{aware}$}
\addlegendimage{smooth,mark=o,black, mark options={scale=0.5}}\addlegendentry{$\Pi_{naive}$}

\addplot[ultra thick,mark=*,ggreen, mark options={scale=0.5}]
  coordinates{
(511.0000,	0.735555268)
(647.0000,	0.756721588)
(797.0000,	0.779386446)
(868.0000,	0.794014017)
(976.0000,	0.829216776)
(1059.000,	0.870447694)
(1072.000,	0.878501068)
(1121.000,	0.907947735)
(1525.000,	1          )
(1548.000,	1          )
(1570.000,	1          )
(1651.000,	1          )
(1676.000,	1          )
(1736.000,	1          )
%(1817.000,	1          )
%(1896.000,	1          )
%(1922.000,	1          )
%(1928.000,	1          )
%(2004.000,	1          )
%(2107.000,	1          )
}; \addlegendentry{\% VoI}
\end{axis}
\end{tikzpicture}
% \end{figure}

%% file: charts/con_cap_s3.tex
\definecolor{ggreen}{HTML}{6bc976}

% \begin{figure}[H]
\centering
\pgfkeys{/pgf/number format/.cd,1000 sep={}}
\begin{tikzpicture}[scale=\resize]
\pgfplotsset{every axis y label/.append style={yshift=-0.5cm, scale=0.75}, every axis x label/.append style={yshift=0.5cm}, every x tick label/.append style={font=\tiny}, every y tick label/.append style={font=\tiny}}
\begin{axis}[
    axis y line*=left,
    xlabel = \(S\),
    ylabel = {\( \Pi \)},
    height=5cm,
    width=8cm,
    x label style={at={(axis description cs:0.5,-0.1)},anchor=north}
]
\addplot[ultra thick,mark=x,red]
  coordinates{
(571.0000,	38154.81582)
(596.0000,	39825.30722)
(691.0000,	46173.17455)
(699.0000,	46707.73179)
(766.0000,	50972.67615)
(804.0000,	52938.83614)
(841.0000,	54757.83442)
(1024.000,	25911.61952)
(1048.000,	21725.10753)
(1118.000,	9514.447565)
(1232.000,	0          )
(1266.000,	0          )
(1275.000,	0          )
(1307.000,	0          )
(1341.000,	0          )
(1346.000,	0          )
(1484.000,	0          )
(1571.000,	0          )
(1672.000,	0          )
(1678.000,	0          )
};

\addplot[ultra thick,mark=o,black, mark options={scale=0.5}]
  coordinates{
(571.0000,	0)
(596.0000,	0)
(691.0000,	0)
(699.0000,	0)
(766.0000,	0)
(804.0000,	0)
(841.0000,	0)
(1024.000,	0)
(1048.000,	0)
(1118.000,	0)
(1232.000,	0)
(1266.000,	0)
(1275.000,	0)
(1307.000,	0)
(1341.000,	0)
(1346.000,	0)
(1484.000,	0)
(1571.000,	0)
(1672.000,	0)
(1678.000,	0)
};

\end{axis}
\pgfplotsset{every axis y label/.append style={rotate=180, yshift=10.5cm}}
\begin{axis}[
    axis y line*=right,
    axis x line=none,
    ylabel = {\(\% VoI\)},
    scaled y ticks={real:0.01},
    every y tick scale label/.style={
            color=white
        },
    height=5cm,
    width=8cm,
    legend style={at={(0.9,0.72)},
	anchor=north east, nodes={scale=0.7, transform shape}},
]
\addlegendimage{smooth,mark=x,red}\addlegendentry{$\Pi_{aware}$}
\addlegendimage{smooth,mark=o,black, mark options={scale=0.5}}\addlegendentry{$\Pi_{naive}$}

\addplot[ultra thick,mark=*,ggreen, mark options={scale=0.5}]
  coordinates{
(571.0000,	1)
(596.0000,	1)
(691.0000,	1)
(699.0000,	1)
(766.0000,	1)
(804.0000,	1)
(841.0000,	1)
(1024.000,	1)
(1048.000,	1)
(1118.000,	1)
(1232.000,	0)
(1266.000,	0)
(1275.000,	0)
(1307.000,	0)
(1341.000,	0)
(1346.000,	0)
(1484.000,	0)
(1571.000,	0)
(1672.000,	0)
(1678.000,	0)
}; \addlegendentry{\% VoI}
\end{axis}
\end{tikzpicture}
% \end{figure}

%% file: charts/con_cap_w.tex
\definecolor{ggreen}{HTML}{6bc976}

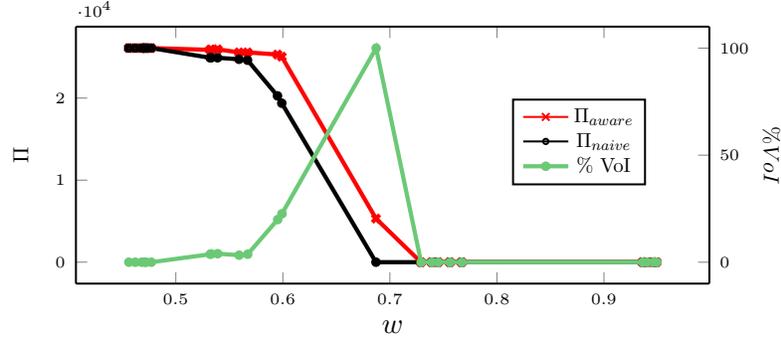
\begin{figure}[H]
\centering
\begin{tikzpicture}
\pgfplotsset{every axis y label/.append style={yshift=-0.5cm, scale=0.75}, every x tick label/.append style={font=\tiny}, every y tick label/.append style={font=\tiny}}
\begin{axis}[
    axis y line*=left,
    xlabel = \(w\),
    ylabel = {\( \Pi \)},
    height=5cm,
    width=10cm,
    x label style={at={(axis description cs:0.5,0.05)},anchor=north}
]
\addplot[ultra thick,mark=x,red]
  coordinates{
%(0.0530,	26098.43189)
%(0.0860,	26098.43189)
%(0.0910,	26098.43189)
%(0.1180,	26098.43189)
%(0.1480,	26098.43189)
%(0.1850,	26098.43189)
%(0.1920,	26098.43189)
%(0.1940,	26098.43189)
%(0.2310,	26098.43189)
%(0.2430,	26098.43189)
%(0.2540,	26098.43189)
%(0.2720,	26098.43189)
%(0.2750,	26098.43189)
%(0.2760,	26098.43189)
%(0.2930,	26098.43189)
%(0.3110,	26098.43189)
%(0.3190,	26098.43189)
(0.4560,	26098.43189)
(0.4620,	26098.43189)
(0.4680,	26098.43189)
(0.4710,	26098.43189)
(0.4720,	26098.43189)
(0.4770,	26098.43189)
(0.5320,	25879.84614)
(0.5330,	25922.7096 )
(0.5390,	25949.70276)
(0.5590,	25573.0349 )
(0.5670,	25573.0349 )
(0.5950,	25310.33641)
(0.5990,	25047.63792)
(0.6870,	5316.645551)
(0.7290,	0          )
(0.7390,	0          )
(0.7450,	0          )
(0.7560,	0          )
(0.7670,	0          )
(0.9360,	0          )
(0.9400,	0          )
(0.9460,	0          )
(0.9490,	0          )
};

\addplot[ultra thick,mark=o,black, mark options={scale=0.5}]
  coordinates{
%(0.0530,	26098.43189)
%(0.0860,	26098.43189)
%(0.0910,	26098.43189)
%(0.1180,	26098.43189)
%(0.1480,	26098.43189)
%(0.1850,	26098.43189)
%(0.1920,	26098.43189)
%(0.1940,	26098.43189)
%(0.2310,	26098.43189)
%(0.2430,	26098.43189)
%(0.2540,	26098.43189)
%(0.2720,	26098.43189)
%(0.2750,	26098.43189)
%(0.2760,	26098.43189)
%(0.2930,	26098.43189)
%(0.3110,	26098.43189)
%(0.3190,	26098.43189)
(0.4560,	26098.43189)
(0.4620,	26098.43189)
(0.4680,	26098.43189)
(0.4710,	26098.43189)
(0.4720,	26098.43189)
(0.4770,	26098.43189)
(0.5320,	24916.28867)
(0.5330,	24916.28867)
(0.5390,	24916.28867)
(0.5590,	24734.92471)
(0.5670,	24617.29114)
(0.5950,	20272.00989)
(0.5990,	19380.31598)
(0.6870,	0          )
(0.7290,	0          )
(0.7390,	0          )
(0.7450,	0          )
(0.7560,	0          )
(0.7670,	0          )
(0.9360,	0          )
(0.9400,	0          )
(0.9460,	0          )
(0.9490,	0          )
};

\end{axis}
\pgfplotsset{every axis y label/.append style={rotate=180, yshift=13.3cm}}
\begin{axis}[
    axis y line*=right,
    axis x line=none,
    ylabel = {\(\% VoI\)},
    scaled y ticks={real:0.01},
    every y tick scale label/.style={
            color=white
        },
    height=5cm,
    width=10cm,
    legend style={at={(0.9,0.72)},
	anchor=north east, nodes={scale=0.7, transform shape}},
]
\addlegendimage{smooth,mark=x,red}\addlegendentry{$\Pi_{aware}$}
\addlegendimage{smooth,mark=o,black, mark options={scale=0.5}}\addlegendentry{$\Pi_{naive}$}

\addplot[ultra thick,mark=*,ggreen, mark options={scale=0.5}]
  coordinates{
%(0.0530,	0          )
%(0.0860,	0          )
%(0.0910,	0          )
%(0.1180,	0          )
%(0.1480,	0          )
%(0.1850,	0          )
%(0.1920,	0          )
%(0.1940,	0          )
%(0.2310,	0          )
%(0.2430,	0          )
%(0.2540,	0          )
%(0.2720,	0          )
%(0.2750,	0          )
%(0.2760,	0          )
%(0.2930,	0          )
%(0.3110,	0          )
%(0.3190,	0          )
(0.4560,	0          )
(0.4620,	0          )
(0.4680,	0          )
(0.4710,	0          )
(0.4720,	0          )
(0.4770,	0          )
(0.5320,	0.037231963)
(0.5330,	0.038823909)
(0.5390,	0.039823735)
(0.5590,	0.0327732  )
(0.5670,	0.037373107)
(0.5950,	0.199062013)
(0.5990,	0.226261732)
(0.6870,	1          )
(0.7290,	0          )
(0.7390,	0          )
(0.7450,	0          )
(0.7560,	0          )
(0.7670,	0          )
(0.9360,	0          )
(0.9400,	0          )
(0.9460,	0          )
(0.9490,	0          )
}; \addlegendentry{\% VoI}
\end{axis}
\end{tikzpicture}
    \caption{Profits and the VoI with respect to the weight of WoM communication $w$ under constant capacity.}
    \label{fig:con_cap_w_prft}
\end{figure}

%% file: charts/con_cap_alpha.tex
\definecolor{ggreen}{HTML}{6bc976}

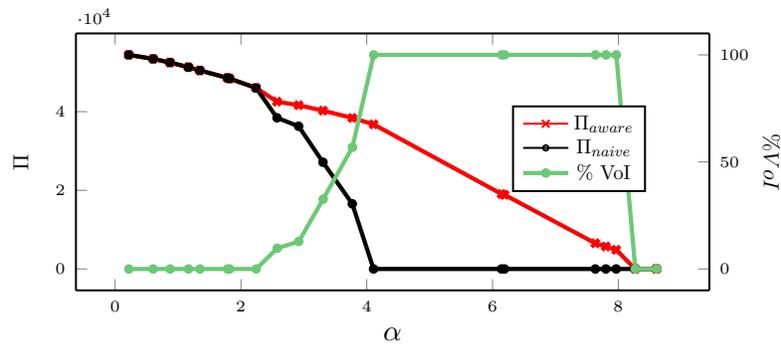
\begin{figure}[H]
\centering
\begin{tikzpicture}
\pgfplotsset{every axis y label/.append style={yshift=-0.5cm, scale=0.75}, every x tick label/.append style={font=\tiny}, every y tick label/.append style={font=\tiny}}
\begin{axis}[
    axis y line*=left,
    xlabel = \(\alpha\),
    ylabel = {\( \Pi \)},
    height=5cm,
    width=10cm,
    x label style={at={(axis description cs:0.5,0.05)},anchor=north}
]
\addplot[ultra thick,mark=x,red]
  coordinates{
(0.2160,	54444.2874 )
(0.6048,	53423.60149)
(0.8725,	52491.6675 )
(1.1644,	51314.10417)
(1.3463,	50481.01287)
(1.7848,	48595.75786)
(1.8188,	48414.34387)
(2.2401,	46034.12786)
(2.5733,	42548.67126)
(2.9168,	41629.68187)
(3.3044,	40262.28502)
(3.7701,	38404.08191)
(4.1084,	36772.027  )
(6.1425,	19060.25453)
(6.1905,	18856.20825)
(7.6334,	6531.502948)
(7.8006,	5706.782652)
(7.9651,	4865.16953 )
(8.2667,	0          )
(8.6045,	0          )
};

\addplot[ultra thick,mark=o,black, mark options={scale=0.5}]
  coordinates{
(0.2160,	54444.2874 )
(0.6048,	53423.60149)
(0.8725,	52491.6675 )
(1.1644,	51314.10417)
(1.3463,	50481.01287)
(1.7848,	48595.75786)
(1.8188,	48414.34387)
(2.2401,	46034.12786)
(2.5733,	38435.95848)
(2.9168,	36282.07438)
(3.3044,	27131.34948)
(3.7701,	16566.53655)
(4.1084,	0          )
(6.1425,	0          )
(6.1905,	0          )
(7.6334,	0          )
(7.8006,	0          )
(7.9651,	0          )
(8.2667,	0          )
(8.6045,	0          )
};

\end{axis}
\pgfplotsset{every axis y label/.append style={rotate=180, yshift=13.3cm}}
\begin{axis}[
    axis y line*=right,
    axis x line=none,
    ylabel = {\(\% VoI\)},
    scaled y ticks={real:0.01},
    every y tick scale label/.style={
            color=white
        },
    height=5cm,
    width=10cm,
    legend style={at={(0.9,0.72)},
	anchor=north east, nodes={scale=0.7, transform shape}},
]
\addlegendimage{smooth,mark=x,red}\addlegendentry{$\Pi_{aware}$}
\addlegendimage{smooth,mark=o,black, mark options={scale=0.5}}\addlegendentry{$\Pi_{naive}$}

\addplot[ultra thick,mark=*,ggreen, mark options={scale=0.5}]
  coordinates{
(0.2160,	0            )
(0.6048,	0            )
(0.8725,	0            )
(1.1644,	0            )
(1.3463,	0            )
(1.7848,	0            )
(1.8188,	0            )
(2.2401,	0            )
(2.5733,	0.096659018  )
(2.9168,	0.128456602  )
(3.3044,	0.326134881  )
(3.7701,	0.568625632  )
(4.1084,	1            )
(6.1425,	1            )
(6.1905,	1            )
(7.6334,	1            )
(7.8006,	1            )
(7.9651,	1            )
(8.2667,	0            )
(8.6045,	0            )
}; \addlegendentry{\% VoI}
\end{axis}
\end{tikzpicture}
    \caption{Profits and the VoI with respect to advertisement resistance $\alpha$ under constant capacity.}
    \label{fig:con_cap_alpha_prft}
\end{figure}

%% file: sections/summary.tex
\subsection{Managerial Insights}
\label{sec:remarks}
In both variable and constant capacity problems, our analysis reveal that the naive firms fail to identify the potential of the market. Even though they can achieve similar profits to the aware firm under some special cases, mostly they end up with realizing less profit than they could if they had invested in discovering the underlying WoM communication. Indeed, they frequently misjudge the market conditions and end up not operating while their aware equivalents seize profits. This is due to the \emph{myopic behavior} of the naive firm.

One of the key decisions in this problem is the advertisement strategy. Under endogenous demand structures, it is vital for companies to avoid attracting demand that may overwhelm the system. Otherwise, the company will have large number of unsatisfied customers, hurting the reputation and demand for the future through WoM communication. In variable capacity problem, fill rate is steady. Therefore, the role of the advertisement is mainly to pull the overall reputation to a desired level. We have shown it to have a simple structure: aware firm starts with high level advertisements and switches to low level once, the naive firm, on the other hand, can start with low level advertisements and have at most two switching points. We have seen that the VoI is low when the naive firm also starts with the high level advertisements since this is similar to the advertisement decision under full information. Therefore, a better strategy for a firm that has a flexible service capacity is to start with a high advertisement/promotion strategy when customers are involved in WoM communication even if the firm does not understand the underlying dynamics of this communication.

%Under variable capacity, the optimal advertisement decision of the aware firm can be characterized by the location of the switching point while the naive firm's decision can be described by two switching points. The gap between the two firms' profits mostly comes from the early periods, when reputation is low. If the naive firm switches from $L$ to $H$ late, this gap widens and we observe high VoI. This behavior (late switching) occurs when a high reputation level is required for picking $H$ to be immediately-profitable. On the contrary, if the naive firm's optimal solution contains single switching point (from $H$ to $L$), the naive firm usually ends up approximating aware firm's solution well. 

%On the other hand, constant capacity scenario presents a more challenging problem, often requiring elaborate advertisement policies. --> general advertisement strategy

For the constant capacity problem, it is noted that the VoI highly depends on how variable firm's demand realization ability is. When the cardinality of $\Omega$ is high, the firm may be at many different states throughout the planning horizon, which increases the complexity of the decision. Otherwise, the problem is easier as shown by Propositions \ref{Prop:ConstCapCEF}, \ref{Prop:AbsorbingA} and Lemma \ref{Lemma:ConstCapABD}. It is observed that VoI is the highest when the problem is complicated (e.g., high cardinality $\Omega$) as there is more room for the WOM to take place.

As a result of the nature of the constant capacity problem, both positive and negative WoM can be observed indicating that the efforts placed in advertisement may immediately get overridden. As a result, reputation fluctuates throughout the planning horizon, unlike the variable capacity problem. This requires the firm to make more rigorous planning to maintain its reputation and maybe even make certain maneuvers. For instance, if the firm has comprehensive understanding of the WoM, it can act anticipating a dose of negative WoM to keep demand at desired levels or keep a low profile when needed and advertise more aggressively at other times. A strong observation on this matter shows that not incorporating WoM in the decision making can yield higher promotion expenditures and yet, less revenue. Note that in the variable capacity problem, the aware firm beats the naive firm by realizing more revenue. However, under constant capacity, the aware firm may collect more revenue, have less expenditure or do both at the same time compared with its naive equivalent. On top of that, the naive firm may end up reducing reputation beyond repair after a few bad decisions, reducing the size of its own playground for the future periods.

We observe that starting reputation ($R_1$) decreases VoI under variable capacity while it is mostly not significant under constant capacity. This results from the fact that reputation converges under variable capacity while fluctuating under constant capacity. For the variable capacity problem, we may say that having information on the endogenous demand is more important for the firms that are just entering the market or the ones that have recently suffered a brand damage. Meanwhile, under constant capacity, the status of the firm is not stable at all; therefore, information is valuable regardless of the current situation of the firm.

Another observation is regarding the weight of WoM. It is seen that under variable capacity, $w$ only speeds up the convergence of the reputation curves, and does not impact VoI much. Meanwhile under constant capacity, it determines how fast firm's past actions may get overridden and is highly significant in terms of VoI.

Generally, in both variable and constant capacity problems, VoI is the lowest when the market conditions are trivial (e.g., very easy or difficult to build up and maintain reputation), and it is the highest when the market is challenging enough for the naive firm to overlook the potential profits but profitable for the aware firm. A challenging market typically involves consumers that are not easily convinced (high enough $\alpha$ and low $m$), notable substitutes for the service leading to high reservation utility ($u$) and moderate-to-high price ($p$).

%% file: sections/conclusion.tex
\section{Conclusion}
\label{sec:5}
In this paper, we consider a service system where customers share information on service availability based on their experience through Word-of-Mouth (WoM), and the service provider's reputation on the availability affects the demand. Therefore, the demand is endogenous to the service provider's operational decisions. The company can make advertisements to alter the customer perception. We study optimal advertisement strategy for the firm to maximize profits while managing the demand on a finite horizon. We consider two capacity scenarios; flexible (or on-demand) capacity and constant capacity. A restaurant using Uber Eats for food delivery is an example for the former, and a restaurant using their own staff for delivery is an example for the latter. We analyze this problem from perspectives of two types of firms based on the available information on the demand structure, namely \emph{aware} and \emph{naive} firms. The aware firm has full information on the endogenous demand structure, while the naive firm has partial information. We partially characterize the optimal advertisement and capacity decisions for both firms and analyze the value of information on the endogenous demand through extensive computational studies.

%In this paper we consider a firm operating in a service environment, where the fill rate of the service is of essence in the eyes of the consumer population. The consumer population is composed of members that have prior opinions on the fill rate of the firm, exchange information through word-of-mouth (WoM) regarding their experience and update their opinions. The ones who find it unlikely to receive the service upon request resort to substitutes in the future. The mentioned opinions and perceptions customers have regarding the fill rate constitute firm's \emph{reputation} and is the main determinant of demand. Therefore, WoM communication is a factor dynamically impacting demand through reputation, which implies that the service capacity of the firm impacts the future demand (endogenous demand). In addition, the firm may place high or low level advertisement efforts to directly improve reputation.

%We study this system in a multi-period planning horizon, where the firm makes advertisement and capacity choices in each period, which keep altering demand. We model and analyze the variable and constant capacity problems in this setting, and compare two types of firms in each problem: (i) aware firm, who has complete information on the endogenous demand structure, (ii) naive firm, who has limited information on the endogenous demand structure.

We have shown that in flexible capacity systems, the service provider meets all demand under certain cost structures to avoid negative WoM. WoM communication motivates firms to provide high quality service. The aware firm follows a more aggressive advertisement strategy in early periods and then switches to moderate advertising. The naive firm fails to foresee the benefits of current advertisement and acts myopically. We see that the aware firm can operate in many conditions (determined by problem parameters) successfully when the naive firm stops operating. In constant capacity systems, it may not be possible to avoid negative WoM. In those cases, WoM has a strong role controlling the operational decisions. The advertisement strategy becomes more compatible with the actual service quality. Therefore, the WoM leads to a more truth-telling behavior. Ultimately, the lack of information on WoM can result in loss of opportunities and unnecessary investment on advertisements especially when the service capacity is not flexible.

%% file: sections/appendices.tex
\section{Constant Capacity Models}
\label{model:const_cap}
\textbf{Aware Firm Under Constant Capacity:}

The optimality equation can be expressed as:
\begin{equation}
  P_i(R_i) =
    \begin{dcases}
      \max_{\substack{A_i\in \{H,\:L\}}} \{\pi_i + P_{i+1}(R_{i+1})\} & \text{if } i<N,\\
      \max_{\substack{A_i\in \{H,\:L\}}} \{\pi_i\} & \text{if } i=N,
    \end{dcases}       
\end{equation}
\noindent where
\begin{align*}
        & && \pi_i = p \: \min \{D_i,\:S\} - c S - c_A A_i,\\
        & && R_{i+1}=
        \begin{dcases}
        (1-w) \left(R_{i} + (1-R_{i}) A_{i}^\alpha \right) + w \frac{\min\{D_i,\:S\}}{D_i} & \text{if} \:\: D_i>0,\\
        R_{i} + (1-R_{i}) A_{i}^\alpha & \text{if} \:\: D_i=0,
        \end{dcases}\\
        & && D_i = \max\left\{0,\: \Lambda\left(1 - \frac{u+p}{m \left(R_i + (1-R_i) A_i^\alpha \right)} \right)\right\}.
\end{align*}    

\textbf{Naive Firm Under Constant Capacity:}

The firm's problem in any period $i$ can be expressed as:
\begin{equation}
    \max_{\substack{A_i\in \{H,\:L\}}} \{\pi_i\}
\end{equation}
\noindent where
\begin{align*}
        & && \pi_i = p \: \min \{D_i,\:S\} - c S - c_A A_i,\\
        & && D_i = \max\left\{0,\: \Lambda\left(1 - \frac{u+p}{m \left(R_i + (1-R_i) A_i^\alpha \right)} \right)\right\}.
\end{align*}

\section{Proof of Proposition \ref{Prop:SiDi}}
\label{Proof:Prop:SiDi}
    We want to show that the optimal capacity decision $S^*_i$ equals the realized demand under the optimal policy, $D^*_i$ for all $i$.
    
    We will first show that this strategy maximizes $\pi_i$ to cover the naive firm's policy. Then, we will provide the rest of the proof to cover the aware firm as well.

\textit{Naive Firm:}

    If we have $S_i = D_i$, then we obtain the following profit in period $i$:
    \begin{equation}
    \label{Eq:pro1}
         \pi_i=(p-c) \: D_i - c_A A_i.
    \end{equation}
    
    Otherwise, if we choose a capacity $S_i$ for period $i$ such that $S_i > D_i$, we obtain the following profit in period $i$:
    \begin{equation}
    \label{Eq:pro2}
         \pi_i=p \: D_i - c \: S_i - c_A A_i.
    \end{equation}

 The profit obtained in Equation \ref{Eq:pro1} is greater than of Equation \ref{Eq:pro2} by $c \: (S_i - D_i) > 0$. Therefore, it is sub-optimal for the naive firm to pick a capacity level greater than demand.
 
 Then, as $S_i^* \leq D^*_i$, we can rewrite the profit function as follows:
\begin{equation}
    \begin{split}
    & \pi_i = (p-c) \: S_i - c_A A_i \\
    & s. t. \\
    & 0 \leq S_i \leq D_i
    \end{split}
\label{variable s leq demand}
\end{equation}
Choosing a capacity level $S_i$ that is strictly less than $D_i$ results in loss of profit opportunity by $(p-c)(D_i-S_i) > 0$ in the current period, and therefore is sub-optimal. This concludes the proof for the naive firm.

\textit{Aware Firm:}

 We have shown that having $S_i=D_i$ maximizes $\pi_i$. However, the aware firm is interested in maximizing the profit of the entire planning horizon, $\Pi$.
 
 Suppose the aware firm has $S_i > D_i$ for some $i$. This results in the same fill rate $F_i$, and leads to the same $R_{i+1}$ as picking the capacity level equal to demand. Therefore, $S_i > D_i$ is sub-optimal.
 
 Now, suppose the aware firm has $S_i < D_i$ for some $i$. This yields less profit in period $i$ than picking the capacity level equal to demand, and results in a lower fill rate $F_i$ leading to a lower $R_{i+1}$. As a result, the firm realizes less demand in period $i+1$, which yields lower $\pi_{i+1}$. Therefore, for any set of advertisement decisions, having $S_i<D_i$ for some $i$ generates a smaller profit than having $S_i=D_i$ for all $i$, which concludes the second part of the proof.
 
\section{Proof of Proposition \ref{Prop:SwitchingPoint}}
\label{Proof:Prop:SwitchingPoint}
To prove the proposition, we will show that $A_{i+1} = H$ cannot be optimal when $A_i = L$.

Let us consider the following advertisement sequence $V$ for an N-period problem.
\begin{equation}
\begin{split}
     & V = \{V_1, \: V_2, \: ...\: V_N \}, \text{ where each } V_i \text{ is either } H \text{ or } L, \\
     & \text{and } V_k = L, \; V_{k+1} = H, \text{ for some } 1 \leq k \leq N - 1.
 \end{split}
\end{equation}
Then, using $V$, let us define a second sequence $Q$ such that:
\begin{equation}
\begin{split}
    & Q_i = V_i, \text{ for all } i \; {\in} \; \{1, \: 2, \: ..., \; k-1, \; k+2, \; ... \: N \}, \text{ and}  \\
    & Q_k = H, \; Q_{k+1} = L.
\end{split}
\end{equation}

We denote the parameters such as reputation, post-advertisement reputation, service capacity, fill rate, demand and profit at period i, using the sequence information (e.g. firm's reputation at sequence $V$ and period $i$ is represented by $R_i^V$).

By Proposition \ref{Prop:SiDi}, the variable capacity profit function is $\pi_i = (p-c) \: D_i - c_A \: A_i$. Since, in the above sequences, all decisions are common before the $k^{th}$ period, we have the following.

\begin{equation}
    \pi^Q_{1} + \pi^Q_{2} + ... \pi^Q_{k-1} = \pi^V_{1} + \pi^V_{2} + ... \pi^V_{k-1}
\end{equation}

The reputation level at the beginning of period $k$ is the same for both sequences, that is $R^V_k = R^Q_k$. The sequence $Q$ has high level advertisement at period $k$ while sequence $V$ has low. Therefore, we have $R^{'Q}_k > R^{'V}_k$, indicating $D^{Q}_k \geq D^{V}_k$.

Then, we have either one of the following cases regarding period $k+1$:
\begin{itemize}
\item[{(i)}] If $D^{Q}_k = 0$, we have $D^{V}_k=0$ as well. As a result, new information on fill rate is not created, and WoM does not happen. Then, we have the following expression for post-advertisement reputation of period $k+1$, which is symmetric around $H$ and $L$ allowing $R^{'Q}_{k+1}$ and $R^{'V}_{k+1}$ to be equal.
\begin{equation}
\begin{split}
R^{'Q}_{k+1}=R^{'V}_{k+1} &= R^V_k + H - HR^V_k  -LR^V_k - HL + HLR^V_k + L \\
    &= H (1 -R^V_k) + L (1 -R^V_k) - HL(1 -R^V_k) + R_k
\end{split}
\end{equation}
In that case, we have $D^{Q}_k + D^{Q}_{k+1} = D^{V}_k + D^{V}_{k+1}$.
\item[{(ii)}] If $D^{Q}_k > 0$ while $D^{V}_k=0$, we have the following expressions for post-advertisement reputation of period $k+1$; thus, we have $R^{'Q}_{k+1} > R^{'V}_{k+1}$.
\begin{equation}
\begin{split}
    & R^{'Q}_{k+1} = R^V_k (1-w) + w + H(1 - R^V_k - w + wR^V_k) \\
    & + L(1 - R^V_k - w + wR^V_k) + HL (-1 + R^V_k + w - wR^V_k)\\ \\
    & R^{'V}_{k+1} = H (1 -R^V_k) + L (1 -R^V_k) - HL(1 -R^V_k) + R^V_k
\end{split}
\end{equation}
In that case, we have $D^{Q}_k + D^{Q}_{k+1} > D^{V}_k + D^{V}_{k+1}$.
\item[{(iii)}] If $D^{Q}_k > 0$ while $D^{V}_k > 0$, we have the following expression for post-advertisement reputation of period $k+1$, which is symmetric around $H$ and $L$ allowing $R^{'Q}_{k+1}$ and $R^{'V}_{k+1}$ to be equal.
\begin{equation}
\begin{split}
    & R^{'Q}_{k+1} = R^{'V}_{k+1}\\ \\
    &= R^V_k + H - HR^V_k -wR^V_k -wH + wHR^V_k +w -LR^V_k - HL + HLR^V_k + \\
    & wR^V_kL +wHL -wHLR^V_k - wL + L \\ \\
    & = R^V_k (1-w) + w + H(1 - R^V_k - w + wR^V_k) + L(1 - R^V_k - w + wR^V_k) + \\
    & HL (-1 + R^V_k + w - wR^V_k)
\end{split}
\end{equation}
In that case, we have $D^{Q}_k + D^{Q}_{k+1} = D^{V}_k + D^{V}_{k+1}$.
\end{itemize}

Since $D^{Q}_k + D^{Q}_{k+1} \geq D^{V}_k + D^{V}_{k+1}$, firm realizes more revenue in sequence $Q$.

We already know that cost of advertisement at $k$ and $k+1$ is $c_A \: (H + L)$ for both sequences. Then, we obtain:
\begin{equation}
    \pi^{Q}_k + \pi^{Q}_{k+1} \geq \pi^{V}_k + \pi^{V}_{k+1}
\end{equation}
Since $R^{'Q}_{k+1} \geq R^{'V}_{k+1}$, we have that $R^{Q}_{k+2} \geq R^{V}_{k+2}$. Then, by Corollary \ref{start w big R}, we have:
\begin{equation}
    \pi^{Q}_{k+2} + \pi^{Q}_{k+3} + ... \pi^{Q}_{N} \geq \pi^{V}_{k+2} + \pi^{V}_{k+3} + ... \pi^{V}_{N}  
\end{equation}
Finally, we conclude that the overall profit of the sequence $Q$ is an upper bound for $V$'s.
\section{Optimal Decision Thresholds for the Naive Firm}
\label{thresholds:naive}
\textbf{Variable Capacity:}
\begin{equation}
\addtolength{\jot}{0.5em}
    \begin{split}
        & \iota = \frac{\frac{u+p}{m} -H^\alpha}{1-H^\alpha}\\
        & \rho = \frac{\frac{u+p}{m} -L^\alpha}{1-L^\alpha} \\
        & \kappa = \frac{\frac{\Lambda (p-c)(u+p)}{m (\Lambda (p-c)-c_A(H-L))} - H^\alpha}{1-H^\alpha} \\
        & \tau = \Lambda (p-c)-c_A(H-L) \\
        & \nu \:\: st. \:\: \frac{1}{\nu+L^\alpha-\nu L^\alpha} - \frac{1}{\nu+H^\alpha-\nu H^\alpha} = \frac{m c_A (H-L)}{\Lambda (p-c)(u+p)}
    \end{split}
\end{equation}
\textbf{Constant Capacity:}
\begin{equation}
\addtolength{\jot}{0.5em}
    \begin{split}
    	& \gamma = \frac{\frac{u+p}{m} \frac{p \Lambda}{p \Lambda - c_A (H-L)} -H^\alpha}{1-H^\alpha}	\\
    	& \omega = \frac{\frac{p \Lambda}{p \Lambda - p S + c_A (H-L)} - L^\alpha}{1 - L^\alpha}\\
    	& \phi = p S-c_A(H-L)\\
    	& \beta, \:\: st. \:\:	\frac{1}{\beta + (1-\beta) L^\alpha} - \frac{1}{\beta + (1-\beta) H^\alpha} = \frac{m c_A (H-L)}{p \Lambda (u+p)}
    \end{split}
\end{equation}
\section{Proof of Proposition \ref{Prop:LBForH}}
\label{Proof:Prop:LBForH}
To begin with, let us introduce the following notation:
\begin{itemize}
	\item $H_{R_1}^C$: the optimal number of high level advertisements for the aware firm when the starting reputation is $R_1$,
	\item $H_{R_1}^I$: the optimal number of high level advertisements for the naive firm when the starting reputation is $R_1$.
\end{itemize}

	We will show that $H_{R_1}^C$ is (i) a non-decreasing function of $N$ and (ii) a non-increasing function of $R_1$. These results regarding the complete-information firm will be used in the upcoming parts of the proof of Proposition \ref{Prop:LBForH}.
	\begin{itemize}
		\item [(i)]
		
		When the optimal switching point of an $N$ period problem is $j^*$, it produces more $\Pi$ than any later switching point. Thus, $\Pi(j^*) \geq \Pi(j^*+1)$, and we have the following.
		\begin{equation}
			(p-c) \left(  \sum_{i=1}^{N}{ D_i^{j^*+1} - D_i^{j^*} }  \right) \leq c_A (H-L)
			\label{switch at j^* rather than j^* + 1}
		\end{equation}
		As N increases, the left-hand side of Equation \ref{switch at j^* rather than j^* + 1} increases as well. When $N$ increases enough so that the above inequality does not hold, the optimal switching point becomes greater than $j^*$, which proves that $H_{R_1}^C$ is a non-decreasing function of $N$.
		
		\item [(ii)]
			Let us again consider an N-period problem with $j^*$ as the optimal switching point. $j^*$ produces more $\Pi$ than any earlier switching point. Thus, $\Pi(j^*) \geq \Pi(j^*-1)$, and we have the following.
		\begin{equation}
			(p-c) \left( \sum_{i=1}^{N}{ D_i^{j^*} - D_i^{j^*-1} } \right) > c_A (H-L)
		\end{equation}
		
		As $R_1$ goes up, $R_{i}$ for all $i$ go up as well, and the firm will not be able to increase reputation with advertisement as much, due to diminishing returns of advertisement. Left-hand side of the above inequality will decrease, and the optimal switching point will be less than $j^*$, which proves that $H_{R_1}^C$ is a non-increasing function of $R_1$.
		
	\end{itemize}
	
	Next, we show that Proposition \ref{Prop:LBForH} holds under every condition:
	 
	\begin{itemize}
		\item [1. ] We first prove that Proposition \ref{Prop:LBForH} holds when $\tau \leq 0 \:\: or \:\: \iota \geq \rho$:
		
		$\iota \geq \rho$ implies that $\frac{u+p}{m} \geq 1$. Since $0 < R'_i < 1$, $\gamma_i = \Lambda \left(1 - \frac{u+p}{m R'_i} \right)$ becomes negative, which means it is not possible to realize demand at all. Then, making high level advertisement becomes meaningless, and the incomplete-information firm makes low level advertisement in all periods, which is optimal for complete-information firm as well. A similar situation occurs when we have $\tau \leq 0$; since it means that advertisement is so costly that making high level advertisement can never yield higher $\pi_i$ than low level advertisement. In both these special cases, firms focus on cutting costs, and none of them makes high level advertisement. Therefore, $H_{R_1}^C = H_{R_1}^I = 0$, and Proposition \ref{Prop:LBForH} holds.
		
		\item [2. ] Next, we will prove that Proposition \ref{Prop:LBForH} holds when $\tau > 0$ and $\iota < \rho$:
		
		\begin{itemize}
			\item [2.1. ] Suppose, $ R_1 \: \in \: \left[ \: max\{\rho, \: \nu\}, \: 1 \: \right]$.
			
			It is established that $R_i \leq R_{i+1}$ for all $i$. Therefore, when the sequence starts with $R_1 \in \left[ \: max\{\rho, \: \nu\}, \: 1 \: \right]$, we say that $ R_i \in \left[ \: max\{\rho, \: \nu\}, \: 1 \: \right]$ for all $i$. Then, the incomplete-information firm's advertisement decision is the same for every period, which is $A_i = L$. 
			
			As a result, $H_{R_1}^C \geq H_{R_1}^I = 0$, and Proposition \ref{Prop:LBForH} holds.
			
			\item [2.2. ] Suppose, $ R_1 \in \left(\:min \: \{ max\{ \iota, \: \kappa \}, \: \rho \}, \: max\{\rho, \: \nu\} \: \right)$.
			
			When we have $ R_1 \: \in \: \left(\:min \: \{ max\{ \iota, \: \kappa \}, \: \rho \}, \: max\{\rho, \: \nu\} \: \right)$; the incomplete-information firm makes high level advertisement. Reputation only increases during the problem horizon; therefore, after a while, the firm reaches a point where $R_i \: \in \: \left[ \: max\{\rho, \: \nu\}, \: 1 \: \right]$. Then, it continues with low level advertisements in the remaining periods.
			
			These imply that the incomplete-information firm adopts a \emph{single switching point solution}.
			
			Every time the incomplete-information firm decides on making high level advertisement at a period $i$, the complete-information firm does so as well, since it increases both $\pi_i$ and $\sum_{k = i+1}^{N} \pi_k$ (by Corollary \ref{start w big R}, since $R_{i+1}$ is higher after making high level advertisement).
			
			When the incomplete-information firm switches to low level advertisement at period $j$, we have $ \pi_j(j) \geq \pi_j(j+1)$, where $\pi_i(j)$ is the profit of period $i$ when switching point is at $j$. However, since $R_{j+1}$ will be greater if switching point is greater than $j$, we have that $\sum_{k = j+1}^{N} \pi_k(j+1) > \sum_{k = j+1}^{N} \pi_k(j)$. Then, if the following holds, the complete-information firm makes high level advertisement at $j$.
			\begin{equation}
				\sum_{k = j+1}^{N} \pi_k(j+1) - \sum_{k = j+1}^{N} \pi_k(j) > \pi_j(j) - \pi_j(j+1)
			\end{equation}
			
			Thus, we have $H_{R_1}^C \geq H_{R_1}^I$, and Proposition \ref{Prop:LBForH} holds.
			
			\item [2.3. ] Finally, suppose, $R_1 \: \in \: \left(\:0, \: min \: \{ max\{ \iota, \: \kappa \}, \: \rho \} \: \right]$.
			
			We consider the case where, we have significantly small reputation level at the first period. Thus, we should observe a pattern with \emph{two switching points}.
			
			Let us call the first switching point $s$, as in $A_s = H$ and $\: A_i = L$, for all $i < s$.
			
			At period $s$, we have $R_s \in \left(\:min \: \{ max\{ \iota, \: \kappa \}, \: \rho \}, \: max\{\rho, \: \nu\} \: \right)$.
			We can consider the periods after $s$ as an $N-s+1$ period problem with $R_s$ as the starting reputation. Then, this new problem has a single switching point, and by result 2.2 of this proof, we have $H_{R_s}^I \leq H_{R_s}^C$.
			
			The actual complete-information case has $N > N-s+1$ periods and starts with $R_1 < R_s$. Then, as $H_{R_1}^C$ is a non-decreasing function of $N$ and a non-increasing function of $R_1$, we have $H_{R_1}^C \geq H_{R_s}^C$. Finally, we obtain $H_{R_1}^C \geq H_{R_s}^C \geq H_{R_s}^I = H_{R_1}^I$, and Proposition \ref{Prop:LBForH} holds.
			
			We have shown that Proposition \ref{Prop:LBForH} holds in all possible cases, which concludes the proof.
		\end{itemize}
	\end{itemize}

\section{Demand Realization Ability for Naive and Aware Firms under Constant Capacity}
\label{App:Omega}
At the beginning of each time period, the firm makes an advertisement decision, affecting the demand generated in that period. We provide a classification for the status of the firm at the beginning of a period, in terms of its demand generating ability. To do so, we first introduce the following thresholds for reputation that help us with the classification:

\begin{itemize}
    \item $X_1 = \frac{\frac{u+p}{m} \frac{\Lambda}{\Lambda - S} - L^\alpha}{1 - L^\alpha}$ (min $R_i$ to have $D_i \geq S$, given $A_i=L$)\\~\
    \item  $X_2 = \frac{\frac{u+p}{m} \frac{\Lambda}{\Lambda - S} - H^\alpha}{1 - H^\alpha}$ (min $R_i$ to have $D_i \geq S$, given $A_i=H$)\\~\
    \item  $Y_1 = \frac{\frac{u+p}{m} - L^\alpha}{1 - L^\alpha}$ (min $R_i$ to have $D_i \neq 0$, given $A_i=L$)\\~\
    \item  $Y_2 = \frac{\frac{u+p}{m} - H^\alpha}{1 - H^\alpha}$ (min $R_i$ to have $D_i \neq 0$, given $A_i=H$)\\~\
\end{itemize}

By definition, it is clear that $X_1 > X_2$, $Y_1 > Y_2$, $X_1 > Y_1$ and $X_2 > Y_2$. However, the relationship between $Y_1$ and $X_2$ depends on the parameters.

Then, we define a set of \emph{states}, by comparing the amount of demand the firm can generate under high and low level advertisements with capacity $S$ and with $0$; in other words, by comparing the reputation with $X_1$, $X_2$, $Y_1$ and $Y_2$. At the beginning of any period $i$, the firm can be in one of the following states:

\textbf{State A:} $D_i \geq S$; meaning $R_i \geq X_1$.
        
\textbf{State B:} $D_i \geq S$ when $A_i = H$, and $S > D_i > 0$ when $A_i = L$; meaning $X_2 \leq R_i < X_1$ and $Y_1 < R_i$.

\textbf{State C:} $S > D_i > 0$; meaning $Y_1 < R_i < X_2$.

\textbf{State D:} $D_i \geq S$ when $A_i = H$, and $D_i = 0$ when $A_i = L$; meaning $X_2 \leq R_i \leq Y_1$.

\textbf{State E:} $S > D_i > 0$ when $A_i = H$, and $D_i = 0$ when $A_i = L$; meaning $Y_2 < R_i < X_2$ and $R_i \leq Y_1$.

\textbf{State F:} $D_i = 0$; meaning $R_i \leq Y_2$.
 
If we have $Y_1 > X_2$, occurrences of the states A, B, C, D, E and F, on a reputation axis, can be shown below. For instance, if $X_2 < R_i < Y_1$ at some period $i$, the firm is at state D.
    \input{shapes/y1_greater_than_x2}
Otherwise, if $Y_1 < X_2$, states A, B, C, D, E and F occur by the following reputation axis.
    \input{shapes/y1_less_than_x2}

Let us denote the set of states that can be observed with a given parameter set by $\Omega$. Then, we can come up with the following results.
\begin{itemize}
	\item[{(i)}] If $Y_1 > X_2$, $A \in \Omega$.
	\begin{proof}
		We have assumed that $\frac{u+p}{m} < 1$; therefore, by definition, we have $Y_1 < 1$ and $Y_2 < 1$.
		Then, when $Y_1 > X_2$, we say $1 > X_2$ as well, which indicates that $\frac{u+p}{m} \frac{\Lambda}{\Lambda - S} < 1$. Therefore, $X_1 < 1$ must hold, meaning that there must be some $1 > \epsilon > 0$ such that $X_1 + \epsilon < 1$, in other words, state A is always observable.
	\end{proof}
	\item[{(ii)}] If $Y_1 < X_2$; and if $\{E,\:F\} \cap \Omega \neq \emptyset$, then $C \in \Omega$.
	\begin{proof}
		We have $\frac{u+p}{m} < 1$; therefore, $Y_2 < Y_1 < 1$.
		$\{E,\:F\} \cap \Omega \neq \empty$ implies that reputation can take smaller values than $Y_1$; therefore, $Y_1 > 0$.
		Then, there must be some $1 > \epsilon > 0$ such that $0 < Y_1 + \epsilon < 1$, allowing state C to occur. We conclude that $\Omega$ cannot include $E$ and/or $F$ while excluding $C$.
	\end{proof}
	\item[{(iii)}] If $Y_1 < X_2$; and if $B \in \Omega$, then $A \in \Omega$.
	\begin{proof}
		If $B \in \Omega$, we should have $1 > X_2$, meaning $\frac{u+p}{m} \frac{\Lambda}{\Lambda - S} < 1$. This indicates that $X_1 < 1$. Then, there must be some $1 > \epsilon > 0$ such that $X_1 + \epsilon < 1$, implying that $A \in \Omega$.
	\end{proof}
\end{itemize}

Given a parameter set, cardinality of $\Omega$ can be reduced by ruling out some of the states if the corresponding reputation intervals of those states do not intersect with $[0,\:1]$. Using the results above, $\Omega$ of a given parameter set can be one of the following.
\begin{itemize}
	\item If $Y_1 > X_2$:
	\begin{itemize}
		\item $\Omega = \{A, \: B, \: D, \: E, \: F\}$, if $\: H^\alpha < \frac{u+p}{m} $ and $\frac{u+p}{m} \frac{\Lambda}{\Lambda - S} < 1$.
		\item $\Omega = \{A, \: B, \: D, \: E\}$, if $\: H^\alpha \geq \frac{u+p}{m} > L^\alpha $ and $H^\alpha < \frac{u+p}{m} \frac{\Lambda}{\Lambda - S} < 1$.
		\item $\Omega = \{A, \: B, \: D\}$, if $\: H^\alpha \geq \frac{u+p}{m} > L^\alpha $ and $\frac{u+p}{m} \frac{\Lambda}{\Lambda - S} \leq H^\alpha$.
		\item $\Omega = \{A, \: B\}$, if $\: L^\alpha \geq \frac{u+p}{m} $ and $L^\alpha < \frac{u+p}{m} \frac{\Lambda}{\Lambda - S} \leq H^\alpha$.
		\item $\Omega = \{A\}$, if $\: L^\alpha \geq \frac{u+p}{m} \frac{\Lambda}{\Lambda - S} $.
	\end{itemize}
	\item If $Y_1 < X_2$:
	\begin{itemize}
		\item $\Omega = \{A, \: B, \: C, \: E, \: F\}$, if $H^\alpha < \frac{u+p}{m} $ and $\frac{u+p}{m} \frac{\Lambda}{\Lambda - S} < 1$.
		\item $\Omega = \{A, \: B, \: C, \: E\}$, if $\: H^\alpha \geq \frac{u+p}{m} > L^\alpha $ and $H^\alpha < \frac{u+p}{m} \frac{\Lambda}{\Lambda - S} < 1$.
		\item $\Omega = \{A, \: B, \: C\}$, if $\: L^\alpha \geq \frac{u+p}{m} $ and $H^\alpha < \frac{u+p}{m} \frac{\Lambda}{\Lambda - S} < 1$.
		\item $\Omega = \{A, \: B\}$, if $\: L^\alpha \geq \frac{u+p}{m} $ and $L^\alpha < \frac{u+p}{m} \frac{\Lambda}{\Lambda - S} \leq H^\alpha$.
		\item $\Omega = \{A\}$, if $\: L^\alpha \geq \frac{u+p}{m} \frac{\Lambda}{\Lambda - S}$.
		\item $\Omega = \{C, \: E, \: F\}$, if $\: H^\alpha < \frac{u+p}{m}$ and $\frac{u+p}{m} \frac{\Lambda}{\Lambda - S} \geq 1$.
		\item $\Omega = \{C, \: E\}$, if $\: H^\alpha \geq \frac{u+p}{m} > L^\alpha $ and $\frac{u+p}{m} \frac{\Lambda}{\Lambda - S} \geq 1$.
		\item $\Omega = \{C\}$, if $\: L^\alpha \geq \frac{u+p}{m} $ and $\frac{u+p}{m} \frac{\Lambda}{\Lambda - S} \geq 1$.
	\end{itemize}
\end{itemize}

Finally, $\Omega$ has to be either one of: $\{A\}$, $\{A, \: B\}$, $\{A, \: B, \: C\}$, $\{A, \: B, \: C, \: E\}$, $\{A, \: B, \: C, \: E, \: F\}$, $\{A, \: B, \: D\}$, $\{A, \: B, \: D, \: E\}$, $\{A, \: B, \: D, \: E, \: F\}$, $\{C\}$, $\{C, \: E\}$, $\{C, \: E, \: F\}$.

\section{Proof of Proposition \ref{Prop:ConstCapCEF}}
\label{Proof:Prop:ConstCapCEF}

	We have $R_{i+1} = R^{'}_i$ when $D_i=0$, and $R_{i+1} = w F_i + (1-w) R'_i$ when $D_i>0$. Therefore, for simplicity, we can use the following function to set fill rate, and stick with the formula $R_{i+1} = w F_i + (1-w) R'_i$ to determine $R_{i+1}$.
	\begin{equation}
		F_i = \begin{cases}
			R^{'}_i, & D_i=0\\
			\frac{min\{S,\:D_i\}}{D_i}, & D_i>0
		\end{cases}
	\end{equation}
	Since $D_i \leq S$ always holds, we have that $F_i = R^{'}_i$ or $F_i=1$ for all $i$. As a result, $R_{i+1} \geq R_i$, which implies that $F_i = R^{'}_i$ is true up to the first period demand is observed and from that point and on, we have $F_i = 1$.
	
	Suppose we have an advertisement sequence $V = \{V_1, \: V_2, \: ... \: V_N \}$, where $V_k = L$ and $V_{k+1} = H$, for some $k$, in an N-period problem. 
	Then, using $V$, let us define a second sequence $Q$ such that $Q_i = V_i$, for all $i \; {\in} \; \{1, \: 2, \: ..., \; k-1, \; k+2, \; ... \: N \}  $ while $Q_k = H, \; Q_{k+1} = L$.
	
	We denote the parameters such as reputation, post-advertisement reputation, fill rate, demand and profit at period $i$, using the sequence information (e.g. firm's reputation at sequence $V$ and period $i$ is represented by $R_i^V$).
	
	In both sequences, the number of times high and low level advertisements are made and the total capacity costs are the same. Therefore, our sole concern in comparing the sequences is the total revenue, which is expressed as $\sum_{i=1}^{N} p \: D_i$, since demand cannot exceed capacity. As $p$ is constant as well, we can reduce the comparison to the sum of demands, that is $\sum_{i=1}^{N} D^V_i$ versus $\sum_{i=1}^{N} D^Q_i$.
	
	Up to period $k$, the two sequences are equivalent, always making the same decisions. Thus, $\sum_{i=1}^{k-1} D^Q_i = \sum_{i=1}^{k-1} D^V_i$.
	
	At period $k$, we have that $Q_k = H$ while $V_k = L$; therefore, $R^{'V}_k \leq R^{'Q}_k$ and consequently $F^{V}_k\leq F^{Q}_k$ and $D^{V}_k \leq D^{Q}_k$.
	
	At period $k+1$, we have the following post-advertisement reputation for sequence $V$, where $F^V_k = 1$ if $D_k>0$ and $F^V_k = R_k + (1-R_k)H^\alpha$ otherwise:
	\begin{equation}
		R^{'V}_{k+1} = H^\alpha + (1-H^\alpha) (w F^V_k + (1-w) (L^\alpha + (1-L^\alpha)R_k)),
	\end{equation}
	and the following post-advertisement reputation for sequence $Q$, where $F^Q_k = 1$ if $D_k>0$ and $F^Q_k = R_k + (1-R_k)L^\alpha$ otherwise:
	\begin{equation}
		R^{'Q}_{k+1} = L^\alpha + (1-L^\alpha) (w F^Q_k + (1-w) (H^\alpha + (1-H^\alpha)R_k)).
	\end{equation}
	We have $R^{'V}_{k+1} < R^{'Q}_{k+1}$ if $F^V_k = R_k + (1-R_k)H^\alpha$ while $F^Q_k = 1$, and otherwise $R^{'V}_{k+1} = R^{'Q}_{k+1}$. Then, $F^{V}_{k+1}\leq F^{Q}_{k+1}$ and $D^{V}_{k+1} \leq D^{Q}_{k+1}$.
	
	Since $R^{'V}_{k+1} \leq R^{'Q}_{k+1}$ and $F^{V}_{k+1}\leq F^{Q}_{k+1}$, we have $R^{V}_{k+2} \leq R^{Q}_{k+2}$. Therefore, for all remaining periods we can say that $D^Q_{i} \geq D^V_{i}$, where $i > k+1$.
	
	Finally, we conclude that $\sum_{i=1}^{N} D^Q_i \geq \sum_{i=1}^{N} D^V_i$, leading to $\Pi^Q \geq \Pi^V$.

\section{Proof of Lemma \ref{Lemma:ConstCapABD}}
\label{Proof:Lemma:ConstCapABD}
	Let us denote the number of high level advertisements this heuristic results in by $h^*$.
	
	Following the given heuristic, firm obtains the maximum possible revenue, that is $p S N$. Making more high level advertisements would not increase the total revenue but it would increase the total cost. Therefore, the above heuristic provides a better solution than any advertising sequence with more high level advertisements than $h^*$.

\section{Proof of Proposition \ref{Prop:AbsorbingA}}
\label{Proof:Prop:AbsorbingA}
Dropping from state A means, next period begins with a reputation level less than $X_1$, as in $R_{i+1} < X_1$, which is expressed as follows.
\begin{equation}
	w \frac{S}{\Lambda \left(1 - \frac{u+p}{m R'_i} \right)} + (1-w) R'_i < X_1
\end{equation}

By working on this expression, we obtain the following.
\begin{equation}
	0 > (1-w) R_i^{'2} + \left(\frac{w S}{\Lambda} - X_1 - (1-w) \frac{u+p}{m}\right) R'_i + X_1 \frac{u+p}{m}
	\label{dropping from A}
\end{equation}

Using the quadratic formula, we obtain the roots for post-advertisement reputation that satisfy Equation \ref{dropping from A}:
\begin{equation}
	\resizebox{.7\hsize}{!}{$r'_{1,2} = \frac{- \left(\frac{w S}{\Lambda} - X_1 - (1-w) \frac{u+p}{m}\right) \pm \sqrt{\left(\frac{w S}{\Lambda} - X_1 - (1-w) \frac{u+p}{m}\right)^2 - 4 (1-w) X_1 \frac{u+p}{m}}}{2 (1-w)}$}
\end{equation}
Within the interval defined by the roots, the firm drops from state A, while outside the interval, it remains at state A. Note that these roots are for the post-advertisement reputation. The corresponding $R_i$ roots can be obtained as follows.
\begin{equation}
	r_{1,2} = \frac{\left(\frac{- \left(\frac{w S}{\Lambda} - X_1 - (1-w) \frac{u+p}{m}\right) \pm \sqrt{\left(\frac{w S}{\Lambda} - X_1 - (1-w) \frac{u+p}{m}\right)^2 - 4 (1-w) X_1 \frac{u+p}{m}}}{2 (1-w)}\right) - A_i^\alpha}{1 - A_i^\alpha}
\end{equation}

Proposition \ref{Prop:AbsorbingA} results from the fact that state A brings the maximum possible revenue ($p S$) regardless of the advertisement level, while low level advertisement has the lowest advertisement cost.

As a result, when the given interval to be at A does not overlap with the interval to leave A after making low level advertisement, firm is able to stay at state A for the lowest possible cost forever. Then, it always makes low level advertisement once state A is observed for the first time.

\small
\section{Sample Experiment for the Optimal Advertising Decisions with Respect to Advertisement Resistance $\alpha$ under Variable Capacity}
\label{exp:varcap:alpha}
$R_1=0.57$, $L=0.28$, $H=0.36$, $w=0.85$, $\Lambda = 1,528$, $c=0.33$, $m=37.03$, $u=2.76$, $p=18.25$, $c_A=623.22$, $N=29$.
\input{tables/voi_var_alpha}

\section{Sample Experiment 1 for Variable Capacity Firms' Optimal Advertising Decisions with Respect to Starting Reputation $R_1$}
\label{exp:varcap:r1_sample1}
$L=0.09$, $H=0.66$, $\alpha=7.12$, $w=0.02$, $\Lambda = 971$, $c=7.69$, $m=29.11$, $u=1.61$, $p=13.37$, $c_A=875.30$, $N=7$.
\input{tables/voi_var_r1_sample1}

%\section{Sample Experiment 2 for Variable Capacity Firms' Optimal Advertising Decisions with Respect to Starting Reputation}
%\label{exp:varcap:r1_sample2}
%$L=0.68$, $H=0.83$, $\alpha=8.24$, $w=0.65$, $\Lambda = 1,041$, $c=7.16$, $m=32.14$, $u=4.04$, $p=16.94$, $c_A=645.69$, $N=29$.
%\input{tables/voi_var_r1_sample2}

\section{Sample Experiment 1 for Variable Capacity Firms' Optimal Advertising Decisions with Respect to Price $p$}
\label{exp:varcap:p_sample1}
$R_1= 0.37$, $L=0.21$, $H=0.28$, $\alpha=3.11$, $w=0.20$, $\Lambda = 1,372$, $c=9.06$, $m=30.23$, $u=4.02$, $c_A=651.60$, $N=46$.
\input{tables/voi_var_p_sample1}

\section{Sample Experiment 2 for Variable Capacity Firms' Optimal Advertising Decisions with Respect to Price $p$}
\label{exp:varcap:p_sample2}
$R_1= 0.75$, $L=0.55$, $H=0.68$, $\alpha=8.23$, $w=0.11$, $\Lambda = 2,494$, $c=3.55$, $m=23.91$, $u=1.41$, $c_A=913.86$, $N=26$.
\input{tables/voi_var_p_sample2}

\section{Sample Experiments for Constant Capacity Firms' Optimal Advertising Decisions with Respect to Service Capacity $S$}
$R_1= 0.18$, $L=0.68$, $H=0.82$, $\alpha=6.04$, $w=0.99$, $\Lambda = 653$, $c=0.70$, $m=39.22$, $u=2.70$, $p=16.13$ $c_A=731.85$, $N=18$.
\label{exp:concap:s1_sample}

\input{tables/voi_con_s1}

$R_1= 0.18$, $L=0.61$, $H=0.80$, $\alpha=6.71$, $w=0.01$, $\Lambda = 2,221$, $c=2.05$, $m=26.99$, $u=4.92$, $p=7.18$ $c_A=246.03$, $N=11$.
%\label{exp:concap:s2_sample}
\input{tables/voi_con_s2}

$R_1= 0.45$, $L=0.02$, $H=0.20$, $\alpha=2.03$, $w=0.78$, $\Lambda = 1,697$, $c=9.18$, $m=34.58$, $u=2.06$, $p=15.08$ $c_A=811.36$, $N=19$.
%\label{exp:concap:s3_sample}
\input{tables/voi_con_s3}

%\section{Sample Experiment 2 for Constant Capacity Firms' Optimal Advertising Decisions with Respect to Service Capacity}
%$R_1= 0.18$, $L=0.61$, $H=0.80$, $\alpha=6.71$, $w=0.01$, $\Lambda %= 2,221$, $c=2.05$, $m=26.99$, $u=4.92$, $p=7.18$ $c_A=246.03$, %$N=11$.
%\label{exp:concap:s2_sample}
%\input{tables/voi_con_s2}

%\section{Sample Experiment 3 for Constant Capacity Firms' Optimal %Advertising Decisions with Respect to Service Capacity}
%$R_1= 0.45$, $L=0.02$, $H=0.20$, $\alpha=2.03$, $w=0.78$, $\Lambda %= 1,697$, $c=9.18$, $m=34.58$, $u=2.06$, $p=15.08$ $c_A=811.36$, %$N=19$.
%\label{exp:concap:s3_sample}
%\input{tables/voi_con_s3}

\section{Sample Experiment for Constant Capacity Firms' Optimal Advertising Decisions with Respect to Weight of WoM $w$}
$R_1= 0.55$, $L=0.26$, $H=0.62$, $\alpha=7.47$, $\Lambda = 2,347$, $c=6.81$, $m=35.39$, $u=4.80$, $p=11.05$ $c_A=369.52$, $N=24$, $S=279$.
\label{exp:concap:w_sample}
\input{tables/voi_con_w}

\section{Sample Experiment for Constant Capacity Firms' Optimal Advertising Decisions with Respect to Advertisement Resistance $\alpha$}
$R_1= 0.01$, $L=0.53$, $H=0.79$, $w=0.42$, $\Lambda = 1,451$, $c=6.30$, $m=28.57$, $u=0.37$, $p=15.52$ $c_A=829.40$, $N=11$, $S=587$.
\label{exp:concap:alpha_sample}
\input{tables/voi_con_alpha}

%% file: shapes/y1_greater_than_x2.tex
\begin{center}
			\tikzset{every picture/.style={line width=0.75pt}} %set default line width to 0.75pt        
			
			\begin{tikzpicture}[x=0.75pt,y=0.75pt,yscale=-1,xscale=1]
				%uncomment if require: \path (0,432); %set diagram left start at 0, and has height of 432
				
				%Straight Lines [id:da5222053419914638] 
				\draw    (50,111) -- (448,111) ;
				%Straight Lines [id:da12268226453737041] 
				\draw    (50,107) -- (50,115.71) ;
				%Straight Lines [id:da49752082281503074] 
				\draw    (208,107) -- (208,115.71) ;
				%Straight Lines [id:da7444958367340222] 
				\draw    (120,107) -- (120,115.71) ;
				%Straight Lines [id:da22519208182918504] 
				\draw    (290,107) -- (290,115.71) ;
				%Straight Lines [id:da6759212733256537] 
				\draw    (377,106) -- (377,114.71) ;
				\draw   (442,108) -- (449,111) -- (442,114) ;
				
				% Text Node
				\draw (369,81) node [anchor=north west][inner sep=0.75pt]   [align=left] {{\small $\displaystyle X_{1}$}};
				% Text Node
				\draw (283,82) node [anchor=north west][inner sep=0.75pt]   [align=left] {{\small $\displaystyle Y_{1}$}};
				% Text Node
				\draw (201,81) node [anchor=north west][inner sep=0.75pt]   [align=left] {{\small $\displaystyle X_{2}$}};
				% Text Node
				\draw (115,82) node [anchor=north west][inner sep=0.75pt]   [align=left] {{\small $\displaystyle Y_{2}$}};
				% Text Node
				\draw (418.6,94.31) node  [color={rgb, 255:red, 222; green, 24; blue, 24 }  ,opacity=1 ] [align=left] {\begin{minipage}[lt]{15.78pt}\setlength\topsep{0pt}
						$\displaystyle A$
				\end{minipage}};
				% Text Node
				\draw (332.6,93.31) node  [color={rgb, 255:red, 222; green, 24; blue, 24 }  ,opacity=1 ] [align=left] {\begin{minipage}[lt]{15.78pt}\setlength\topsep{0pt}
						$\displaystyle B$
				\end{minipage}};
				% Text Node
				\draw (257.6,93.31) node  [color={rgb, 255:red, 222; green, 24; blue, 24 }  ,opacity=1 ] [align=left] {\begin{minipage}[lt]{15.78pt}\setlength\topsep{0pt}
						$\displaystyle D$
				\end{minipage}};
				% Text Node
				\draw (163.6,95.31) node  [color={rgb, 255:red, 222; green, 24; blue, 24 }  ,opacity=1 ] [align=left] {\begin{minipage}[lt]{15.78pt}\setlength\topsep{0pt}
						$\displaystyle E$
				\end{minipage}};
				% Text Node
				\draw (81.6,96.31) node  [color={rgb, 255:red, 222; green, 24; blue, 24 }  ,opacity=1 ] [align=left] {\begin{minipage}[lt]{15.78pt}\setlength\topsep{0pt}
						$\displaystyle F$
				\end{minipage}};
				% Text Node
				\draw (454,101) node [anchor=north west][inner sep=0.75pt]   [align=left] {{\small $\displaystyle R_{i}$}};
			\end{tikzpicture}
		\end{center}

%% file: shapes/y1_less_than_x2.tex
\begin{center}
			\tikzset{every picture/.style={line width=0.75pt}} %set default line width to 0.75pt        
			
			\begin{tikzpicture}[x=0.75pt,y=0.75pt,yscale=-1,xscale=1]
				%uncomment if require: \path (0,432); %set diagram left start at 0, and has height of 432
				
				%Straight Lines [id:da5222053419914638] 
				\draw    (50,111) -- (448,111) ;
				%Straight Lines [id:da12268226453737041] 
				\draw    (50,107) -- (50,115.71) ;
				%Straight Lines [id:da49752082281503074] 
				\draw    (208,107) -- (208,115.71) ;
				%Straight Lines [id:da7444958367340222] 
				\draw    (120,107) -- (120,115.71) ;
				%Straight Lines [id:da22519208182918504] 
				\draw    (290,107) -- (290,115.71) ;
				%Straight Lines [id:da6759212733256537] 
				\draw    (377,106) -- (377,114.71) ;
				\draw   (442,108) -- (449,111) -- (442,114) ;
				
				% Text Node
				\draw (369,81) node [anchor=north west][inner sep=0.75pt]   [align=left] {{\small $\displaystyle X_{1}$}};
				% Text Node
				\draw (199,81) node [anchor=north west][inner sep=0.75pt]   [align=left] {{\small $\displaystyle Y_{1}$}};
				% Text Node
				\draw (281,82) node [anchor=north west][inner sep=0.75pt]   [align=left] {{\small $\displaystyle X_{2}$}};
				% Text Node
				\draw (115,82) node [anchor=north west][inner sep=0.75pt]   [align=left] {{\small $\displaystyle Y_{2}$}};
				% Text Node
				\draw (418.6,94.31) node  [color={rgb, 255:red, 222; green, 24; blue, 24 }  ,opacity=1 ] [align=left] {\begin{minipage}[lt]{15.78pt}\setlength\topsep{0pt}
						$\displaystyle A$
				\end{minipage}};
				% Text Node
				\draw (332.6,93.31) node  [color={rgb, 255:red, 222; green, 24; blue, 24 }  ,opacity=1 ] [align=left] {\begin{minipage}[lt]{15.78pt}\setlength\topsep{0pt}
						$\displaystyle B$
				\end{minipage}};
				% Text Node
				\draw (257.6,93.31) node  [color={rgb, 255:red, 222; green, 24; blue, 24 }  ,opacity=1 ] [align=left] {\begin{minipage}[lt]{15.78pt}\setlength\topsep{0pt}
						$\displaystyle C$
				\end{minipage}};
				% Text Node
				\draw (163.6,95.31) node  [color={rgb, 255:red, 222; green, 24; blue, 24 }  ,opacity=1 ] [align=left] {\begin{minipage}[lt]{15.78pt}\setlength\topsep{0pt}
						$\displaystyle E$
				\end{minipage}};
				% Text Node
				\draw (81.6,96.31) node  [color={rgb, 255:red, 222; green, 24; blue, 24 }  ,opacity=1 ] [align=left] {\begin{minipage}[lt]{15.78pt}\setlength\topsep{0pt}
						$\displaystyle F$
				\end{minipage}};
				% Text Node
				\draw (454,101) node [anchor=north west][inner sep=0.75pt]   [align=left] {{\small $\displaystyle R_{i}$}};

			\end{tikzpicture}   
        \end{center}

%% file: tables/voi_var_alpha.tex
	\begin{longtable} {c|c|c}
 \hline\hline $\alpha$   &   \textbf{ aware firm's $SP_0$}   &  \textbf{naive firm's $SP_1$ and $SP_2$}\\
 \hline  0.52      &    2	                &  $SP_2 = 2$ (single switching action)\\
 \hline  0.69      &    2	                &  $SP_2 = 2$ (single switching action)\\
 \hline  1.16      &    3	                &  $SP_2 = 3$ (single switching action)\\
 \hline  1.27      &    3	                &  $SP_2 = 3$ (single switching action)\\
 \hline  1.39      &    3	                &  $SP_2 = 3$ (single switching action)\\
 \hline  1.63      &    3	                &  $SP_2 = 3$ (single switching action)\\
 \hline  2.06      &    3	                &  $SP_2 = 2$ (single switching action)\\
 \hline  2.08      &    3	                &  $SP_2 = 2$ (single switching action)\\
 \hline  2.30      &    2	                &  $SP_2 = 2$ (single switching action)\\
 \hline  2.43      &    2	                &  $SP_2 = 2$ (single switching action)\\
 \hline  2.52      &    2	                &  $SP_2 = 2$ (single switching action)\\
 \hline  2.65      &    2	                &  $SP_2 = 2$ (single switching action)\\
 \hline  2.89      &    2	                &  $SP_2 = 2$ (single switching action)\\
 \hline  3.12      &    2	                &  $SP_2 = 2$ (single switching action)\\
 \hline  3.24      &    2	                &  $SP_2 = 2$ (single switching action)\\
 \hline  3.44      &    2	                &  $SP_2 = 2$ (single switching action)\\
 \hline  3.67      &    2	                &  $SP_2 = 2$ (single switching action)\\
 \hline  3.67      &    2	                &  $SP_2 = 2$ (single switching action)\\
 \hline  3.68      &    2	                &  $SP_2 = 2$ (single switching action)\\
 \hline  4.09      &    2	                &  $SP_2 = 2$ (single switching action)\\
 \hline  4.12      &    2	                &  $SP_2 = 2$ (single switching action)\\
 \hline  4.39      &    2	                &  $SP_2 = 2$ (single switching action)\\
 \hline  5.00      &    2	                &  $SP_1 = 3$, $SP_2 = 4$ (two switching points)\\
 \hline  5.09      &    2	                &  $SP_1 = 3$, $SP_2 = 4$ (two switching points)\\
 \hline  5.10      &    2	                &  $SP_1 = 3$, $SP_2 = 4$ (two switching points)\\
 \hline  5.12      &    2	                &  $SP_1 = 3$, $SP_2 = 4$ (two switching points)\\
 \hline  5.20      &    3	                &  $SP_1 = 4$, $SP_2 = 5$ (two switching points)\\
 \hline  5.33      &    3	                &  $SP_1 = 5$, $SP_2 = 6$ (two switching points)\\
 \hline  5.44      &    3	                &  $SP_1 = 6$, $SP_2 = 7$ (two switching points)\\
 \hline  5.55      &    3	                &  $SP_1 = 7$, $SP_2 = 8$ (two switching points)\\
 \hline  6.11      &    4	                &  always $L$\\
 \hline  6.31      &    5	                &  always $L$\\
 \hline  6.59      &    6	                &  always $L$\\
 \hline  6.61      &    6	                &  always $L$\\
 \hline  7.15      &    9	                &  -\\
 \hline  7.23      &    10	                &  -\\
 \hline  7.32      &    11	                &  -\\
 \hline  7.39      &    11	                &  -\\
 \hline  7.53      &    13	                &  -\\
 \hline  7.78      &    16	                &  -\\
 \hline  7.94      &    18	                &  -\\
 \hline  8.11      &    22	                &  -\\
 \hline  8.47      &    -	&  -\\
 \hline  ...       &    -	&  -\\
% \hline  9.24      &    -	&  -\\
% \hline  9.25      &    -	&  -\\
% \hline  9.42      &    -	&  -\\
% \hline  9.71      &    -	&  -\\
 \hline \hline
	\end{longtable}

%% file: tables/voi_var_r1_sample1.tex
	\begin{longtable} 
		 {c|c|c}
		\hline\hline $R_1$     &    \textbf{aware firm's $SP_0$}   &  \textbf{naive firm's $SP_1$ and $SP_2$}\\
		\hline ...      & -         &  - \\
	%	\hline 0.4098   & -         &  - \\
	%	\hline 0.4352   & -         &  - \\
	%	\hline 0.4504   & -         &  - \\
		\hline 0.45   & 6	                       &  - \\
		\hline 0.46   & 6	                       &  - \\
		\hline 0.47   & 5	                       &  - \\
		\hline 0.51   & 5	                       &  - \\
		\hline 0.54   & 4	                       &  always $L$ \\
		\hline 0.56  & 4	                       &  always $L$ \\
		\hline 0.59  & 4	                       &  always $L$ \\
		\hline 0.59   & 4	                       &  always $L$ \\
		\hline 0.59   & 4	                       &  always $L$ \\
		\hline 0.60   & 4	                       &  always $L$ \\
		\hline 0.63   & 3	                       &  always $L$ \\
		\hline 0.64   & 3	                       &  always $L$ \\
		\hline 0.70   & 2	                       &  always $L$ \\
		\hline 0.75   & always $L$         &  always $L$ \\
		\hline ...      & always $L$         &  always $L$ \\
%		\hline 0.7655   & always $L$         &  always $L$ \\
%		\hline 0.8175   & always $L$         &  always $L$ \\
%		\hline 0.8181   & always $L$         &  always $L$ \\
%		\hline 0.8354   & always $L$         &  always $L$ \\
%		\hline 0.8421   & always $L$         &  always $L$ \\
%		\hline 0.8462   & always $L$         &  always $L$ \\
%		\hline 0.8629   & always $L$         &  always $L$ \\
%		\hline 0.8795   & always $L$         &  always $L$ \\
%		\hline 0.9304   & always $L$         &  always $L$ \\
%		\hline 0.9441   & always $L$         &  always $L$ \\
%		\hline 0.9955   & always $L$         &  always $L$ \\
		\hline\hline
	\end{longtable}

%% file: tables/voi_var_p_sample1.tex
	\begin{longtable} {c|c|c}
		\hline\hline $p$     &    \textbf{aware firm's $SP_0$}   &  \textbf{naive firm's $SP_1$ and $SP_2$}\\
		\hline  9.11   &    -	& -\\
		\hline  9.23   &    -	& -\\
		\hline  9.34   &    7					& -\\
		\hline  9.39   &    7	                & always $L$ \\
		\hline  9.44   &    8	                & always $L$ \\
		\hline  10.03  &    10	                & always $L$ \\
		\hline  10.12  &    10	                & always $L$ \\
		\hline  10.67  &    12	                & always $L$ \\
		\hline  11.10  &    14	                & always $L$ \\
		\hline  11.19  &    15	                & always $L$ \\
		\hline  11.78  &    18	                & always $L$ \\
		\hline  11.89  &    18	                & always $L$ \\
		\hline  13.12  &    24	                & - \\
		\hline  13.90  &    27	                & - \\
		\hline  14.02  &    27	                & - \\
		\hline  14.30  &    29	                & - \\
		\hline  14.42  &    30	                & - \\
		\hline  14.63  &    31	                & - \\
		\hline  15.04  &    33	                & - \\
		\hline  15.38  &    34	                & - \\
		\hline  15.56  &    35	                & - \\
		\hline  15.60  &    35	                & - \\
		\hline  15.72  &    36	                & - \\
		\hline  16.00  &    37	                & - \\
		\hline  16.37  &    39	                & - \\
		\hline  16.54  &    40	                & - \\
		\hline  16.69  &    41	                & - \\
		\hline  16.92  &    41	                & - \\
		\hline  16.99  &    42	                & - \\
		\hline  17.15  &    42	                & - \\
		\hline  17.17  &    42	                & - \\
		\hline  17.60  &    44	                & - \\
		\hline  17.70  &    44	                & - \\
		\hline  17.73  &    -	& - \\
		\hline  ...    &    -	& - \\
		% \hline  17.84  &    -	& - \\
		% \hline  18.02  &    -	& - \\
		% \hline  18.69  &    -	& - \\
		% \hline  19.27  &    -	& - \\
		% \hline  19.36  &    -	& - \\
		% \hline  19.53  &    -	& - \\
		% \hline  20.00  &    -	& - \\
		% \hline  20.47  &    -	& - \\
		% \hline  20.64  &    -	& - \\
		% \hline  20.75  &    -	& - \\
		% \hline  20.95  &    -	& - \\
		% \hline  21.79  &    -	& - \\
		% \hline  22.40  &    -	& - \\
		% \hline  22.66  &    -	& - \\
		% \hline  23.64  &    -	& - \\
		% \hline  24.23  &    -	& - \\
		% \hline  24.29  &    -	& - \\
		% \hline  24.66  &    -	& - \\
		\hline\hline
	\end{longtable}

%% file: tables/voi_var_p_sample2.tex
	\begin{longtable} 
		 {c|c|c}
		\hline\hline $p$     &    \textbf{aware firm's $SP_0$}   &  \textbf{naive firm's $SP_1$ and $SP_2$}\\
		\hline   3.60	&     -	  &  -\\
		\hline   3.86	&     always $L$	  &  always $L$\\
		\hline   4.18	&     always $L$	  &  always $L$\\
		\hline   4.53	&     always $L$	  &  always $L$\\
		\hline   5.16	&     always $L$	  &  always $L$\\
		\hline   5.30	&     2	                  &  always $L$\\
		\hline   5.99	&     4	                  &  always $L$\\
		\hline   7.23	&     7	                  &  always $L$\\
		\hline   7.55	&     7	                  &  always $L$\\
		\hline   8.93	&     10                  &  always $L$\\
		\hline   8.98	&     10                  &  always $L$\\
		\hline   9.14	&     10                  &  always $L$\\
		\hline   9.16	&     10                  &  always $L$\\
		\hline   9.89	&     11                  &  always $L$\\
		\hline   10.53	&     12                  &  $SP_2 = 2$ (single switching action)\\
		\hline   10.86	&     12                  &  $SP_2 = 2$ (single switching action)\\
		\hline   11.16	&     13                  &  $SP_2 = 3$ (single switching action)\\
		\hline   11.23	&     13                  &  $SP_2 = 3$ (single switching action)\\
		\hline   11.23	&     13                  &  $SP_2 = 3$ (single switching action)\\
		\hline   11.43	&     13                  &  $SP_2 = 3$ (single switching action)\\
		\hline   11.64	&     13                  &  $SP_2 = 3$ (single switching action)\\
		\hline   12.63	&     14                  &  $SP_2 = 4$ (single switching action)\\
		\hline   12.82	&     14                  &  $SP_2 = 4$ (single switching action)\\
		\hline   13.43	&     15                  &  $SP_2 = 4$ (single switching action)\\
		\hline   13.71	&     15                  &  $SP_2 = 5$ (single switching action)\\
		\hline   13.78	&     15                  &  $SP_2 = 5$ (single switching action)\\
		\hline   13.90	&     15                  &  $SP_2 = 5$ (single switching action)\\
		\hline   13.96	&     15                  &  $SP_2 = 5$ (single switching action)\\
		\hline   14.28	&     16                  &  $SP_2 = 5$ (single switching action)\\
		\hline   14.48	&     16                  &  $SP_2 = 5$ (single switching action)\\
		\hline   14.59	&     16                  &  $SP_2 = 5$ (single switching action)\\
		\hline   15.79	&     17                  &  $SP_2 = 6$ (single switching action)\\
		\hline   16.55	&     17                  &  $SP_2 = 7$ (single switching action)\\
		\hline   16.86	&     18                  &  $SP_1 = 5$, $SP_2 = 11$ (two switching points)\\
		\hline   16.91	&     18                  &  $SP_1 = 6$, $SP_2 = 12$ (two switching points)\\
		\hline   17.00	&     18                  &  $SP_1 = 8 $, $SP_2 = 14$ (two switching points)\\
		\hline   17.06	&     18                  &  $SP_1 = 10$, $SP_2 = 16$ (two switching points)\\
		\hline   17.18	&     19                  &  $SP_1 = 13$, $SP_2 = 19$ (two switching points)\\
		\hline   17.30	&     19                  &  $SP_1 = 16$, $SP_2 = 22$ (two switching points)\\
		\hline   18.15	&     21                  &  -\\
		\hline   18.19	&     21                  &  -\\
		\hline   18.32	&     21                  &  -\\
		\hline   18.64	&     22                  &  -\\
		\hline   18.86	&     23                  &  -\\
		\hline   19.15	&     24                  &  -\\
		\hline   19.66	&     25                  &  -\\
		\hline   20.47	&     -	  &  -\\
		\hline   ...	&     -	  &  -\\
		%\hline   20.69	&     -	  &  -\\
		%\hline   21.14	&     -	  &  -\\
		%\hline   21.54	&     -	  &  -\\
		\hline \hline
	\end{longtable}

%% file: tables/voi_con_s1.tex
	\newcommand{\compExtraB}{
    \makecell[l]{2H-L-H-L-H-L-H-10L or \\
                 2H-L-H-L-H-3L-H-8L or \\
                 2H-L-H-L-H-5L-H-6L or \\
                 2H-L-H-L-H-7L-H-4L or \\
                 2H-L-H-L-H-9L-H-2L or \\
                 2H-L-H-3L-H-L-H-8L or \\
                 2H-L-H-3L-H-3L-H-6L or \\
                 2H-L-H-3L-H-5L-H-4L}}

	\begin{longtable} {p{0.1\columnwidth}|p{0.2\columnwidth}|p{0.26\columnwidth}|p{0.23\columnwidth}}
		\hline $S$     &    $\Omega$            &    \textbf{aware firm's \newline advertising policy}   &  \textbf{naive firm's  \newline advertising policy}\\
		%	\hline 14       & \{A, B, D, E, F\} &  -    										&    -\\
		%	\hline 20       & \{A, B, D, E, F\} &  -    										&    -\\
		%	\hline 26       & \{A, B, D, E, F\} &  -    										&    -\\
		\hline ...    & \{A, B, D, E, F\} &  -    										&    -\\
			\hline 35       & \{A, B, D, E, F\} &  -    										&    -\\
			\hline 59       & \{A, B, D, E, F\} &  4H-L-H-2L-5H-2L-H-2L		            &    -\\
			\hline 81       & \{A, B, D, E, F\} &  2H-L-H-L-H-2L-6H-2L-H		        &    2L H alternating %-2L-H-2L-H-2L-H-2L-H-2L 
			\\
			\hline 88       & \{A, B, D, E, F\} &  2H-2L-H-L-2H-4L-2H-2L-H		        &    2L-H L alternating %-H-L-H-L-H-L-H-L-H-L-H-L-H 
			\\
			\hline 93       & \{A, B, D, E, F\} &  H-L-H-4L-H-L-6H-3L		            &    2L-H L alternating %-H-L-H-L-H-L-H-L-H-L-H-L-H 
			\\
			\hline 100     	& \{A, B, D, E, F\} &  H-L-H-L-6H-4L-H-3L		            &    2L-H L alternating %-H-L-H-L-H-L-H-L-H-L-H-L-H 
			\\
			\hline 142     	& \{A, B, C, E, F\} &  2H-L H alternating %-L-H-L-H-L-H-L-H-L-H-L-H-L		
			&    2L-H L alternating %-H-L-H-L-H-L-H-L-H-L-H-L-H
			\\
			\hline 228     	& \{A, B, C, E, F\} &  \compExtraB                                              &	2L-H-L-H-3L-H-5L-H-3L\\
			\hline 307     	& \{A, B, C, E, F\} &  2H-16L	    &   2L-H-15L  \\
		%	\hline 314     	& \{A, B, C, E, F\} &  2H-16L	    &   2L-H-15L  \\
		%	\hline 316     	& \{A, B, C, E, F\} &  2H-16L	    &   2L-H-15L  \\
		%	\hline 318     	& \{A, B, C, E, F\} &  2H-16L	    &   2L-H-15L  \\
	%		\hline 324     	& \{A, B, C, E, F\} &  2H-16L	    &   2L-H-15L  \\
		%	\hline 335     	& \{A, B, C, E, F\} &  2H-16L	    &   2L-H-15L  \\
			\hline ...    	& \{A, B, C, E, F\} &  2H-16L	    &   2L-H-15L  \\
			\hline 336     	& \{A, B, C, E, F\} &  2H-16L	    &   2L-H-15L  \\
			\hline 380     	& \{C, E, F	\} &       2H-16L	    &   2L-H-15L  \\
				\hline ...     	& \{C, E, F	\} &       2H-16L	    &   2L-H-15L  \\
		%	\hline 400     	& \{C, E, F	\} &       2H-16L	    &   2L-H-15L  \\
		%	\hline 420     	& \{C, E, F	\} &       2H-16L	    &   2L-H-15L  \\
		%	\hline 426     	& \{C, E, F	\} &       2H-16L	    &   2L-H-15L  \\
		%	\hline 449     	& \{C, E, F	\} &       2H-16L	    &   2L-H-15L  \\
		%	\hline 454     	& \{C, E, F	\} &       2H-16L	    &   2L-H-15L  \\
		%	\hline 472     	& \{C, E, F	\} &       2H-16L	    &   2L-H-15L  \\
		%	\hline 474     	& \{C, E, F	\} &       2H-16L	    &   2L-H-15L  \\
		%	\hline 480     	& \{C, E, F	\} &       2H-16L	    &   2L-H-15L  \\
		%	\hline 488     	& \{C, E, F	\} &       2H-16L	    &   2L-H-15L  \\
		%	\hline 494     	& \{C, E, F	\} &       2H-16L	    &   2L-H-15L  \\
		%	\hline 495     	& \{C, E, F	\} &       2H-16L	    &   2L-H-15L  \\
		%	\hline 509     	& \{C, E, F	\} &       2H-16L	    &   2L-H-15L  \\
		%	\hline 541     	& \{C, E, F	\} &       2H-16L	    &   2L-H-15L  \\
		%	\hline 543     	& \{C, E, F	\} &       2H-16L	    &   2L-H-15L  \\
		%	\hline 545     	& \{C, E, F	\} &       2H-16L	    &   2L-H-15L  \\
		%	\hline 556     	& \{C, E, F	\} &       2H-16L	    &   2L-H-15L  \\
		%	\hline 585     	& \{C, E, F	\} &       2H-16L	    &   2L-H-15L  \\
		%	\hline 588     	& \{C, E, F	\} &       2H-16L	    &   2L-H-15L  \\
		%	\hline 600     	& \{C, E, F	\} &       2H-16L	    &   2L-H-15L  \\
		%	\hline 627     	& \{C, E, F	\} &       2H-16L	    &   2L-H-15L  \\
		%	\hline 639     	& \{C, E, F	\} &       2H-16L	    &   2L-H-15L  \\
		%	\hline 652     	& \{C, E, F	\} &       2H-16L	    &   2L-H-15L  \\
		\hline\hline
	\end{longtable}

%% file: tables/voi_con_s2.tex
	\begin{longtable} {p{0.08\columnwidth}|p{0.2\columnwidth}|p{0.23\columnwidth}|p{0.22\columnwidth}}
		\hline\hline $S$     &    $\Omega$            &    \textbf{aware firm's \newline advertising policy}   &  \textbf{naive firm's \newline advertising policy}\\
		\hline 511	    & \{A, B, C, E, F\}	&	3H-8L	&	5L-2H-4L\\
		\hline 647	    & \{A, B, C, E, F\}	&	4H-7L	&	5L-3H-3L\\
		\hline 797	    & \{A, B, C, E, F\}	&	4H-7L	&	5L-4H-2L\\
		\hline 868	    & \{A, B, C, E, F\}	&	5H-6L	&	5L-4H-2L\\
		\hline 976	    & \{A, B, C, E, F\}	&	6H-5L	&	5L-5H   \\
		\hline 1,059	& \{A, B, C, E, F\}	&	7H-4L	&	5L-6H\\
		\hline 1,072	& \{A, B, C, E, F\}	&	7H-4L	&	5L-6H\\
		\hline 1,121	& \{A, B, C, E, F\}	&	9H-2L	&	5L-6H\\
		\hline 1,525	& \{C, E, F\}		&	11H	&	-               \\
		\hline ...	& \{C, E, F\}		&	11H	&	-               \\
	%	\hline 1,548	& \{C, E, F\}		&	11H	&	-               \\
	%	\hline 1,570	& \{C, E, F\}		&	11H	&	-               \\
	%	\hline 1,651	& \{C, E, F\}		&	11H	&	-               \\
	%	\hline 1,676	& \{C, E, F\}		&	11H	&	-               \\
	%	\hline 1,736	& \{C, E, F\}		&	11H	&	-               \\
	%	\hline 1,817	& \{C, E, F\}		&	11H	&	-               \\
	%	\hline 1,896	& \{C, E, F\}		&	11H	&	-               \\
	%	\hline 1,922	& \{C, E, F\}		&	11H	&	-               \\
	%	\hline 1,928	& \{C, E, F\}		&	11H	&	-               \\
	%	\hline 2,004	& \{C, E, F\}		&	11H	&	-               \\
	%	\hline 2,107	& \{C, E, F\}		&	11H	&	-               \\
		\hline\hline
	\end{longtable}

%% file: tables/voi_con_s3.tex
	\begin{longtable} {p{0.1\columnwidth}|p{0.2\columnwidth}|p{0.23\columnwidth}|p{0.23\columnwidth}}
		\hline\hline $S$     &    $\Omega$            &    \textbf{aware firm's \newline advertising policy}   &  \textbf{naive firm's \newline advertising policy}\\
		\hline 571  & \{A, B, C, E, F\}	&	3H-16L	& -\\
	%	\hline 596  & \{A, B, C, E, F\}	&	3H-16L	& -\\
	%	\hline 691  & \{A, B, C, E, F\}	&	3H-16L	& -\\
	%	\hline 699  & \{A, B, C, E, F\}	&	3H-16L	& -\\
	%	\hline 766  & \{A, B, C, E, F\}	&	3H-16L	& -\\
	%	\hline 804  & \{A, B, C, E, F\}	&	3H-16L	& -\\
	\hline ...  & \{A, B, C, E, F\}	&	3H-16L	& -\\
		\hline 841  & \{A, B, C, E, F\}	&	3H-16L	& -\\
		\hline 1024 & \{C, E, F\}		&	3H-16L	& -\\
		\hline 1048 & \{C, E, F\}		&	3H-16L	& -\\
		\hline 1118 & \{C, E, F\}		&	3H-16L	& -\\
		\hline 1232 & \{C, E, F\}		&	-	& -                                    \\
		\hline ... & \{C, E, F\}		&	-	& -                                    \\
		%\hline 1,266 & \{C, E, F\}		&	-	& -                                    \\
	%	\hline 1,275 & \{C, E, F\}		&	-	& -                                    \\
	%	\hline 1,307 & \{C, E, F\}		&	-	& -                                    \\
		%\hline 1,341 & \{C, E, F\}		&	-	& -                                    \\
	%	\hline 1,346 & \{C, E, F\}		&	-	& -                                    \\
	%	\hline 1,484 & \{C, E, F\}		&	-	& -                                    \\
	%	\hline 1,571 & \{C, E, F\}		&	-	& -                                    \\
	%	\hline 1,672 & \{C, E, F\}		&	-	& -                                    \\
	%	\hline 1,678 & \{C, E, F\}		&	-	& -                                    \\
		\hline\hline
	\end{longtable}

%% file: tables/voi_con_w.tex
	\begin{longtable} {p{0.1\columnwidth}|p{0.2\columnwidth}|p{0.3\columnwidth}|p{0.3\columnwidth}}
		\hline \hline $w$     &    $\Omega$            &    \textbf{aware firm's \newline advertising policy}   &  \textbf{naive firm's \newline advertising policy}\\
		%\hline 0.053	& \{A, B, C, E, F\}	&   always $L$	&   always $L$	\\
		%\hline 0.086	& \{A, B, C, E, F\}	&   always $L$	&   always $L$	\\
		%\hline 0.091	& \{A, B, C, E, F\}	&   always $L$	&   always $L$	\\
		%\hline 0.118	& \{A, B, C, E, F\}	&   always $L$	&   always $L$	\\
		%\hline 0.148	& \{A, B, C, E, F\}	&   always $L$	&   always $L$	\\
		%\hline 0.185	& \{A, B, C, E, F\}	&   always $L$	&   always $L$	\\
		%\hline 0.192	& \{A, B, C, E, F\}	&   always $L$	&   always $L$	\\
		%\hline 0.194	& \{A, B, C, E, F\}	&   always $L$	&   always $L$	\\
		%\hline 0.231	& \{A, B, C, E, F\}	&   always $L$	&   always $L$	\\
		%\hline 0.243	& \{A, B, C, E, F\}	&   always $L$	&   always $L$	\\
		%\hline 0.254	& \{A, B, C, E, F\}	&   always $L$	&   always $L$	\\
		%\hline 0.272	& \{A, B, C, E, F\}	&   always $L$	&   always $L$	\\
		%\hline 0.275	& \{A, B, C, E, F\}	&   always $L$	&   always $L$	\\
		%\hline 0.276	& \{A, B, C, E, F\}	&   always $L$	&   always $L$	\\
		%\hline 0.293	& \{A, B, C, E, F\}	&   always $L$	&   always $L$	\\
		%\hline 0.311	& \{A, B, C, E, F\}	&   always $L$	&   always $L$	\\
		%\hline 0.319	& \{A, B, C, E, F\}	&   always $L$	&   always $L$	\\
	%	\hline 0.46	& \{A, B, C, E, F\}	&   always $L$	&   always $L$	\\
		\hline 0.46	& \{A, B, C, E, F\}	&   always $L$	&   always $L$	\\
	%	\hline 0.47	& \{A, B, C, E, F\}	&   always $L$	&   always $L$	\\
		\hline 0.47	& \{A, B, C, E, F\}	&   always $L$	&   always $L$	\\
		\hline 0.47	& \{A, B, C, E, F\}	&   always $L$	&   always $L$	\\
	%	\hline 0.48	& \{A, B, C, E, F\}	&   always $L$	&   always $L$	\\
		\hline 0.53	& \{A, B, C, E, F\}	&    H-23L          	                                                    &   6L-H L alternating %-H-L-H-L-H-L-H-L-H-L-H-L-H-L-H-L	
		\\
		\hline 0.53	& \{A, B, C, E, F\}	&    H-23L          	                                                    &   6L-H L alternating %-H-L-H-L-H-L-H-L-H-L-H-L-H-L-H-L
		\\
		\hline 0.54	& \{A, B, C, E, F\}	&   2L-H-21L                                                           	    &   6L-H L alternating %-H-L-H-L-H-L-H-L-H-L-H-L-H-L-H-L	
		\\
		\hline 0.56	& \{A, B, C, E, F\}	&    H-4L-2H-8L-H-8L                                                      	&   4L-H L alternating %-H-L-H-L-H-L-H-L-H-L-H-L-H-L-H-L-H-L
		\\
		\hline 0.57	& \{A, B, C, E, F\}	&    H-6L-H-5L-H-3L-H-6L                                             	    &   4L-H L alternating %-H-L-H-L-H-L-H-L-H-L-H-L-H-L-H-L-H-L
		\\
		\hline 0.60	& \{A, B, C, E, F\}	&    H-4L-4H-10L-H-4L          	                                            &   4L-H L alternating %-H-L-H-L-H-L-H-L-H-L-H-L-H-L-H-L-H-L	
		\\
	%	\hline 0.60	& \{A, B, C, E, F\}	&   H-2L-2H-L-H-L-2H-2L-H-3L-H-7L          	                                &   4L-H L alternating %-H-L-H-L-H-L-H-L-H-L-H-L-H-L-H-L-H-L\\
		\hline 0.61	& \{A, B, C, E, F\}	&   H-2L-2H-L H alternating % -L-H-L-H-L-H-L-H-L-H-L-H-L-H-L-H-L  
		&   - 														\\
		\hline 0.73	& \{A, B, C, E, F\}	&   -			&   - 		\\
		\hline
		...	& \{A, B, C, E, F\}	&   -														&   - 														 \\
	
	%	\hline 0.74	& \{A, B, C, E, F\}	&   -														&   - 														 \\
	%	\hline 0.75	& \{A, B, C, E, F\}	&   -														&   - 														 \\
%		\hline 0.76	& \{A, B, C, E, F\}	&   -														&   - 														 \\
	%	\hline 0.77	& \{A, B, C, E, F\}	&   -														&   - 														 \\
%		\hline 0.94	& \{A, B, C, E, F\}	&   -														&   - 														 \\
	%	\hline 0.94	    & \{A, B, C, E, F\}	&   -														&   - 														 \\
%		\hline 0.95	& \{A, B, C, E, F\}	&   -														&   - 														 \\
	%	\hline 0.95	& \{A, B, C, E, F\}	&   -														&   - 														 \\
		\hline\hline
	\end{longtable}

%% file: tables/voi_con_alpha.tex
	\begin{longtable} {p{0.1\columnwidth}|p{0.18\columnwidth}|p{0.23\columnwidth}|p{0.23\columnwidth}}
		\hline\hline $\alpha$     &    $\Omega$            &    \textbf{aware firm's \newline advertising policy}   &  \textbf{naive firm's \newline advertising policy}\\
	\hline 0.22  &	\{A, B\}			&   H-10L		& 	H-10L \\
	\hline 0.60  &	\{A, B, C\}			&   H-10L 		&	H-10L \\
	\hline 0.87  &	\{A, B, C\}			&   H-10L 		&	H-10L \\
	\hline 1.16  &	\{A, B, C, E\}		&   H-10L 		&	H-10L \\
	\hline 1.35 &	\{A, B, C, E\}		&   2H-9L		&	2H-9L \\
	\hline 1.78  &	\{A, B, C, E\}		&   2H-9L		&	2H-9L \\
	\hline 1.82  &	\{A, B, C, E\}		&   2H-9L		&	2H-9L \\
	\hline 2.24  &	\{A, B, C, E\}		&   2H-9L		&	2H-9L \\
	\hline 2.57  &	\{A, B, C, E, F\}	&	3H-8L		&	L-2H-8L \\
	\hline 2.92  &	\{A, B, C, E, F\}	&	3H-8L		&	L-2H-8L \\
	\hline 3.30  &	\{A, B, C, E, F\}	&	3H-8L		&	2L-2H-7L \\
	\hline 3.77  &	\{A, B, C, E, F\}	&	3H-8L		&	3L-2H-6L \\
	\hline 4.11 &	\{A, B, C, E, F\}	&	3H-8L		&	- \\
	\hline 6.14 &	\{A, B, C, E, F\}	&	6H-5L		&	- \\
	\hline 6.19  &	\{A, B, C, E, F\}	&	6H-5L		&	- \\
	\hline 7.63  &	\{A, B, C, E, F\}	&	7H-4L		&	- \\
	\hline 7.80  &	\{A, B, C, E, F\}	&	7H-4L		&	- \\
	\hline 7.97  &	\{A, B, C, E, F\}	&	7H-4L		&	- \\
	\hline 8.27  &	\{A, B, C, E, F\}	&	-            		&	- \\
	\hline 8.60  &	\{A, B, C, E, F\}	&	-            		&	- \\
		\hline\hline
	\end{longtable}